\documentclass[onefignum,onetabnum]{siamart190516}

\usepackage{amsfonts}
\usepackage{amsmath}
\usepackage{epstopdf}
\usepackage{graphicx}
\usepackage{subfigure}
\usepackage{subfigmat}

\newcommand{\R}{\mathbb{R}}

\newcommand{\bmat}[1]{\begin{bmatrix}#1\end{bmatrix}}

\let\mathopfont=\mathrm

\newcommand{\diag}{\mathop{\mathopfont{diag}}}

\newcommand{\mcl}[1]{\mathcal{ #1}}
\newcommand{\mbf}[1]{\mathbf{ #1}}

\newlength{\eqnspace}
\newlength{\parspace}
\newtheorem{defn}{Definition}

\usepackage{lipsum}
\usepackage{algorithmic}
\ifpdf
  \DeclareGraphicsExtensions{.eps,.pdf,.png,.jpg}
\else
  \DeclareGraphicsExtensions{.eps}
\fi

\newsiamremark{remark}{Remark}
\newsiamremark{hypothesis}{Hypothesis}
\crefname{hypothesis}{Hypothesis}{Hypotheses}
\newsiamthm{claim}{Claim}

\headers{PIE-Galerkin Framework for PDE Solution}{Y. T. Peet, and M. M. Peet}

\title{A New Treatment of Boundary Conditions in PDE Solution with Galerkin Methods via Partial Integral Equation Framework\thanks{Submitted to the editors DATE.
\funding{This work was supported by grants NSF CMMI-1935453 and NSF CAREER-1944568.}}}

\author{Yulia T. Peet\thanks{School for Engineering of Matter, Transport and Energy, Arizona State University, Tempe, AZ, 85287
  (\email{ypeet@asu.edu}, \email{mpeet@asu.edu}, \url{http://isim.asu.edu}, \url{http://control.asu.edu}).}
\and Matthew M. Peet\footnotemark[2]}

\usepackage{amsopn}

\ifpdf
\hypersetup{
  pdftitle={A New Treatment of Boundary Conditions in PDE Solution with Galerkin Methods via Partial Integral Equation Framework},
  pdfauthor={Y. T. Peet, and M. M. Peet}
}
\fi

\begin{document}

\maketitle

\begin{abstract}
  We present a new analytical and numerical framework for solution of Partial Differential Equations (PDEs) that is based on an analytical transformation that moves the boundary constraints into the dynamics of the corresponding governing equation. The framework is based on a Partial Integral Equation (PIE) representation of PDEs, where a PDE equation is transformed into an equivalent PIE representation that does not require boundary conditions on its solution state. The PDE-PIE framework allows for a development of a generalized PIE-Galerkin approximation methodology for a broad class of linear PDEs with non-constant coefficients governed by non-periodic boundary conditions, including, e.g., Dirichlet, Neumann and Robin boundaries. The significance of this result is that solution to almost any linear PDE can now be constructed in a  form of an analytical approximation based on a series expansion using a suitable set of basis functions, such as, e.g., Chebyshev polynomials of the first kind, irrespective of the boundary conditions. In many cases involving homogeneous or simple time-dependent boundary inputs, an analytical integration in time is also possible. We present several PDE solution examples in one spatial variable implemented with the developed PIE-Galerkin methodology using both analytical and numerical integration in time. The developed framework can be naturally extended to multiple spatial dimensions and, potentially, to nonlinear problems.
\end{abstract}


\textit{Science is a Differential Equation. Religion is a Boundary Condition}. -- Alan Turing (1912--1954).
\section{Introduction}
 The need to enforce boundary conditions has been a major challenge in developing analytical and numerical tools for finding solutions of Partial Differential Equations (PDEs) ever since the concept of PDEs emerged in the 18th century following the works of  Euler, d'Alembert, Lagrange and Laplace, who recognized their central role in the description of the laws of nature~\cite{brezis1998partial}. To enforce boundary conditions, a solution is typically split into a homogeneous part that satisfies homogeneous boundary conditions, and an inhomogeneous part~\cite{deville2002high, grigoryan2010partial}. For the inhomogeneous part, one must typically find a general appropriately smooth function defined on the solution domain that satisfies specified constraints on the boundary, a task that is daunting by itself.  However, it is the search for a homogeneous solution, which is required to satisfy \textit{both} the PDE and the homogeneous boundary condition, that represents the utmost challenge and has hindered a development of a unifying theoretical framework for solving PDE equations for more than two centuries.

The easiest way of handling boundary conditions would be to seek a solution to a PDE in terms of the functions that already satisfy the boundary conditions, which is done in the so-called Galerkin methods~\cite{canuto1988spectral}. Unfortunately, such basis functions are readily available only for a limited class of problems, e.g., the ones with periodic boundary conditions, for which Fourier methods based on harmonic function expansions offer an elegant, efficient, and generalizable approach to the solution of PDEs with periodic boundaries~\cite{gottlieb1977numerical}. 
For boundary conditions other than periodic, the picture is more obscure. An unfortunate fact to accept is that there are no convenient basis functions (viz. harmonic functions or classical orthogonal polynomials) that satisfy general, non-periodic boundary conditions. This yields, in a classical PDE analysis framework, three options: 1) construct more sophisticated basis functions from the primary ones that do satisfy boundary conditions~\cite{shen1994efficient}--\cite{shen2003new}, 2) enforce boundary conditions discretely on the expansion coefficients~\cite{haidvogel1979accurate, canuto1986boundary, siyyam1997accurate}, 3) enforce boundary conditions in a weak form, by introducing penalty terms or Lagrange multipliers into the variational form of the equations~\cite{nitsche1971variationsprinzip, bazilevs2007weakimp, jovanovic2016ritz}. The problem with the first approach is that it leads to a complicated basis that depends on the order of equations and on the boundary conditions~\cite{shen1994efficient, shen2003new, guo2009generalized, yu2019jacobi}, limiting the generalizability of approach. The second option, which is typically used in conjunction with either tau methods~\cite{haidvogel1979accurate, siyyam1997accurate} or nodal/collocation methods~\cite{canuto1986boundary, fischer1997overlapping, karniadakis2005spectral}, is inherently tied to a discretization, and thus has limited options for providing generalized close-form solutions that are useful for analysis and control of continuous models~\cite{smyshlyaev2005backstepping, fridman2009lmi,peet2020partial}. Additionally, it requires an ad-hoc modification of the discrete matrix operators,  which can lead to ill-conditioned matrices and effect stability and accuracy of the methods~\cite{gottlieb1977numerical,  lehotzky2016pseudospectral, bostrom2017boundary}. The weak enforcement of the boundary conditions attempts to circumvent the above deficiencies~\cite{ruess2013weakly, vymazal2019weak}. However, it introduces a tunable penalty parameter, which is not known from the first principles, problem-dependent, and leads to a lack of robustness of the solution~\cite{vymazal2019weak, freund1995weakly, juntunen2009nitsche}. Moreover, a weak imposition of boundary conditions forfeits the possibility of exactly satisfying the conservation laws, which, in some cases, e.g. for hyperbolic systems, is highly desirable~\cite{leveque1992numerical, tadmor2012review, bansal2020structure}. 

In this paper, we present a conceptually new approach to address the problems associated with the enforcement of boundary conditions in the solution of PDEs. Specifically, we exploit a novel Partial Integral Equation (PIE) framework for representation of Partial Differential Equations~\cite{peet2020partial}.  In this framework, PIEs can be used to equivalently represent the solution of PDEs, yet require no boundary conditions. This is due to the fact that solutions of the PIE equations are expressed using a so-called ``fundamental state'', which consists of specially constructed functions that include derivatives of the primary solution. In a PIE, the fundamental state solution function lies in a space of $L_2$ square-integrable functions and requires no boundary conditions. Instead, the effect of boundary conditions is incorporated directly into the PIE dynamics through the construction of the corresponding partial-integral operators. This integral representation essentially acts to move the boundary conditions from the realm of ``religion'' (artificial constraints on a solution) to the realm of ``science'' (integro-differential equations). Significantly, by solving PIEs, we are now free to represent the solution using any choice of approximation space without the need to impose the boundary conditions on the basis functions for that space! This means that we can now use Galerkin method based on a native set of orthogonal polynomials~\cite{gottlieb1977numerical, canuto1988spectral} for a large class of PDEs with non-periodic boundary conditions, extending the benefits of classical Galerkin methods to a broad range of PDE systems. In this paper, the corresponding PIE-Galerkin formulation is derived and implemented for linear PDEs with non-constant coefficients in one spatial variable, governed by a general set of boundary constraints that can include, e.g., Dirichlet, Neumann and Robin boundary conditions.

Since the idea of solving boundary value problems by relating the boundary condition functions to the interior solution resonates with several other techniques in mathematics, here we contrast our approach with the popular methods of Green functions~\cite{stakgold1979green, roach1982green, greenlibrary} and boundary integral equations (BIE)~\cite{atkinson1997numerical, marin2012highly, carvalho2020asymptotic}. Both Green functions and BIE approaches require a knowledge of the fundamental solutions of the corresponding differential operator, while no such a-priori knowledge is required in the current approach. Note that the ``fundamental state'' in a PIE is completely different from the ``fundamental solution'', which is a response of a linear differential operator to an impulse forcing~\cite{encyclopaedia, kythe2012fundamental}. In a classical Green function approach, these functions also require to satisfy homogeneous boundary conditions.
 In a BIE formulation, this requirement is relaxed, and solution satisfying the desired boundary conditions is formulated as a continuous superposition of arbitrary fundamental solutions, giving rise to an integral equation for the distribution density on the boundary of the domain~\cite{atkinson1997numerical, marin2012highly, carvalho2020asymptotic}. Both these approaches are fundamentally different from the  methodology presented in this paper, since, first of all, the integral operators act on the domain boundary in BIEs, and they do on the domain interior in PIEs, and, second, the PIE formulation does not require any a-priori knowledge of the fundamental solutions, which are only available for certain equations~\cite{stakgold1979green, roach1982green, greenlibrary}, and, for the case of non-constant coefficients, only approximately~\cite{johnson2010notes, manen2005modeling, sesma2001approximate}.
 
Several other approaches utilized a spatial integration of PDEs to eliminate function derivatives from a solution as a means to arrive at better-conditioned and more compact discrete matrix operators attributed to an integration as opposed to a differentiation procedure~\cite{greengard1991spectral, hiegemann1997chebyshev, driscoll2010automatic}. However, these approaches do not eliminate the boundary conditions and still have to enforce them on a solution, which is typically done at a discrete level by modifying the corresponding rows of discrete matrix operators to represent the algebraic constraints on the expansion coefficients~\cite{greengard1991spectral, hiegemann1997chebyshev}, similar to the corresponding differentiation tau or collocation techniques. 
 
 In this regard, it is also useful to mention the Fokas method~\cite{fokas1997unified}, which seeks to propose a unified transform procedure for solving initial-boundary value problems. The method involves performing joint Fourier-type integral transforms of the PDE together with initial and boundary conditions in space and time, solving for a global relation, and performing an inverse Fourier transform, which involves taking an indefinite integral over specified contours in a complex half-plane. This approach, however, is associated with certain difficulties as applied to a general case: first, it relies on the existence of a Lax pair~\cite{fokas1998lax}, which can only be formulated for certain equations~\cite{fokas1997unified, treharne2007initial}; second, extension to a finite-size interval is challenging in that it yields an integrand which is no longer analytic, and requires evaluation of the residues at the complex poles, which may lack convergence and complicate the computation~\cite{deconinck2014method, kesici2018numerical}. As opposed to the Fokas method which predominantly seeks to provide an integral solution to an initial-boundary value problem (IBVP), the current PIE framework reformulates the PDEs into an equivalent set of governing equations, which is suitable not only for the IBVP solution, but also for analysis and control of PDEs~\cite{das2019h, peet2020partial}, as well as for coupling PDEs with auxiliary models, such as ODEs~\cite{shivakumar2021extension, shivakumar2020duality}, or other PDEs interfacing through a joint boundary.  
 
 The current paper is organized as follows. In \cref{sec:pie}, we present a general formulation of the PIE framework for linear PDEs with non-constant coefficients and extend the original representation of~\cite{peet2020partial} to include inhomogeneous boundary conditions. In \cref{sec:numsol}, we introduce Galerkin approach based on Chebyshev polynomials of the first kind for a solution of the PDE equations in the PIE framework, and present the corresponding stability and convergence proofs for the PIE-Galerkin approach.  In \cref{sec:example}, we show numerical examples, followed by conclusions in \cref{sec:conclude}.

\section{Partial Integral Equations Framework}\label{sec:pie}

\subsection{Standardized PDE Representation}\label{sec:standard}
We first define some notations. We solve a Partial Differential Equation (PDE), or a coupled system of PDEs, on a spatio-temporal domain $(x,t)\in([a,b]\times\mathbb R^+)$. Let $L_2[a,b]^n$ be a space of $\mathbb{R}^n$-valued square-integrable functions in a  Lebesgue sense defined on $[a,b]$, with a suitable inner product. We adopt the notation $H_{k}[a,b]^n$ to denote a Sobolev subspace of $L_2[a,b]^n$ defined as $\big\{\mathbf{u}\in L_2[a,b]^n: \frac{\partial^q \mathbf{u}}{\partial x^q}\in L_2[a,b]^n,\: \forall \:q\leq k\big\}$. $I_n\in\mathbb R^{n\times n}$ is used to denote the identity matrix, while $0_n$ denotes a zero vector of size $n$. It is implied that, for all the solution states $u(x,t)$, a partial derivative with respect to time exists for $t\in \mathbb R^+$.


We now consider a class of linear Partial Differential Equations in one spatial dimension given in its ``state-space'' representation~\cite{peet2018new, peet2020partial}

\begin{align}\label{eq:primary}
\bmat{\mathbf{u}_0(x,t)\\\mathbf{u}_1(x,t)\\\mathbf{u}_2(x,t)}_t&= A_0(x)\bmat{\mathbf{u}_0(x,t)\\\mathbf{u}_1(x,t)\\\mathbf{u}_2(x,t)}+A_1(x)\bmat{\mathbf{u}_1(x,t)\\\mathbf{u}_2(x,t)}_x +A_2(x) \bmat{\mathbf{u}_2(x,t)}_{xx}+\mbf f(x,t),
\end{align}
boundary conditions,
\begin{equation}\label{eq:bcon}
B {\scriptsize \bmat{ \mathbf{u}_1(a,t) \\ \mathbf{u}_1(b,t) \\ \mathbf{u}_2(a,t) \\ \mathbf{u}_2(b,t)\\ \mathbf{u}_{2x}(a,t) \\ \mathbf{u}_{2x}(b,t)}}=\mathbf{h}(t)\in C^1(\mathbb R^+)^{n_1+2n_2},
\end{equation}
and initial conditions
\begin{equation}\label{eq:incon}
\bmat{\mathbf{u}_0(x,0)\\\mathbf{u}_{1}(x,0)\\\mathbf{u}_{2}(x,0)}=\pmb \beta^h(x).
\end{equation}

Here, $A_0(x): R\rightarrow R^{\,ns\times ns} , A_1(x): R\rightarrow R^{\,ns\times (n_1+n_2)}, A_2(x): R\rightarrow R^{\,ns\times n_2}$ are bounded matrix-valued real functions. We introduce a functional space $X$ of dimension $ns=n_0+n_1+n_2$, such that
\begin{equation}
X:=\left\{\bmat{\mathbf{u}_0(x,t)\\\mathbf{u}_{1}(x,t)\\\mathbf{u}_{2}(x,t)}\in \bmat{L_2[a,b]^{n_0}\\H_1[a,b]^{n_1}\\H_2[a,b]^{n_2}},\,t\in \mathbb R^{+}\right\}.
\end{equation}
Furthermore, we denote a subset of functions $X^h\subset X$ satisfying the boundary conditions (\ref{eq:bcon}) as 
\begin{equation}
X^{h}:=\left\{\bmat{\mathbf{u}_0(x,t)\\\mathbf{u}_{1}(x,t)\\\mathbf{u}_{2}(x,t)}\in X \cap B {\scriptsize \bmat{ \mathbf{u}_1(a,t) \\ \mathbf{u}_1(b,t) \\ \mathbf{u}_2(a,t) \\ \mathbf{u}_2(b,t)\\ \mathbf{u}_{2x}(a,t) \\ \mathbf{u}_{2x}(b,t)}}=\mathbf{h}(t),\,t\in \mathbb R^+\right\}.
\end{equation} 

We say that a solution
\begin{equation}\label{eq:pstate} 
\mathbf{u}^{h} (x,t)=\bmat{\mathbf{u}^h_0(x,t)\\\mathbf{u}^h_1(x,t)\\\mathbf{u}^h_2(x,t)}\in X^h,
\end{equation}
to the equation (\ref{eq:primary}) with boundary (\ref{eq:bcon}) and initial (\ref{eq:incon}) conditions is in its \textit{primary} state. Here, a superscript $h$ denotes a dependency of the solution on the boundary conditions. Note that, for well-posedness, we demand that initial conditions (\ref{eq:incon}) satisfy boundary conditions at $t=0$, i.e. $\pmb \beta^h(x)\in X^h, t=0$. 

To arrive at an equation (\ref{eq:primary}), a set containing an original scalar-valued dependent variable $v(x,t)$ of a PDE (or a vector-valued dependent variable $\mathbf{v}(x,t)$ for a system of coupled PDEs) and their partial derivatives must be transformed into its corresponding state-space form, where the functions $\mathbf{u}_0(x,t)\in L_2[a,b]^{n_0}$ admit no partial spatial derivatives, the functions $\mathbf{u}_1(x,t)\in H_1[a,b]^{n_1}$ admit only first-order partial spatial derivatives, and the functions $\mathbf{u}_2(x,t)\in H_2[a,b]^{n_2}$ admit up to second-order spatial partial derivatives. Note that the functions $\{\mathbf{u}_0, \mathbf{u}_1, \mathbf{u}_2\}$ in a state-space form are generally vector-valued, even if the original dependent variable $v(x,t)$ was a scalar~\cite{peet2018new, peet2020partial}. Matrix $B\in\mathbb{R}^{(n_1+2n_2)\times(2n_1+4n_2)}$ is the boundary conditions matrix, and $\mathbf{h}(t)\in\mathbb{R}^{n_1+2n_2}$ is the vector of the boundary condition values. According to a decomposition of the functions into its state-space form, the functions $\mathbf{u}_0(x,t)$ admit no boundary conditions, functions $\mathbf{u}_1(x,t)$ admit one boundary condition per each scalar component, and functions $\mathbf{u}_2(x,t)$ admit two boundary conditions per each scalar component. Since these boundary conditions can be prescribed either on the left or the right end of the domain, or, in general, contain boundary constraints that couple the two ends, a boundary conditions matrix $B$  has $2n_1+4n_2$ number of columns. Most 1D PDEs can be formulated using this standardized representation, with multiple examples on how to accomplish this transformation for various linear PDE models given in our previous work~\cite{peet2018new, peet2020partial}, and in the numerical examples below.

\subsection{Conversion to a Partial Integral Equation (PIE) Representation}
\subsubsection{Some Useful Preliminaries}

Peet~\cite{peet2020partial} have introduced a framework for converting PDE equations in the form of (\ref{eq:primary}) to a Partial Integral Equation (PIE) form. The original formulation is, however, restricted to a homogeneous case, i.e. a zero forcing function $\mbf f(x,t)$, and homogeneous boundary conditions (\ref{eq:bcon}) given by $\mbf h(t)=0$. Here, we extend the previous result to inhomogeneous boundary conditions in (\ref{eq:bcon}) defined by an arbitrary vector $\mathbf{h}(t) \in C^{1}(\mathbb{R}^+)^{2n_1+4n_2}$, and an arbitary forcing function $\mbf f(x,t)\in L_2^{n_0+n_1+n_2}$ in the equation (\ref{eq:primary}). We will try to minimize the repetition of the proofs that already appeared in~\cite{peet2020partial, peet_arxiv_PDE}, and will refer the reader to these two manuscripts, whenever possible. 

For the homogeneous boundary conditions, we have the following lemma.


\begin{lemma}\label{lemma:space}
If $\mathbf{h}(t)=0$, i.e. boundary conditions are homogeneous, $X^{0}$ is a linear subspace of $X$.
\begin{proof}
We show the following properties of $X^{0}$ that makes it a linear subspace:
\begin{enumerate}
\item The zero element $0_{ns}\in X^{0}$, since $0_{ns}\in X$, and it satisfies (\ref{eq:bcon}) with $\mbf h(t)=0$. 
\item $X^{0}$ is closed under addition and scalar multiplication, since $X$ is closed under addition and scalar multiplication, and these operations preserve homogeneous boundary conditions.  
\end{enumerate}
\end{proof}
\end{lemma}
Note that, for inhomogeneous boundary conditions,  $\mbf h(t)\ne 0$, $X^{h}$ is not a linear subspace, since, for one, it does not contain a zero vector. Instead, it corresponds to an affine space isomorphic to $X^{0}$ that is obtained from $X^{0}$ by a translation transformation, as will be dicussed later.

Given a primary state defined by (\ref{eq:pstate}), we now introduce a \textit{fundamental state} as
\begin{equation}\label{eq:fvector}
\mathbf{u}_f(x,t)=\bmat{\mathbf{u}_{f0}(x,t)\\\mathbf{u}_{f1}(x,t)\\\mathbf{u}_{f2}(x,t)}=\bmat{\mathbf{u}_0(x,t)\\\mathbf{u}_{1x}(x,t)\\\mathbf{u}_{2xx}(x,t)}\in \bmat{(L_2[a,b])^{n_0}\\(L_2[a,b])^{n_1}\\(L_2[a,b])^{n2}},\,t\in\mathbb{R}^+.
\end{equation}
Note that the fundamental state solution is in $L_2[a,b]^{n_0+n_1+n_2}$ space, and thus, it does not admit boundary constraints, which is reflected in the fact that the superscript $h$ is now omitted from the notation. It can be seen, that the fundamental state is related to the primary state by the following differentiation operation
\begin{equation}\label{eq:fstate}
\mbf u_f(x,t)=\mathcal{D}\, \mbf u^{ h}(x,t),
\end{equation}
where the differentiation operator $\mathcal{D}$ has the form
\begin{equation}\label{eq:difopp}
\mathcal{D}:=\bmat{I_{n_0}&&\\&I_{n_1}\,\partial_x&\\&& I_{n_2}\, \partial_x^2}.
\end{equation}
Note that, in general, a map $\mathcal{D}:X\rightarrow L_2^{ns}$ is non-injective, since there can be multiple elements of $X$ mapped into the same fundamental state $\mbf u_f(x,t)$, differing by boundary conditions.

We now proceed with invoking the following lemma proven in~\cite{peet2020partial}.
\begin{lemma}\label{lem:identity1}Suppose that $u\in H_2[a,b]$. Then for any $x \in [a,b]$, 
\begin{align}
u(x)&=u(a)+\int_a^x u_{x}(s)d s  \\
u_x(x)&=u_x(a)+\int_a^x u_{xx}(s)d s \label{tran:ux}\\
u(x)&=u(a)+u_x(a)(x-a)+\int_a^x (x-s)  u_{xx}(s)d s \label{tran:uxx}
\end{align}
\end{lemma}
\begin{proof}
See the manuscript~\cite{peet2020partial} for a proof.
\end{proof}

Next, we define the boundary conditions vectors as
\begin{equation}\label{eq:ubfubc}
 \mbf u_{bf}(t)=\bmat{ \mathbf{u}_1(a,t) \\ \mathbf{u}_1(b,t) \\ \mathbf{u}_2(a,t) \\ \mathbf{u}_2(b,t)\\ \mathbf{u}_{2x}(a,t) \\ \mathbf{u}_{2x}(b,t)},\: \mbf u_{bc}(t)=\bmat{ \mathbf{u}_1(a,t)  \\ \mathbf{u}_2(a,t) \\ \mathbf{u}_{2x}(a,t)},
\end{equation}
where $\mbf u_{bf}(t)$ corresponds to a full set of boundary conditions, and $\mbf u_{bc}(t)$ corresponds to a ``core'' set of boundary conditions~\cite{peet_arxiv_PDE}. Note that, under this definition, boundary constraint (\ref{eq:bcon}) reads as $B\mbf u_{bf}(t)=\mbf h(t)$.

We now have to introduce the notation to define a partial-integral operator of a specific form, which will be referred to as a 3-PI operator. 

\begin{defn}
If $N_0: [a,b]\rightarrow \mathbb R^{n\times n}$, $N_1: [a,b]^2\rightarrow \mathbb R^{n\times n}$, $N_2: [a,b]^2\rightarrow \mathbb R^{n\times n}$ are bounded matrix-valued functions, we define a 3-PI operator $\mcl P:L_2^n[a,b]\rightarrow L_2^n[a,b]$ as
\begin{align}
&(\mcl P \mbf u)(x):=\left(\mcl P_{\{N_0,N_1,N_2\}}\mbf u\right)(x):= N_0(x) \mbf u(x) \\ 
& \qquad +\int_a^x N_1(x,s)\mbf u(s)\,d s+\int_a^bN_2(x,s)\mbf u(s)d s,  \notag
\end{align}
where $N_0$ defines a multiplier operator and $N_1,N_2$ define the kernels of the integral operators.
\end{defn}

Our definition is slightly different from the one presented in~\cite{peet2020partial} in that a last term here is defined as an integration from $a$ to $b$, while it is defined as an  integration from $x$ to $b$ in \cite{peet2020partial}, however, with the appropriate modification of the integral kernels, the two definitions are equivalent. It is proven in~\cite{peet2020partial} that 3-PI operators are closed under addition, scalar multiplication and composition, and thus form an algebra. For a reference, a composition rule for 3-PI operators with the current definition is included in the~\cref{sec:composition}.

We now define two specific 3-PI operators, which will be instrumental for conversion of the PDEs into the PIE framework, as will be seen below. 
\begin{align}\label{eq:3-pi}
\mcl T&:=\mcl{P}_{\{G_0,G_1,G_2\}},\qquad \mcl A:=\mcl P_{\{H_0,H_1,H_2\}},\notag\\
H_0(x)&=A_0(x)G_0+A_1(x)G_3+A_{20}(x),\notag\\
H_1(x,s)&=A_0(x)G_1(x,s)+A_1(x)G_4(s),\\
H_2(x,s)&=A_0(x)G_2(x,s)+A_1(x)G_5(s),\notag\\
A_{20}(x)&=\bmat{0&0&A_2(x)}\vspace{\eqnspace},\notag
\end{align}
where $A_i(x),\:i=0\ldots2$, are as defined in equation~(\ref{eq:primary}), $G_i(x,s), i=0\ldots 5$, are given in the \cref{sec:appendix}.

\subsubsection{PIE Representation}
We are now ready to prove the following theorem.
\begin{theorem}\label{thm:transform}
 If the matrix 
 \begin{equation}
 B_T=BT
 \end{equation}
  is invertible, where $T$ is given by
 \begin{equation}\label{eq:bmat}
T:=\bmat{I_{n_1}&0 &0\\I_{n_1}&0&0 \\0&I_{n_2} &0\\0&I_{n_2}&(b-a)I_{n_2}\\0&0&I_{n_2}\\0&0 &I_{n_2}},
\end{equation}
then for any $\mbf u^h(x,t)\in X^h$ there exists a fundamental state $\mbf u_f(x,t)\in L_2^{ns}$ given by (\ref{eq:fstate}), such that 
$\mbf u^{h}(x,t)$ can be obtained from $\mbf u_f(x,t)$ by a transformation
\begin{equation}\label{eq:maph}
\mathbf{u}^{h}(x,t)=K(x)B_T^{-1}\mbf{h}(t)+\mcl {T}\mbf u_f (x,t),
\end{equation}
with $\mcl T$ as defined in (\ref{eq:3-pi}), and
$K(x)$ given in~\cref{sec:appendix}. Furthermore, for any $\mbf u_f(x,t)\in L_2^{ns}$, $\mathbf{u}^{h}(x,t)$ obtained via (\ref{eq:maph}) is in $X^h$.
\end{theorem}
\begin{proof}
Suppose $\mbf u^h(x,t)\in X^h$. Define the corresponding fundamental state $\mbf u_f(x,t)$ via (\ref{eq:fstate}). Clearly, $\mbf u_f(x,t)\in L_2^{ns}$. Using~\cref{lem:identity1}, we can express  $\mbf u_{bf}(t)$ through  $\mbf u_{bc}(t)$ (see equation (\ref{eq:ubfubc})) and the fundamental state $\mbf u_f (x,t)$ given by~(\ref{eq:fstate}) as
\begin{equation}\label{eq:ubf}
\mbf u_{bf}(t)=T\mbf u_{bc}(t)+\mcl P_{\{0,0,Q\}} \mbf u_f (x,t),
\end{equation}
where $T$ is given by (\ref{eq:bmat}), and $Q$ is defined in~\cref{sec:appendix}. 
Analogously, the primary state $\mbf u^{h}(x,t)$ can be expressed through $\mbf u_{bc}(t)$ and $\mbf u_f (x,t)$ as
\begin{equation}\label{eq:up}
\mbf u^{h}(x,t)=K(x)\mbf u_{bc}(t)+\mcl P_{\{G_0,G_1,0\}} \mbf u_f (x,t),
\end{equation}
where
$G_0, G_1$ are as defined in~\cref{sec:appendix}. Using~(\ref{eq:ubf}), the boundary constraint (\ref{eq:bcon}) can be expressed as
\begin{equation}\label{eq:identity}
B\mbf u_{bf}(t)=B_T\mbf u_{bc}(t)+B\mcl P_{\{0,0,Q\}} \mbf u_f (x,t),
\end{equation}
from where, since $B\mbf u_{bf}(t)=\mbf h(t)$, we have
\begin{equation}
B_T\mbf u_{bc}(t)+B\mcl P_{\{0,0,Q\}} \mbf u_f (x,t)=\mbf h(t).
\end{equation}
Using the assumption of invertibilty of $B_T$, we may now express the core boundary condition vector as 
\begin{align}\label{eq:ubc}
\mbf u_{bc}(t)&=B_T^{-1}\mbf h(t)-B_T^{-1}B\mcl P_{\{0,0,Q\}} \mbf u_f (x,t)  \\ &
=B_T^{-1}\mbf h(t) -\mcl P_{\{0,0,B_T^{-1}BQ\}}\mbf u_f (x,t) \notag. 
\end{align}
Substituting (\ref{eq:ubc}) into (\ref{eq:up}), we get
\begin{align}
\mbf u^{h}(x,t)= K(x)B_T^{-1}\mbf h(t)  -\mcl P_{\{K,0,0\}}\mcl P_{\{0,0,B_T^{-1}BQ\}}\mbf u_f (x,t)+ \\ \mcl P_{\{G_0,G_1,0\}} \mbf u_f (x,t) 
=K(x)B_T^{-1}\mbf h(t)+\mcl P_{\{G_0,G_1,G_2\}}\mbf u_f (x,t) \notag,
\end{align}
which concludes the proof of the first part of the theorem. Note that the addition rule, scalar multiplication rule and the composition rule for the 3-PI operators, given in~\cref{sec:composition}, were used in this proof.

Conversely, let $\mbf u_f(x,t)$ be in $L_2^{ns}$. It is proven in~\cite{peet_arxiv_PDE} that $\mcl T \mbf u_f(x,t)\in X^0$. Therefore,  $\mcl T\mbf u_f(x,t)\in X$, since $X^0\subset X$. It is easy to see that $K(x)B_T^{-1}\mbf{h}(t)\in H^{\infty\,ns}$, therefore  $K(x)B_T^{-1}\mbf{h}(t)\in X$, and $\mbf u^h(x,t)\in X$. We now only need to show that $\mbf u^h(x,t)$ satisfies boundary conditions (\ref{eq:bcon}). 
 We may evaluate the value of components $\mbf u_1^h(x,t)$, $\mbf u_2^h(x,t)$ from (\ref{eq:maph}) using the definition of $K(x)$ and $\mcl T$. Correspondingly, we have
\begin{align}
\mbf u_1^h(x,t)&=\bmat{I_{n_1} &0 &0} B_T^{-1} \mbf h(t)- \bmat{0 & I_{n_1} &0} \mcl P_{\{0,G_1,G_2\}} \mbf u_{f} (x,t),\label{eq:u1}\\
\mbf u_2^h(x,t)&=\bmat{0 &I_{n_2} &(x-a) I_{n_2}}   B_T^{-1}\mbf h(t)-\bmat{0 & 0 &I_{n_2} }  \mcl P_{\{0,G_1,G_2\}} \mbf u_{f} (x,t)\label{eq:u2}.
\end{align}
Furthermore, differentiating (\ref{eq:u2}) with respect to $x$, we get
\begin{equation}
\mbf u_{2x}^h(x,t)=\bmat{0 &0 & I_{n_2}}   B_T^{-1}\mbf h(t)-\frac{ \partial}{{\partial x}}\, \left(\bmat{0 & 0 & I_{n_2}}\mcl P_{\{0,G_1,G_2\}} \mbf u_{f} (x,t)\right)\label{eq:u2x}.
\end{equation}
Now, evaluating (\ref{eq:u1}), (\ref{eq:u2}), (\ref{eq:u2x}) at $x=a$ nullifies the contribution of $\mcl P_{\{0,G_1,0\}}$ operator and gives us the boundary conditions vector $\mbf u_{bc}(t)$ as
\begin{equation}\label{xbc2}
 \mbf u_{bc}(t)=\bmat{ \mathbf{u}^h_1(a,t)  \\ \mathbf{u}^h_2(a,t) \\ \mathbf{u}_{2x}^h(a,t)}=\bmat{I_{n_1} &0 &0 \\0 &I_{n_2} & 0\\0 &0 & I_{n_2}}B_T^{-1}\mbf h(t)-B_T^{-1}B\,\mcl P_{\{0,0,Q\}}\mbf u_f(x,t),
 \end{equation}
 see also~\cite{peet_arxiv_PDE}. Now, multiplying both sides of (\ref{xbc2}) by $B_T$ shows that 
 $B_T\mbf u_{bc}(t)+B\mcl P_{\{0,0,Q\}} \mbf u_f (x,t)=\mbf h(t)$, which, by identity (\ref{eq:identity}) proves that the primary state $\mbf u^h (x,t)$ constructed via the transformation (\ref{eq:maph}) satisfies the boundary conditions.
\end{proof}
We also have the following corollary that further establishes the properties of the transformation (\ref{eq:maph}). 
\begin{corollary}\label{cor:surjection}
A transformation $L_2^{ns}\rightarrow  X^{h}$ defined by equation (\ref{eq:maph}) is a surjection. 
\end{corollary}
\begin{proof}
Since, by~\cref{thm:transform}, for every $\mbf u^h(x,t)\in X^h$ there exists $\mbf u_f(x,t)\in L_2^{ns}$ that can be mapped into $\mbf u^h(x,t)$, this shows that (\ref{eq:maph}) is a surjection.
\end{proof}
Another corollary allows to view the transformation (\ref{eq:maph}) as a sequence of a linear and an affine transformation.
\begin{corollary}
A transformation $L_2^{ns}\rightarrow  X^{h}$ defined by equation (\ref{eq:maph}) can be viewed as a sequence of transformations $L_2^{ns}\underbrace{\rightarrow}_\mcl T  X^{0}\underbrace{\rightarrow}_\mcl R X^{h}$, where the transformation $\mcl T: L_2^{ns}\rightarrow  X^{0}$ is a unitary map, and a transformation $\mcl R: X^{0}\rightarrow  X^{h}$ is an affine isomorphism defined by a translation.
\end{corollary}
\begin{proof} 
Denote $\mbf u^0(x,t)=\mcl T \mbf u_f(x,t)$. From~\cite{peet2020partial, peet_arxiv_PDE}, we see that $\mbf u^0(x,t)\in X^0$. Since $X^0$ is a special case of $X^h$ with $\mbf h(t)=0$, \cref{cor:surjection} shows that $\mcl T: L_2^{ns}\rightarrow  X^{0}$ is a surjection (an alternative proof can be found in (\cite{peet_arxiv_PDE}).  Since, by~\cref{lemma:space}, $X^0$ is a linear subspace, an inner product can be defined. References~\cite{peet2020partial, peet_arxiv_PDE} further show that $\mcl T$ preserves the inner products, and thus is a unitary map. 

Now, we define $R(x,t)=K(x)B_T^{-1}\mbf{h}(t)$, such that $\mcl R:  X^{0}\rightarrow  X^{h}$ is given by $\mbf u^h(x,t)=\mbf u^0(x,t)+R(x,t)$, which is an affine transformation of translation. Given a specific vector of boundary conditions $\mbf h(t)$ that fixes $X^h$, a translation function $R(x,t)$ is uniquely defined. We now show that $\mcl R$ is isomorphism. Let $\mbf u^h(x,t)$ be in $X^h$. \cref{thm:transform} shows that $\mbf u^0(x,t)=\mbf u^h(x,t)-R(x,t)$ is in $X^0$, and thus $\mcl R:  X^{0}\rightarrow  X^{h}$ is a surjection. Now, we have to show that $R$ is also an injection. Suppose there are two elements in $X^0$, $\mbf u^0_1(x,t)$ and $\mbf u^0_2(x,t)$ that are mapped into a single element $\mbf u^h(x,t)$. We then have $\mbf u^0_1(x,t)=\mbf u^h(x,t)-R(x,t)$, and $\mbf u^0_2(x,t)=\mbf u^h(x,t)-R(x,t)$. Since $R(x,t)$ is a unique function for every $X^h$, this shows that $\mbf u^0_1(x,t)=\mbf u^0_2(x,t)$, and thus $\mcl R$ is an injection. Hence, $\mcl R$ is an isomorphism, as desired. 
\end{proof}
We are now ready to state the final result concerning the conversion of PDEs with inhomogeneous boundary conditions to the PIE framework.
\begin{theorem}\label{th:pie}
The function $\mbf u^h(x,t)\in X^{h}$ satisfies the PDE equation (\ref{eq:primary}) with  boundary conditions (\ref{eq:bcon}) and initial conditions $\mbf u^h(x,0)=\pmb \beta^h(x)$, $\pmb \beta^h(x)\in X^h$, if and only if the corresponding fundamental state function $\mbf u_f(x,t)=\mcl D\,\mbf u^h(x,t)\in L_2^{ns}$ satisfies the following PIE equation 
\begin{equation}\label{eq:piein}
\mcl T\: \frac{\partial \mbf u_f(x,t)}{\partial \,t}=\mcl A \:\mbf u_f(x,t)+\mbf g(x,t),
\end{equation}
with $\mbf g(x,t)$ given by
\begin{align}\label{eq:g}
\mbf g(x,t)&=A_0(x){K}(x)B_T^{-1}\, \mbf{h}(t) \\ &+A_1(x) \bmat{0_{n_1\times n_1} & 0_{n_1\times n_2} &0 \\0 & 0 & I_{n_2}}B_T^{-1}\mathbf{h}(t)-{K}(x)B_T^{-1}\frac{d\,\mbf{h}(t)}{d\,t}+\mbf f(x,t)\notag,
\end{align}
initial conditions $\mbf u_f(x,0)=\pmb \beta_f(x)$, where $\pmb \beta_f(x)=\mcl D \,\pmb \beta^h(x)$, and the 3-PI operators $\mcl T$, $\mcl A$ as defined by (\ref{eq:3-pi}). Moreover, $\mbf u^h(x,t)$ is related to $\mbf u_f(x,t)$ by the transformation (\ref{eq:maph}), and $\pmb {\beta}^{h}(x)=K(x)B_T^{-1}\mbf{h}(0)+\mcl {T}\pmb \beta_f (x)$.
\end{theorem}
\begin{proof}
Suppose $\mbf u^h(x,t)\in X^h$ satisfies the PDE (\ref{eq:primary}) with boundary conditions (\ref{eq:bcon}) and initial conditions (\ref{eq:incon}). Since $\mbf u_f(x,t)=\mcl D\,\mbf u^h(x,t)$, it immediately follows that $\mbf u_f(x,0)=\mcl D\,\mbf u^h(x,0)$, i.e. $\pmb \beta_f(x)=\mcl D\,\pmb \beta^h(x)$. Using the definition of the PDE (\ref{eq:primary}) and defining an auxiliary differentiation operator $\mcl D_1$ as
\begin{equation}
\mathcal{D}_1:=\bmat{0_{n_1\times n_0}&I_{n_1} \partial_x&0\\0&0 &I_{n_2}\,\partial_x},
\end{equation}
we get
\begin{align}\label{pdeproof1}
\frac{\partial \mbf u^h(x,t)}{\partial \,t}&=\mcl P_{\{A_0,0,0\}}\mbf u^h(x,t)+\mcl P_{\{A_1,0,0\}}\,\mcl D_1\,\mbf u^h(x,t)\\&+\mcl P_{\{A_{20},0,0\}}\,\mcl D\,\mbf u^h(x,t)+\mbf f(x,t)\notag.
\end{align}
To evaluate $\mcl D\, \mbf u^h(x,t)$, equation (\ref{eq:fstate}) can be used, while $\mcl D_1 \mbf u^h(x,t)$ can be obtained from
\begin{equation}\label{eq:d1}
\mcl D_1 \mbf u^h(x,t)=\mcl D_1 \mcl{P}_{\{\tilde{K},0,0\}}\mbf{h}(t)+\mcl D_1\mcl T \mathbf{u}_f(x,t),
\end{equation}
where the notation $\tilde{K}(x)=K(x)B_T^{-1}$ is used.
Substituting (\ref{eq:maph}), (\ref{eq:fstate}) and (\ref{eq:d1}) into (\ref{pdeproof1}), we obtain
\begin{align}
\frac{\partial \mbf u^h(x,t)}{\partial \,t}=\mcl P_{\{A_0,0,0\}}\mcl{P}_{\{\tilde{K},0,0\}}\mbf{h}(t)+\mcl P_{\{A_0,0,0\}}\mcl T \mathbf{u}_f(x,t)\notag\\+\mcl P_{\{A_1,0,0\}}\mcl D_1 \mcl{P}_{\{\tilde{K},0,0\}}\mbf{h}(t)+\mcl P_{\{A_1,0,0\}}\mcl D_1\mcl T \mathbf{u}_f(x,t)\\+\mcl P_{\{A_{20},0,0\}}\, \mcl D \, \mcl{P}_{\{\tilde{K},0,0\}}\mbf{h}(t)+\mcl P_{\{A_{20},0,0\}}\mcl D \,\mcl T \mathbf{u}_f(x,t)+\mbf f(x,t)\notag.
\end{align}
Separating homogeneous and non-homogeneous terms in the right-hand side, we have
\begin{equation}\label{eq:separ}
\frac{\partial \mbf u^h(x,t)}{\partial \,t}=H(x,t)+I(x,t), 
\end{equation}
where
\begin{equation}\label{eq:homog}
H(x,t)=\mcl P_{\{A_0,0,0\}}\mcl T \mathbf{u}_f(x,t)+\mcl P_{\{A_1,0,0\}}\mcl D_1 \mcl T \mathbf{u}_f(x,t)+\mcl P_{\{A_{20},0,0\}} \mathbf{u}_f(x,t),
\end{equation}
\begin{equation}\label{eq:inhomog}
I(x,t)=\mcl P_{\{A_0,0,0\}}\mcl{P}_{\{\tilde{K},0,0\}}\,\mbf{h}(t)+\mcl P_{\{A_1,0,0\}}\mcl D_1 \mcl{P}_{\{\tilde{K},0,0\}}\mbf{h}(t)+\mbf f(x,t).
\end{equation}
Homogeneous term, as shown in \cite{peet2020partial}, reduces to
\begin{equation}\label{eq:h}
H(x,t)=\mcl P_{\{H_0,H_1,H_2\}}\mbf u_f(x,t)=\mcl A\,\mbf u_f(x,t).
\end{equation} 
Finally, taking a partial derivative with respect to time of equation (\ref{eq:maph}), we have
\begin{equation}\label{eq:partialtime}
\frac{\partial \mbf u^h(x,t)}{\partial\,t}=K(x)B_T^{-1}\frac{d\,\mbf h(t)}{d\,t}+\mcl T\, \frac{\partial \mbf u_f(x,t)}{\partial\,t}.
\end{equation}
Combining equations (\ref{eq:separ})--(\ref{eq:partialtime}) leads to (\ref{eq:piein})--(\ref{eq:g}).

Conversely, suppose $\mbf u_f(x,t)\in L_2^{n_s}$ satisfies the PIE equation (\ref{eq:piein})--(\ref{eq:g}) with initial conditions $\mbf u_f(x,0)=\pmb \beta_f(x)$. Define $\mbf u^h(x,t)$ according to the transformation (\ref{eq:maph}). By~\cref{thm:transform}, $\mbf u^h(x,t)\in X^h$, and thus satisfies boundary conditions (\ref{eq:bcon}). Furthermore, evaluating (\ref{eq:maph}) evaluated at $t=0$ gives $\pmb {\beta}^{h}(x)=K(x)B_T^{-1}\mbf{h}(0)+\mcl {T}\pmb \beta_f (x)$. Rearrange the PIE equation as 
\begin{equation}\label{eq:pieproof}
\mcl T\: \frac{\partial \mbf u_f(x,t)}{\partial \,t}+K(x)B_T^{-1}\frac{d\,\mbf{h}(t)}{d\,t}=\mcl A \:\mbf u_f(x,t)+I(x,t),
\end{equation}
with $I(x,t)$ as defined in (\ref{eq:inhomog}). The left-hand side of the equation (\ref{eq:pieproof}) equals to a partial time derivative of $\mbf u^h(x,t)$, $\partial \,\mbf u^h(x,t)/\partial\,t$, according to (\ref{eq:partialtime}). Recognizing that, by (\ref{eq:h}), $\mcl A\,\mbf u_f(x,t)=H(x,t)$, and using equations (\ref{eq:homog}) and    (\ref{eq:inhomog}), the right-hand side of (\ref{eq:pieproof}) becomes
\begin{align}\label{eq:hplusi}
H(x,t)+I(x,t)=\mcl P_{\{A_0,0,0\}}\left(\mcl T \mathbf{u}_f(x,t)+\mcl{P}_{\{\tilde{K},0,0\}}\,\mbf{h}(t)\right)\notag \\+\mcl P_{\{A_1,0,0\}}\left(\mcl D_1 \mcl T \mathbf{u}_f(x,t)+\mcl D_1 \mcl{P}_{\{\tilde{K},0,0\}}\mbf{h}(t)\right)\\+\mcl P_{\{A_{20},0,0\}} \mathbf{u}_f(x,t)+\mbf f(x,t)\notag.
\end{align}
Using (\ref{eq:fstate}), (\ref{eq:maph}) and (\ref{eq:d1}), the right-hand side of (\ref{eq:hplusi}) reduces to
\begin{align}
H(x,t)+I(x,t)\\=
\mcl P_{\{A_0,0,0\}}\mbf u^h(x,t)+\mcl P_{\{A_1,0,0\}}\,\mcl D_1\,\mbf u^h(x,t)+\mcl P_{\{A_{20},0,0\}}\,\mcl D\,\mbf u^h(x,t)+\mbf f(x,t)\notag,
\end{align}
which is equivalent to the right-hand side of the PDE equation (\ref{eq:primary}), showing that $\mbf u^h(x,t)$ indeed satisfies the original PDE. 
\end{proof}

\subsubsection{Note on invertibility of $B_T$}
\cref{thm:transform} relies on the condition of invertibility of the $B_T$ matrix. 
It was proven in~\cite{peet2020partial} that the necessary and sufficient condition for the inverse of $B_T$ to exist is for $B$ to: 1) have a row rank of $n_1+2\,n_2$, and 2) have a row space that has a trivial intersection with the row space of $T^{\perp}$, where $T^{\perp}$ defines an orthogonal complement to a column space of $T$. This leads to an exclusion of the boundary conditions that 
are a linear combination of 
\begin{align}
\mbf u_1(a,t)-\mbf u_1(b,t) &=\mbf h_1(t),\label{eq:cond1}\\ \mbf u_2(a,t)+(b-a)\mbf u_{2x}(a,t)-\mbf u_2(b,t)&=\mbf h_2^{(1)}(t), \\ \mbf u_{2x}(a,t)-\mbf u_{2x}(b,t)&=\mbf h_2^{(2)}(t),
\end{align} 
from the set of the boundary conditions, for which $B_T$ is invertible. Here, 
$\mbf h(t)=[\mbf h_1(t)^T\: \:\mbf h_2^{(1)}(t)^T\: \:\mbf h_2^{(2)}(t)^T]^T$,  $\mbf h_1(t)\in \mathbb R^{n_1}$,  $\mbf h_2^{(1)}(t)\in \mathbb R^{n_2}$,  $\mbf h_2^{(2)}(t)\in \mathbb R^{n_2}$.
Note that the excluded boundary conditions involve periodic boundary conditions on the state $\mbf u_1(x,t)$, periodic boundary conditions on derivatives of the state $\mbf u_2(x,t)$, and  Neumann-Neumann conditions for the state $\mbf u_{2}(x,t)$, among others. In general, such boundary conditions are ill-posed for the boundary value problems, however, they typically result in unique solutions to initial-boundary value problems due to a regularization by initial conditions. In a PIE framework, the problems with $B_T$ invertibility for these boundary conditions arise from the fact that now a fundamental state needs to have an additional constraint in order to satisfy these boundary conditions, implying that the fundamental state is no longer minimal. For example, with the periodic boundary condition on a function, we have a constraint that the integral of its derivative over the domain must be equal to zero. If this derivative enters the fundamental state, as would be the case for $\mbf u_{1x}$ with a periodic state $\mbf u_1$, this additional constraint, since it is not embedded into the PIE dynamics, can not be satisfied. 

To remedy this situation, it is possible to redefine a fundamental state to be free of constraints, and embed the corresponding constraints into the PIE operators. This can be formally accomplished by performing an SVD decomposition of the $B_T$ matrix, introducing an auxiliary state vector $\mbf u_n(t)\in \mathbb R^r$, where $r$ is the rank deficiency of $B_T$, and modifying the PIE equations accordingly~\cite{shivakumar2021extension}. While this is generally possible, such modification will not be considered here, and we will assume that $B_T$ matrix is invertible, with the use of apposite boundary conditions.

\section{Solution of the PDEs in the PIE Framework: PIE-Galerkin approximation}\label{sec:numsol}
\subsection{Spatial treatment}
We are now interested in finding a solution $\mbf u_f(x,t)\in L_2[a,b]^{ns}$ to the PIE equation~(\ref{eq:piein}) with the initial conditions  $\mbf u_f(x,0)=\pmb{\beta}_f(x)$, which, according to \cref{th:pie}, satisfies the original PDE equation~(\ref{eq:primary}). Since  $\mbf u_f(x,t)\in L_2[a,b]^{ns}$, we are free to choose any approximation space without needing to worry about satisfying boundary conditions. We choose Chebyshev polynomials of the first kind as the approximation functions. Since Chebyshev polynomials are defined on the $[-1,1]$ domain, we need to map our original PDE from $x=[a,b]$ onto a computational domain $x^{(c)}=[-1,1]$, which can be readily accomplished  by a linear transformation 
\begin{equation}\label{eq:mapxi}
x^{(c)}=\frac{2x-(b+a)}{b-a},
\end{equation}
with the inverse map 
\begin{equation}\label{eq:mapfromxi}
x=\frac{b-a}{2}\,x^{(c)}+\frac{b+a}{2}. 
\end{equation}
With a slight abuse of notation, 
in what follows, we will assume 
that the corresponding PIE equations are defined on $x\in [-1,1]$ domain,
acknowledging that necessary transformations might had to be done to the original PDE in order to accomplish this.

In accordance with (\ref{eq:fvector}), (\ref{eq:fstate}), and (\ref{eq:difopp}),  we can write for each sub-component $\mbf u_{fp}(x,t)$ of $\mbf u_{f}(x,t)$, $p=0, 1  2$,
\begin{equation}
\mbf u_{fp}(x,t)=\frac{\partial^{\,p} \mbf u_{p}(x,t)}{\partial x^p}.
\end{equation}
Therefore, with each component $u_{fi}(x,t)$, $i=1\ldots n_s$, of the vector $\mbf u_{f}(x,t)$, we can associate an index
\begin{equation}\label{eq:ds}
 p=p(i),
 \end{equation}
 defined as a ``minimum smoothness'' required from the original $u_{i}(x,t)$  function to enter the PDE (\ref{eq:primary}).
We now look for solutions $u_{fi}(x,t)\in \mathbb P[-1,1]^{N-p(i)}$ for each corresponding $u_{fi}(x,t)$ component, where $\mathbb P[-1,1]^{N-p(i)}$ is the space of all polynomial functions of degree $N-p(i)$ or less on $[-1,1]$ domain, i.e. we approximate
\begin{equation}\label{eq:chebapprox}
\hat{u}_{fi} (x,t)=\sum_{k=0}^{N-p(i)} a_{ik}(t) T_k(x), 
\end{equation}
where $T_k(x)$ are the Chebyshev polynomials of the first kind~\cite{canuto1988spectral}, and $a_{ik}(t)\in C^1(R^+)$ are the corresponding time-dependent Chebyshev coefficients, where the subscript $i$ denotes their affiliation with a particular solution component $\hat{u}_{fi}(x,t)$. 
The approximation for the vector-valued function $\hat{\mbf u}_{f} (x,t)$ can then be compactly written as
\begin{equation}\label{eq:chebex}
\hat{\mbf u}_{f} (x,t)=\sum_{i=1}^{ns} \sum_{k=0}^{N-p(i)} a_{ik}(t)\, \pmb{\phi}_{ik}(x), 
\end{equation}
where the vector-valued Chebyshev basis functions $\pmb{\phi}_{ik}(x):R\rightarrow \mathbb R^{n_s}$ can be defined as 
\begin{equation}\label{eq:basisfunction}
\pmb{\phi}_{ik}(x)=\underbrace{\bmat{ 0 &\cdots & \cdots & T_k(x) & \cdots  & 0}^T}_{n_s},
\end{equation}
where $T_k(x)$ is in the $i^{th}$ position of the vector $\pmb{\phi}_{ik}(x)$, $i=1\ldots n_s$, $k=0 \ldots N-p(i)$. We denote the polynomial space spanned by the vector-valued basis functions  $\pmb{\phi}_{ik}(x)$ as $Y^{N_p}:=\mathbb P[-1,1]^{N_p}$, where $N_p=n_0N\times n_1(N-1)\times n_2(N-2)$, so that the composite vector-valued approximation $\hat{\mbf u}_{f} (x,t)\in Y^{N_p}$.

We introduce the same approximation for the lumped inhomogeneous term $\mbf g (x,t)$, see (\ref{eq:g}), i.e. we write
\begin{equation}\label{eq:lumpedapprox}
\hat{\mbf g} (x,t)=\sum_{i=1}^{ns} \sum_{k=0}^{N-p(i)} {b}_{ik} (t) \,\pmb{\phi}_{ik}(x),
\end{equation}
where $b_{ik}(t)$ are the corresponding Chebyshev coefficients associated with the inhomogeneous term, $\hat{\mbf g} (x,t)\in Y^{N_p}$. 

With the expansion (\ref{eq:chebex}), the action of the 3-PI operator  $\mathcal{T}$ on the function approximation $\hat{\mbf u}_{f} (x,t)\in Y^{N_p}$ can be written as
\begin{equation}\label{eq:chebcol}
\mcl T \hat{\mbf u}_{f} (x,t)= \sum_{i=1}^{ns} \sum_{k=0}^{N-p(i)} a_{ik}(t)\, \mcl T \pmb{\phi}_{ik}(x)=\sum_{i=1}^{ns} \sum_{k=0}^{N-p(i)} a_{ik}(t)\, \textrm{Col}_{\,i} (\mcl T) T_k(x), 
\end{equation}
where the notation $\textrm{Col}_{\,i} (\mcl T)$ stands for the $i^{th}$ column of the matrix operator $\mathcal{T}$. We now have the following lemma.

\begin{lemma}\label{lem:cheb}
The action of any element $\mathcal{T}_{mn}$ of the 3-PI operator $\mathcal{T}$ on a Chebyshev polynomial function $T_k(x)$ can be established according to the following rules:
\begin{enumerate}
\item If $\mathcal{T}_{mn}$ is such that $m\leq n_0$, the action is calculated as 
\begin{equation}\label{eq:mult}
\mcl T_{mn} T_{k}(x)=\delta_{mn} T_k(x).
\end{equation}
\item If $\mathcal{T}_{mn}$ is such that $n_0 < m\leq n_0+n_1$, the action is calculated as 
\begin{equation}\label{eq:n1tran}
\mcl T_{mn} T_{k}(x)=b^{(1)}_{0kmn} T_0(x)+b^{(1)}_{1kmn} T_1(x)+\delta_{mn}\left(c^-_{k-1}T_{k-1}(x)+c^+_{k+1}T_{k+1}(x)\right),
\end{equation}
where
\noindent\begin{minipage}{.45\linewidth}
\begin{align}\label{eq:n1coef}
c^-_k=\begin{cases}
0,& k\le 1 \\
-\frac{1}{2k},& k\ge 2
\end{cases}
\end{align}
\end{minipage}%
\begin{minipage}{.45\linewidth}
\begin{align}
c^+_k=\begin{cases}
0,& k\le 1 \\
\frac{1}{2k},& k\ge 2\notag
\end{cases}
\end{align}
\end{minipage}
\item If $\mathcal{T}_{mn}$ is such that $m>n_0+n_1$, the action is calculated as 
\begin{align}\label{eq:n2tran}
\mcl T_{mn} T_{k}(x)=b^{(2)}_{0kmn} T_0(x)+b^{2}_{1kmn} T_1(x) \\ +\delta_{mn} \left(d^{-}_{k-2}T_{k-2}(x)+d_{k}T_{k}(x)+d^{+}_{k+2}T_{k+2}(x)\right),\notag
\end{align}
where
\noindent\begin{minipage}{.45\linewidth}
\begin{align}\label{eq:n2coef}
d_k^-=\begin{cases}
0,& k\le 1 \\
\frac{1}{4\,k(k+1)}, & k\ge 2
\end{cases}
\end{align}
\end{minipage}%
\begin{minipage}{.45\linewidth}
\begin{align}
d_k^+=\begin{cases}
0,& k\le 1 \\
\frac{1}{2\,k(k-1)}, & k= 2\notag\\
\frac{1}{4\,k(k-1)}, & k\ge 3\notag
\end{cases}
\end{align}
\end{minipage}
\begin{align}
d_k=\begin{cases}
0,& k\le 1\\
-\frac{1}{2(k^2-1)}, & k\ge 2\notag,
\end{cases}
\end{align}
\end{enumerate}
where $\delta_{mn}$ is a Kronecker delta function, $b^{(i)}_{jkmn}$, $i=1,2,\,j=0,1$ are real constants, generally dependent on boundary conditions, and $c_k^-, c_k^+, d_k^-, d_k, d_k^+$ are real constants not dependent on boundary conditions.
\end{lemma}
\begin{proof}
The proof of this lemma is included in the~\cref{sec:app-proof}.
\end{proof}
As a consequence of this result, it can be concluded that the action of $\mcl T$ on functions that belong to polynomial subspaces, keeps them in polynomial subspaces, which allows us to evaluate the action of a partial-integral operator $\mcl T$ on the polynomial functions analytically, using the formulas presented in~\cref{lem:cheb}.
In fact, it allows us to prove the following lemma.
\begin{lemma}\label{lemma:chebyshevmap}
If $\hat{\mbf u}_f (x,t)\in Y^{N_p}, \,N_p=n_0N\times n_1(N-1)\times n_2(N-2), \,t\in\mathbb{R}^+$, $N\ge 2$, the corresponding function approximation $\hat{\mbf u}^{h} (x,t)$ to the primary solution 
\begin{equation}\label{eq:discreteop}
\hat{\mathbf{u}}^{h}(x,t)=K(x)B_T^{-1}\mbf{h}(t)+\mcl {T}\hat{\mbf u}_f (x,t)
\end{equation}
is in the space $\mathbb{P}^{Nn_s}, t\in \mathbb{R}^+$, i.e. all the components of the primary vector-valued solution are in $\mathbb{P}^{N}$. Furthermore, for $\hat{\mbf u}^{h} (x,t)\in \mathbb{P}^{Nn_s}$, the corresponding fundamental state approximation 
\begin{equation}\label{duhat}
\hat{\mbf u}_{f} (x,t)=\mcl D \,\hat{\mbf u}^{h} (x,t)
\end{equation}
 is in $Y^{N_p}$.
\end{lemma}
\begin{proof}
Suppose $\hat{\mbf u}_f(x,t)\in Y^{N_p}$. We first note that $K(x)B_T^{-1}\mbf{h}(t)\in \mathbb{P}^{1\cdot n_s}$, which, for $\mcl {T}\hat{\mbf u}_f (x,t)\in\mathbb{P}^{Nn_s}$, $N\ge 2$, keeps the composite function in $\mathbb{P}^{Nn_s}$. We now proceed to show that  $\mcl {T}\hat{\mbf u}_f (x,t)\in\mathbb{P}^{Nn_s}$. Denote 
\begin{equation}
\hat{\mbf u}_f(x,t)=\bmat{ \hat{\mathbf{u}}_{f0}(x,t)  \\ \hat{\mathbf{u}}_{f1}(x,t) \\ \hat{\mathbf{u}}_{f2}(x,t)}, \hat{\mbf u}^{0}(x,t)=\mcl T \hat{\mbf u}_f(x,t)=\bmat{ \hat{\mathbf{u}}^{0}_{0}(x,t)  \\ \hat{\mathbf{u}}^{0}_{1}(x,t)  \\ \hat{\mathbf{u}}^{0}_{2}(x,t)}, \end{equation}
where $\hat{\mathbf{u}}_{fp}(x,t)$ is the polynomial approximation of  $\mathbf{u}_{f p}(x,t)\in L_2^{n_p}$, and $\hat{\mbf u}^0_{p}(x,t)$ is the polynomial approximation of $\mbf u^{0}_{p}(x,t)$, respectively, $p=0,1, 2$, $\mbf u^{0}(x,t)\in X^0$. Noting the structure of the matrix operators $G_0$, $G_1$ and $G_2$, it is easily seen that
\begin{equation}\label{eq:hats}
\bmat{ \hat{\mathbf{u}}^{0}_{0}(x,t)  \\ \hat{\mathbf{u}}^{0}_{1}(x,t)  \\ \hat{\mathbf{u}}^{0}_{2}(x,t)}=\mcl P_{\{G_0,0,0\}}\bmat{ \hat{\mathbf{u}}_{f0}(x,t)  \\ 0 \\ 0}+\mcl P_{\{0,G_1,0\}}\bmat{ 0 \\ \hat{\mathbf{u}}_{f1}(x,t)  \\ \hat{\mathbf{u}}_{f2}(x,t) }+\mcl  P_{\{0,0,G_2\}}\bmat{ 0 \\ \hat{\mathbf{u}}_{f1}(x,t)  \\ \hat{\mathbf{u}}_{f2}(x,t)}.
\end{equation}
The first term in the right-hand side of equation (\ref{eq:hats}) shows that the first $n_0$ components of $\hat{\mbf u}_{f}(x,t)$ are mapped into  the first $n_0$ components of $\hat{\mbf u}^0(x,t)$, with the corresponding $\mathbb{P}^N\rightarrow \mathbb{P}^N$ mapping according to (\ref{eq:mult}). Since the matrix $G_1$ is block-diagonal, and according to (\ref{eq:n1tran}), (\ref{eq:n2tran}), the second term of (\ref{eq:hats}) corresponds to $\mathbb{P}^{N-1}\rightarrow \mathbb{P}^{N}$, and $\mathbb{P}^{N-2}\rightarrow \mathbb{P}^{N}$ mappings of the second $n_1$ and the third $n_2$ components between the vectors $\hat{\mbf u}_f(x,t)$ and $\hat{\mbf u}^{0}(x,t)$, respectively. The last entry of equation (\ref{eq:hats}) corresponds to an integral over an entire domain, and thus, as shown in the proof of~\cref{lem:cheb}, produces only the outputs in $\mathbb{P}^{\,0}$ or $\mathbb{P}^{\,1}$.

Now, let  $\hat{\mathbf{u}}^{h}(x,t)$ be in $\mathbb P^{Nns}$. According to the structure of the differentiation operator $\mcl D$, see equation (\ref{eq:difopp}), it is easy to see that $\mcl D \,\hat{\mathbf{u}}^{h}(x,t)\in Y^{n_0N\times n_1(N-1)\times n_2(N-2)}$, which concludes the proof.
\end{proof}

Define the polynomial space $\mathbb P^{\,h}$ such that the functions that are in $\mathbb P^{\,h}$ are also in $\mathbb{P}^{Nn_s}$, and satisfy the boundary conditions~(\ref{eq:bcon}) mapped onto $[-1,1]$ domain, i.e.
\begin{equation}
\mathbb P^{\,h}:=\left\{\bmat{\mathbf{\hat{u}}_0(x,t)\\\mathbf{\hat{u}}_{1}(x,t)\\\mathbf{\hat{u}}_{2}(x,t)}\in \mathbb P^{Nn_s}\cap B {\scriptsize \bmat{ \mathbf{\hat{u}}_1(-1,t) \\ \mathbf{\hat{u}}_1(1,t) \\ \mathbf{\hat{u}}_2(-1,t) \\ \mathbf{\hat{u}}_2(1,t)\\ \mathbf{
\hat{u}}_{2x}(-1,t) \\ \mathbf{\hat{u}}_{2x}(1,t)}}=\mathbf{h}(t),\,t\in \mathbb R^+\right\}
\end{equation}
The following important theorem allows us to establish the approximation properties of the primary solution $\hat{\mbf u}^{h} (x,t)$ of the PDE (\ref{eq:primary}), given by (\ref{eq:discreteop}).
\begin{theorem}\label{th:discrete}
For every $\hat{\mbf u}^{h} (x,t)\in \mathbb P^{\,h}$, with $N\ge 2$, there exists a corresponding approximation to a fundamental state $\hat{\mbf u}_{f} (x,t)=\mcl D\,\hat{\mbf u}^{h} (x,t)$, \\ $\hat{\mbf u}_{f} (x,t)\in Y^{n_0N\times n_1(N-1)\times n_2(N-2)}, t\in\mathbb{R}^+$, that is mapped into $\hat{\mbf u}^{h} (x,t)$ according to the transformation (\ref{eq:discreteop}). Moreover, for every  $\hat{\mbf u}_{f} (x,t)\in Y^{n_0N\times n_1(N-1)\times n_2(N-2)}$, $\hat{\mbf u}^{h} (x,t)$ defined by (\ref{eq:discreteop}) is in $\mathbb P^{\,h}$.
\end{theorem}
\begin{proof}
Let $\hat{\mbf u}^{h} (x,t)\in \mathbb P^{\,h}$. Therefore, $\hat{\mbf u}^{h} (x,t)\in \mathbb P^{Nn_s}$. Suppose $\hat{\mbf u}_{f} (x,t)$ satisfies Eq.~(\ref{duhat}). 
By~\cref{lemma:chebyshevmap},  $\hat{\mbf u}_{f} (x,t)\in Y^{N_p}$, where $N_p=n_0N\times n_1(N-1)\times n_2(N-2)$. Moreover due to a~\cref{thm:transform}, we have that, since $\mathbb P^{\,h}\subset X^h$, and $Y^{N_p}\subset L_2^{ns}$, $\hat{\mbf u}_{f} (x,t)$ defined by (\ref{duhat}) is mapped into $\hat{\mbf u}^{h} (x,t)$ according to the transformation (\ref{eq:maph}), which is equivalent to (\ref{eq:discreteop}). 

Now, consider any $\hat{\mbf u}_{f} (x,t)\in Y^{N_p}$. Again, by~\cref{lemma:chebyshevmap}, $\hat{\mbf u}^{h} (x,t)$ defined by the transformation (\ref{eq:discreteop}) is in $\mathbb P^{Nn_s}$. We are left to prove that $\hat{\mbf u}^{h} (x,t)$ satisfies the boundary conditions (\ref{eq:bcon}) with $a=-1, b=1$. Since $\hat{\mbf u}_{f} (x,t)\in L_2^{ns}[-1,1]$,~\cref{thm:transform} ensures that $\hat{\mbf u}^{h} (x,t)$ obtained via (\ref{eq:discreteop}), which is equivalent to (\ref{eq:maph}), is in $X^h[-1,1]$, i.e. satisfies the aforementioned boundary conditions, which concludes the proof. 
\end{proof}

We note that the  property given by \cref{th:discrete} could be established due to the fact that 
 the 3-PI operator $\mcl T$ is invariant under a projection onto the polynomial subspace $\mathbb{P}^N$, thus guaranteeing the equivalence of the transformations (\ref{eq:maph}) and (\ref{eq:discreteop}). It would not necessarily hold true for another choice of an approximation space, such as, e.g., with harmonic functions.

To represent the operator $\mcl A=\mcl P_{\{H_0,H_1,H_2\}}$ in the right-hand side of equation~(\ref{eq:piein}), which contains the functions $A_0(x)$, $A_1(x)$, and $A_2(x)$,  in the Chebyshev Galerkin approximation framework, we decompose the functions $A_j(x)$, $j=0,1,2$, into the Chebyshev series as
\begin{equation}
A_{j} (x)=\sum_{m=0}^{\infty} {A}_{jm} T_m(x), 
\end{equation}
where $A_{jm}$ are the matrix-valued coefficients for a particular function $A_j(x)$. Correspondingly, the kernel functions $H_j$, $j=0, 1, 2$, in $\mcl P_{\{H_0, H_1, H_2\}}$ can be decomposed into the  matrix-valued Chebyshev expansion series as
\begin{align}
H_{0} (x)= &\sum_{m=0}^{\infty} {H}_{jm} T_m(x),\\
H_{j} (x,s)= &\sum_{m=0}^{\infty} \sum_{i=0}^{1} {A}_{im} T_m(x) G_{j+3i}(x,s),\:\:j=1,2.
\end{align}

To apply the operator $\mcl A=\mcl P_{\{H_0,H_1,H_2\}}$ to  $\hat{\mbf u}_f(x,t)$ given by (\ref{eq:chebex}), we first note that
\begin{equation}\label{eq:aproduct}
H_{0} (x)T_k(x) =\sum_{m=0}^{\infty} {H}_{0m} T_m(x) T_k(x) =\sum_{m=0}^{\infty} \frac{1}{2}{H}_{0m} \left(T_{m+k} (x)+ T_{|m-k|}(x)\right).
\end{equation}
For the integrative kernels, we note that
\begin{equation}\label{eq:a1product}
\int H_{j}(x,s) \,T_k(s)\,ds =\sum_{m=0}^{\infty}  \sum_{i=0}^{1}  {A}_{im}T_m(x)\int G_{j+3i}(s) T_k(s)\,ds,
\end{equation}
where $j=1,2$, upon which the integrals in the right-hand side of Eq. (\ref{eq:a1product})
can be computed according to the formulas developed in~\cref{sec:app-proof}.

We proceed with applying a method of weighted residuals to the equation~(\ref{eq:piein}), i.e., we introduce a space of test functions $\hat{\mbf v}(x)\in Z^{N_p}$ and demand that, for $\hat{\mbf u}_f(x,t)\in Y^{N_p}$, $t\in R^+$,
\begin{equation}\label{eq:test}
\left(\mcl T\: \frac{\partial \hat{\mbf u}_f(x,t)}{\partial \,t},\hat{\mbf v}(x)\right)=\left(\mcl A \:\hat{\mbf u}_f(x,t)+\hat{\mbf g} (x,t),\hat{\mbf v} (x)\right),\:\forall \hat{\mbf v}(x)\in Z^{N_p},
\end{equation}
with $\left(\hat{\mbf u}_{f}(x,t),\hat{\mbf v}(x)\right),\,t\in \mathbb{R}^+$, denoting an inner product on a Hilbert space defined as 
\begin{equation}\label{eq:inner}
\left(\hat{\mbf u}_{f}(x,t),\hat{\mbf v}(x)\right)=\int_{-1}^1 \hat{\mbf u}_{f}^T(x,t) \hat{\mbf v}(x) w(x) d\,x,\:w(x)=\frac{1}{\sqrt{1-x^2}},
\end{equation} 
where $w(x)$ is the weight function.
Following Galerkin approach, we set $Z^{N_p}=Y^{N_p}$.
Taking an inner product in (\ref{eq:test}) with each basis function $\pmb \phi_{mn}\in Y^{N_p}$, $m=1\ldots n_s$, $n=0 \ldots N-p(m)$,
and using the orthogonality of the Chebyshev polynomials with respect to this weight function~\cite{canuto1988spectral}, a set of $N_d$ linear ordinary differential equations (ODEs) is obtained for $N_d$ unknown Chebyshev coefficients $a_{ik}(t)$ in (\ref{eq:chebex}), $N_d=n_0(N+1)\times n_1N\times n_2(N-1)$, which can be written in a matrix form as
\begin{equation}\label{eq:ode}
M \frac{d\,{\mbf a}(t)}{d\,t}=A\,\mbf{a}(t)+\mbf{b}(t),
\end{equation}
with initial conditions 
\begin{equation}\label{eq:inconode}
\mbf{a}(0)=\pmb a_0.
\end{equation}

Here, $\mbf a(t)\in  \mathbb{R}^{N_d}$ is the vector of the Chebyshev expansion coefficients of the unknown function $\hat{\mbf u}_f(x,t)$ via (\ref{eq:chebex}), and  $\mbf b(t)\in  \mathbb{R}^{N_d}$ is the vector of known Chebyshev coefficients coming from the series expansion of the lumped inhomogeneous term (\ref{eq:g}) via (\ref{eq:lumpedapprox}). To form the $\mbf{a}(t)$ and $\mbf{b}(t)$ vectors, we stack $N-p(i)$ Chebyshev coefficients $a_{ik}(t)$, $b_{ik}(t)$, corresponding to each component $i$,  prior to proceeding to the next component, i.e. the entries $a_j(t)$, $b_j(t)$ of $\mbf{a}(t)$, $\mbf{b}(t)$ can be expressed as $a_{(i-1)ns+k+1}(t)=a_{ik}(t)$, $i=1\ldots n_s, k=0\ldots N-p(i)$, same for $b_j(t)$. Matrices $M\in\mathbb{R}^{N_d\times N_d}$,  $A\in\mathbb{R}^{N_d\times N_d}$ are the matrices consisting of the entries of the discretized $\mcl T$ and $\mcl{A}$ operators, respectively, multiplying the corresponding components of $\mbf a(t)$ vector.  To obtain initial conditions (\ref{eq:inconode}), the corresponding initial conditions $\mbf u_f(x,0)=\pmb \beta_f(x)$ of the PIE equation are projected onto the $Y^{N_p}$ polynomial space as
\begin{equation}
\hat{\pmb \beta}_f(x)=\sum_{i=1}^{ns} \sum_{k=0}^{N-p(i)} a_{0\,ik}\,\pmb{\phi}_{ik}(x),
\end{equation}
$\hat{\pmb \beta}_f(x)\in Y^{N_p}$, and the coefficient vector $\mbf a_0$ of initial conditions is constructed from $a_{0\,ik}$ coefficients in accordance with the procedure outlined above.

The following lemma establishes a sparsity structure of the matrix $M$.
\begin{lemma}\label{lem:chebmat}
Matrix $M$ in the equation (\ref{eq:ode}) has the following structure:
\begin{enumerate}
\item The first $n_0(N+1)$ rows of $M$ are defined by an upper-left  $I_{n_0(N+1)}$ identity matrix, with the rest of the entries being zero both to the right of $I_{n_0(N+1)}$   (i.e. in the first $n_0(N+1)$ rows of $M$) and below $I_{n_0(N+1)}$ (i.e. in the first $n_0(N+1)$ columns of $M$). 
\item The subsequent $n_1\,N$ rows of $M$ consist of $n_1$ tridiagonal blocks of size  $N\times N$, with zeroes on the main diagonal, and the coefficients $c^+_n$, $c^-_{n}$ from (\ref{eq:n1coef}) on the subdiagonal and the superdiagonal in the row $l=(m-1)n_s+n+1$, respectively. The exceptions are the first two rows of each block, which, in general, can be full rows with the real entries in the column positions between $n_0(N+1)+1$ and $n_s$, representing a coupling across states due to the boundary conditions.
\item The last $n_2\,(N-1)$ rows of $M$ consist of $n_2$ pentadiagonal blocks of size  $(N-1)\times (N-1)$, with $d_n$ on the main diagonal, zeroes on the subdiagonal and superdiagonal, and $d^+_{n}$, $d^-_{n}$ from (\ref{eq:n2coef}) on the 2-subdiagonal and 2-superdiagonal in the row $l=(m-1)n_s+n+1$, respectively. The exceptions are the first two rows of each block, which, in general, can be full rows with the real entries in the column positions between $n_0(N+1)+1$ and $n_s$, representing a coupling across states due to the boundary conditions.
\end{enumerate}
\end{lemma}
\begin{proof}
The proof of this lemma is included in the~\cref{sec:app-proofmat}.
\end{proof}

As a consequence of this lemma, it can be seen that the influence of the boundary conditions is felt only in the first two rows in each of the corresponding solution component block of the matrix $M$, which is reminiscent of the characteristics of the Chebyshev tau differentiation and integration methods, albeit the boundary condition structure is embedded into the matrix analytically in the current method, as opposed to discretely in Chebyshev tau methods. The dependence of the matrix $A$ on the boundary conditions is more complex. Since only $H_2(x,s)$ operator in $\mcl A$ contains the matrix $B$, there is no dependence if $A_0(x)=A_1(x)=0$. When  $A_0(x)$ and $A_1(x)$ are present but constant, the topology of the boundary conditions influence in $A$ is the same as in $M$, i.e. only the first two rows in each solution block are effected. However, if   $A_0(x)$ and $A_1(x)$ have variable coefficients, the influence of $B$ propagates into the interior of the matrix $A$ through nonlinear products in (\ref{eq:a1product}), effecting as many additional rows as the degree of nonlinearity of $A_0(x)$, $A_1(x)$. 


\subsection{Stability and Convergence of a semi-discrete approximation}

This section concerns with the stability and convergence estimates of a semi-discrete PIE-Galerkin formulation, namely, when a temporal variable is not discretized. For a sake of brevity, we will consider a scalar case, while extension to a vector-valued case is straightforward. 
Since Eq.~(\ref{eq:piein}) can represent both parabolic and hyperbolic systems, we consider the most conservative situation and, instead of assuming coercivity~\cite{canuto1988spectral, bressan1986analysis}, simply assume a non-positivity property associated with the integral operators $\mcl A, \mcl T$ as
\begin{equation}\label{eq:nonpos}
(\mcl A \,u_f,\mcl T u_f)\le 0\:\: \textrm{for all} \:u_f\in L_2[-1,1]
\end{equation}
with the inner product defined as in (\ref{eq:inner}), and its discrete counterpart
\begin{equation}\label{eq:nonposdis}
(\mcl A \,\hat{u}_f,\mcl T \hat{u}_f)_N\le 0\:\: \textrm{for all} \:\hat{u}_f\in \mathbb{P}[-1,1]^N\, \textrm{and for all}\:N>0,
\end{equation}
where the inner product in the left-hand side of (\ref{eq:nonposdis}) is defined as $(\mcl A \,\hat{u}_f,\mcl T \hat{u}_f)_N=(R_N(\mcl A \,\hat{u}_f),R_N(\mcl T \hat{u}_f))$, with $R_{N}: L_2\rightarrow \mathbb{P}^{N}$ being a projection operator.
The following theorem concerns with a stability of Galerkin approximation of the PIE equation~(\ref{eq:piein}).
\begin{theorem}\label{thm:stability}
Denote $\widehat{\mcl T u_f}=R_{N-p} (\mcl T u_f)$, where $p=0, 1$ or $2$ is defined in (\ref{eq:ds}).
Under the assumption (\ref{eq:nonposdis}), the following inequality holds
\begin{equation}\label{eq:energy}
||\widehat{\mcl T \hat{u}_f(t)}||^2\le C(t)\left(||\widehat{\mcl T \hat{u}_f(0)}||^2+\int_0^t ||\hat{g}(s)||^2\,ds\right)\:\:\textrm{for all}\:\:t\ge0,
\end{equation}
with the constant $C(t)$ independent of $N$, which yields stability of approximation (\ref{eq:test}).
\end{theorem}
\begin{proof}
Estimate (\ref{eq:energy}) is readily obtained from (\ref{eq:test}) by using $\hat{v}=\widehat{\mcl T \hat{u}_f(t)}$ as a test function, assumption (\ref{eq:nonposdis}), Cauchy-Schwarz inequality to bound the inner product $(\hat{g},\widehat{\mcl T \hat{u}_f})\le||\hat{g}|| \,||\widehat{\mcl T \hat{u}_f}||$, algebraic inequality $ab\le 1/(4\epsilon)\,a^2+\epsilon\, b^2$ with $\epsilon=1/2$, and, subsequently, invoking Gronwall's lemma~\cite{canuto1988spectral, canuto1982error, tadmor1994spectral}, yielding $C(t)=\exp(t)$. 
\end{proof}

The following theorem establishes the convergence properties of the PIE-Galerkin methodology.
\begin{theorem}
If (\ref{eq:nonposdis}) is satisfied, the following convergence estimate holds
\begin{align}\label{eq:converg}
||u^h(t)-\widehat{\hat{u}^h(t)}||\le C(N-p)^{p-m}\bigg\{||u_f(t)||+\\ exp\left(\frac{t}{2}\right)\left(\int_0^t \left( ||\dot{u}_f(s)||^2+||u_f(s)||^2+||g(s)||^2\right)ds\right)^{1/2}\bigg\}\:\:\textrm{for all}\:\:t\ge0\notag,
\end{align}
where $p$ is  a minimum smoothness of the primary solution as in \cref{thm:stability}, $m$ is the actual number of square-integrable spatial derivatives of the primary solution, and a dot symbol denotes a partial derivative with respect to time. 
\end{theorem}
\begin{proof}
From (\ref{eq:maph}), (\ref{eq:discreteop}), we have $||u^h(t)-\widehat{\hat{u}^h(t)}||=||\mcl T u_f(t)-\widehat{\mcl T \hat{u}_f(t)}||$. To obtain a convergence estimate, we define an error function $e(x,t)= R_{N-p}u_f(x,t)-\hat{u}_f(x,t)$. Taking an inner product of (\ref{eq:piein}) with $\widehat{\mcl T e}$, substituting $\widehat{\mcl T e}$ as a test function in (\ref{eq:test}), and a subsequent manipulation, the following evolution equation for the error can be obtained:
\begin{equation}\label{eq:error}
\frac{1}{2}\frac{d}{dt}\,\widehat{ \mcl T e(t)}^2=\left(\widehat{\mcl A\, e(x,t)},\,\widehat{\mcl T\,e(x,t)}\right)+\left(R(x,t),\,\widehat{\mcl T\,e(x,t)}\right),
\end{equation}
where the residual term $R(x,t)$ is given by
\begin{align}\label{eq:residual}
R(x,t)= -\left(\mcl T\, \dot{u}_f(x,t)-\widehat{\mcl T\, \dot{u}_f(x,t)}\right)+
\left(\mcl A\, u_f(x,t)-\widehat{\mcl A\, u_f(x,t)}\right)+\left(g(x,t)-\hat{g}(x,t)\right) \notag \\
-\left(\widehat{\mcl T\, \dot{u}_f(x,t)}-\mcl T R_{N-p}(\dot{u}_f(x,t)) \right)+\left(\widehat{\mcl A\, u_f(x,t)}-\mcl A R_{N-p}(u_f(x,t))\right), 
\end{align}
where the last two terms in (\ref{eq:residual}) are errors due to non-commutativity of the integration and projection operators. Applying assumption (\ref{eq:nonposdis}) to the first term in the right-hand side of (\ref{eq:error}), 
bounding the inner product $\left(R(x,t),\,\widehat{\mcl T\,e(x,t)}\right)$ the same way we bounded $(\hat{g},\widehat{\mcl T \hat{u}_f})$ in \cref{thm:stability} and using the Gronwall's lemma, we obtain
\begin{equation}
||\widehat{\mcl T e(t)}||^2\le \exp(t)\left(||\widehat{\mcl T e(0)}||^2+\int_0^t ||R(s)||^2 \,ds \right)\:\:\textrm{for all}\:\:t\ge0.
\end{equation}
We can bound the residual term by noting that, by the properties of the Chebyshev approximation~\cite{canuto1988spectral}, $ ||\mcl T \dot{u}_f(t)-\widehat{\mcl T \dot{u}_f(t)}||\le C_1(N-p)^{-m}||\mcl T \dot{u}_f(t)||\le C_T(N-p)^{-m}||\dot{u}_f(t)||$, $ ||\mcl \,A u_f(t)-\widehat{\mcl \,A u_f(t)}||\le C_2(N-p)^{-m}||\mcl A u_f(t)||\le C_A(N-p)^{-m}|| u_f(t)||$ due to a boundedness of the integral operators $\mcl T,\,\mcl A$. Additionally, $ ||g(t)-\hat{g}(t)||\le C_3(N-p)^{-m}||g(t)||$. For the commutation error, we have 
\begin{align}
||\widehat{\mcl T\, \dot{u}_f(x,t)}-\mcl T R_{N-p}\dot{u}_f(x,t)||\le \\ ||\widehat{\mcl T\, \dot{u}_f(x,t)}-\mcl T \dot{u}_f(x,t)||+||\mcl T\, \left(\dot{u}_f(x,t)- R_{N-p}\dot{u}_f(x,t)\right)||\le C_4(N-p)^{p-m}||\dot{u}_f(t)||\notag ,
\end{align}
and, analogously, for the  $\left(\widehat{\mcl A\, u_f(x,t)}-\mcl A R_{N-p}(u_f(x,t))\right)$ term.

Since $\mcl T u_f-\widehat{\mcl T \hat{u}_f}=\mcl T \left(u_f- R_{N-p} u_f\right) +\left(\mcl T R_{N-p} u_f- \widehat{\mcl T R_{N-p} u_f}\right)+\widehat{\mcl T e}$, and noting that $e(0)=0$ in the current definition, we obtain the desired estimate (\ref{eq:converg}). 
\end{proof} 
Note that the estimate (\ref{eq:converg}) implies an exponential convergence for smooth solutions. 

\subsection{Temporal treatment}
If $M$ is invertible, Eq.~(\ref{eq:ode}) can be rewritten as
\begin{equation}\label{eq:odetran}
 \frac{d\, \mbf a(t)}{d\,t}=\tilde{A}\,\mbf{a}(t)+\tilde{B}\,\mbf b(t),
\end{equation}
where $\tilde{A}=M^{-1}A$, and $\tilde{B}=M^{-1}$. Invertibility of $M$ generally follows from its block-diagonal structure and well-posedness of the boundary conditions. If $M$ is not invertible, Eq.~(\ref{eq:ode})  would  admit linear in time egiensolutions irrespective of the right-hand side, and this situation will not be considered here. 

We now consider several approaches to the time integration of (\ref{eq:odetran}).

\subsubsection{Exact integration}\label{sec:analtime}

The following lemma establishes an exact solution to the matrix equation (\ref{eq:odetran}).
\begin{lemma}
The solution to the matrix equation (\ref{eq:odetran})
with initial conditions $\mbf{a}(0)=\pmb a_0$
is given by
\begin{equation}\label{timegen}
\mbf a(t)=e^{\tilde{A}\,t}\, \mbf a_0+\int_0^t e^{\tilde{A}(t-s)}\,\tilde{B} \,\mbf b(s)\,ds.
\end{equation}
\end{lemma}
\begin{proof}
It is immediately seen that (\ref{timegen}) satisfies initial conditions at $t=0$. To show that (\ref{timegen}) is a solution to (\ref{eq:odetran}), we differentiate (\ref{timegen}) with respect to time:
\begin{align}
\frac{d\,\mbf a(t)}{d\,t}&=\tilde{A}\,e^{\tilde{A}\,t}\mbf a_0+\int_0^t \tilde{A} \,e^{\tilde{A}(t-s)} \tilde{B} \,\mbf b(s)\,d\,s\\&+e^{\tilde{A}(t-t)}\tilde{B} \,\mbf b(t)=\tilde{A}\,\mbf a(t)+\tilde{B} \,\mbf b(t)\notag,
\end{align}
which satisfies (\ref{eq:odetran}). To show uniqueness, we assume that there exists another solution $\mbf a_1(t)$ that satisfies equation (\ref{eq:odetran}) and initial conditions $\mbf a_1(0)=\mbf a_0$. Denote $\Delta \mbf a(t)=\mbf a_1(t)-\mbf a(t)$. It is easy to verify that $\Delta \mbf a (t)$ satisfies homogeneous equation
\begin{equation}\label{eq:odehomog}
\frac{d\,{\Delta \mbf a} (t)}{d\,t}=\tilde{A}\,\Delta \mbf a(t)
\end{equation}
with homogeneous initial conditions $\Delta \mbf a(0)=0$, from which it immediately follows that $\Delta \mbf a (t)=0$ and $\mbf a_1(t)=\mbf a(t)$.

\end{proof}

Upon substitution $\tilde{A}=M^{-1}A$, and $\tilde{B}=M^{-1}$ into (\ref{timegen}), we recover an exact solution to Equation (\ref{eq:ode}) in our original notation
\begin{equation}\label{timegenorig}
\mbf a(t)=e^{M^{-1}A\,t}\, \mbf a_0+\int_0^t e^{M^{-1}A\,(t-s)}\,M^{-1}\,\mbf b(s)\,ds.
\end{equation}

While a general close-form solution to (\ref{eq:ode}) in the form of (\ref{timegenorig}) exists (provided $M$ is invertible), its analytical evaluation, in practice, is often challenging, since it involves the computation of the matrix exponentials. It can, however, be evaluated easily if the matrix $\tilde{A}=M^{-1}A$ is diagonalizable as $\tilde{A}=S\,\Lambda\,S^{-1}$, in which case the equation (\ref{timegenorig}) simplifies to
\begin{equation}\label{eq:timediag}
\mbf a(t)=S\,e^{\Lambda\,t}S^{-1}\, \mbf a_0+S\,\int_0^t e^{\Lambda(t-s)}  S^{-1} M^{-1} \,\mbf b(s)\,d\,s.
\end{equation}
If, additionally, the inputs $\mbf b(t)$ are such that the integrals 
\begin{equation}\label{intsep}
I_{kl}=\int_0^t e^{\lambda_k(t-s)}b_l(s) d s
\end{equation}
can be evaluated analytically, where $\lambda_k$,   $b_l(t)$ for  $\{\,k,\,l\}=\{1\ldots N_d\}$, are the eigenvalues of $\tilde{A}$ and components of the vector $\mbf b(t)$, respectively, the entire vector-valued integral $\mbf I=\int_0^t e^{\Lambda(t-s)}  S^{-1} M^{-1} \,\mbf b(s)\,d\,s$ in (\ref{eq:timediag}) can be evaluated componentwise as $I_k=\sum_{l=1}^{N_d}\,I_{kl} \{S^{-1} M^{-1} \}_{kl}$, where $I_k$ is the $k^{th}$ component of $\mbf I$, $\{S^{-1} M^{-1}\}_{kl}$ is the corresponding entry of the matrix $S^{-1}  M^{-1}$ in the $k^{th}$ row and $l^{th}$ column,  and summation over $k$ is not implied. Furthermore, if inputs $\mbf b(t)$ are separable into $m$ time-dependent entries $\mbf b(t)=\sum_{l=1}^m \boldsymbol{\alpha}_l b_l(t)$, $m<N_d$, $\boldsymbol{\alpha}_l$ are the vectors independent of time, the evaluation of the integral in (\ref{eq:timediag}) reduces to a computation of $m\,N_d$ integrals of the form (\ref{intsep}), and the reconstruction process yields $\int_0^t e^{\Lambda(t-s)}  S^{-1} M^{-1} \,\mbf b(s)\,d\,s= \sum_{l=1}^m D_l\, S^{-1} M^{-1}\boldsymbol{\alpha}_l$, where $D_l$ is a diagonal matrix that, for each $l$, consists of the corresponding $I_{kl}$ values, such that $D_l=\diag(I_{kl})$.

\subsubsection{Alternative exact integration}

While Equation (\ref{eq:timediag}) and its analytical evaluation via an approach described above provides a robust solution whenever $M$ is invertible, the ODE system (\ref{eq:odetran}) is stable, and matrix $\tilde{A}=M^{-1}A$ is diagonalizable, in some cases, we can further reduce the errors associated with the inversion of the matrix $M$
 by employing the following alternative form of the solution to  (\ref{eq:ode}) given by the following lemma.


\begin{lemma}
It matrix $M$ is diagonalizable as $M=S\,\Lambda\,S^{-1}$, and it does not have any zero eigenvalues, a solution to the equation (\ref{eq:ode}) with initial conditions $\mbf{a}(0)=\pmb a_0$
is given by
\begin{equation}\label{timediag}
\mbf a(t)=S\,e^{\,\Lambda^{-1}\,S^{-1}\,A\,S\,t}\, S^{-1}\mbf a_0+S\int_0^t e^{\Lambda^{-1}\,S^{-1}\,A\,S\,(t-s)} \Lambda^{-1}\,S^{-1} \mbf b(s)\,d\,s.
\end{equation}
\end{lemma}
\begin{proof}
Upon substituting $M=S\,\Lambda\,S^{-1}$ into equation (\ref{eq:ode}), multiplying both sides of it by  $S^{-1}$, and defining $\mbf z=S^{-1}\mbf a$, equation (\ref{eq:ode}) reduces to
\begin{equation}\label{eq:odediag}
 \Lambda \frac{d\,{\mbf z}(t)}{d\,t}=S^{-1}A\,S\,\mbf{z}(t)+S^{-1}\mbf{b}(t).
\end{equation}
Upon multiplying equation (\ref{eq:odediag}) by the inverse of $\Lambda$, the solution given by (\ref{timediag}) follows immediately from (\ref{timegen}) and substitution $\mbf a=S\,\mbf z$ .  
\end{proof}

Note that for the PDEs with constant coefficients, $A$ would be a multiple of an identity matrix, so that $\Lambda^{-1}\,S^{-1}\,A\,S$ is by itself diagonal. Alternatively, its diagonalization similar to a procedure described in \cref{sec:analtime} needs to be perrormed for an analytical evaluation of (\ref{timediag}).

Unfortunately, the eigenvalues of $M^{-1} A$ are different from the eigenvalues of $\Lambda^{-1}\,S^{-1}\,A\,S$, which can render the evaluation of the integral in (\ref{timediag}) unstable, especially if the eigenvalues of $M^{-1} A$ are purely imaginary, as in the hyperbolic problems. This approach, therefore, can not be advocated as a general-purpose solution. However, for diffusive problems, it significantly reduces approximation errors associated with the evaluation of the matrix exponentials in (\ref{eq:timediag}). Since the purpose is to demonstrate strong spatial convergence properties of the PIE-Galerkin approximation decoupled from the temporal errors, we intend to use (\ref{timediag}) whenever possible.
\subsubsection{Gauss integration}

An analytical integration procedure described above would fail if
\begin{itemize}
\item Inhomogeneous inputs $\mbf b(t)$ have a functional form that does not allow for an analytical evaluation of the integral in (\ref{eq:timediag}) or (\ref{timediag}).
\item Either $M$ is not diagonalizable, or eigenvalues of  $\Lambda^{-1}\,S^{-1}\,A\,S$ are such that evaluation of (\ref{timediag}) is unstable.
\item $\tilde{A}=M^{-1}A$ is not diagonalizable, so that (\ref{timegenorig}) can not be reduced to (\ref{eq:timediag}).
\end{itemize}
In this case, the integral in (\ref{timegenorig}) can be approximated numerically. In this work, the total time interval is partitioned into $N_{int}$ sub-intervals, and a Gauss-Lobatto quadrature of a specified order $Ng$ is used for each time interval. This approach alleviates the problems associated with the analytical integration mentioned above, and also avoids some difficulties attributed to the classical time stepping procedures. First, it does not suffer from the CFL-type instabilities and the associated time step restrictions of the classical time stepping schemes.  As long as the ODE system is physically stable (that is, it does not possess any eigenvalues with positive real parts), the Gauss integration approach will succeed. Second, it does not require a sequential approach and can, in principle, be evaluated parallelly in time, since the integral at each segment can be independently evaluated and subsequently added to form a final solution. However, there are also some drawbacks associated with this approach. Since it requires evaluation of the matrix exponentials, this becomes sensitive to the value of time step. Since the power of the exponential term in (\ref{timegenorig}) is proportional to $(t-s)$, discretization close to the end of the time period $t$ is especially important. It was found that clustering of the time intervals towards the end of the time period $t$, so that the time discretization is finer as the values of $s$ approach the final time, is helpful for some problems. In these cases, a geometric progression was used to determine the value of the time intervals with a specified ratio $r$. Within each time interval, the Gauss-Lobatto (GL) integration with the nodes specified by GL quadrature are used.

\subsubsection{Backward differentiation formula}

While the above approaches associated with the approximation of the exact solution in the form  (\ref{timegenorig}) typically provide the lowest errors in the current one-dimensional situation, its applicability to mutliple dimensions and to larger matrices might be limited. To compare the two approaches to the classical time stepping procedures and to show the effect of the temporal discretization errors on the spatial convergence, we also implement a Backward differentiation formula (BDF) for the time integration. Backward differentiation formula (BDF) is an implicit time integration scheme, which, as applied to (\ref{eq:odetran}) is given by
\begin{equation}
\sum_{p=0}^k  \beta_p \:\mbf a^{n-p}=\Delta\,t \,(\tilde{A}\,\mbf{a}^n+\tilde{\mbf b}^n),
\end{equation}
where $k$ is the order of accuracy of the scheme, $\Delta t$ is the time step, and the vectors with the superscript $n$ correspond to their value at the discrete time level $t^n$. BDF schemes of the order 3 and 4 are considered here, and the corresponding BDF coefficients $\beta_p$ for these two schemes are provided in Table~\ref{table:bdf}.
\begin{table}
\caption{Coefficients $\beta_{p}$ of the BDF$k$ scheme, $k=3,4$.}
\begin{center} \footnotesize
\begin{tabular}{|c|c|c|c|c|c|} \hline
 $p$  & 0 & 1 & 2&3 & 4\\ 
   \hline
$k=3$& 11/6 & -3 & 3/2 &-1/3 & \\ \hline
 $k=4$& 25/12 & -4  & 3 & -4/3  & 1/4 \\ \hline
\end{tabular}
\end{center}
\label{table:bdf}
\end{table}
Since BDF3/BDF4 schemes can be used only starting from the $3^{rd}/4^{th}$ time steps respectively, to get a nominal temporal order of convergence, we initialize the required number of initial time steps with the exact solution in the subsequent examples. In practice, where exact solution is not available, lower order BDF schemes could be used for the initial time steps. 

\subsection{Software}

The computational methods described above were implemented within a general-purpose open-source PDE solver PIESIM available for download at \verb+http://control.asu.edu/pietools+. PIESIM, which is based on a MATLAB package, is fully integrated with PIETOOLS~\cite{shivakumar2020pietools}, an open-source software previously developed by the authors for construction, manipulation and optimization of the PI operators. For the purposes of the presented methodology, PIETOOLS handles the conversion of a given PDE problem into a PIE framework and constructs the corresponding 3-PI operators, while PIESIM computes a numerical solution of the PIE problem using the PIE-Galerkin methodology, and transforms the PIE solution back to represent a required solution of the original PDE problem. All numerical examples described below were solved using PIESIM.

\section{Numerical Examples}\label{sec:example}
This section demonstrates the application of the presented numerical methodology to several canonical PDE equation problems. 

\subsection{Parabolic Problems}

\subsubsection{Example 1: Diffusion Equation with constant viscosity}

We begin with the consideration of a classical diffusion equation, given by
\begin{equation}\label{eq:heat}
u_t=\nu\,u_{xx},
\end{equation}
whit $\nu$ a scalar,
defined on a domain $x\in[-1,1]$. In lieu of a standardized representation given in~\cref{sec:standard}, Eq.~(\ref{eq:heat}) corresponds to $A_0(x)=A_1(x)=0$, $A_2(x)=\nu$, $n_0=n_1=0$, $n_2=1$, $u_2(x,t)=u(x,t)$ is a primary state, while from (\ref{eq:fvector}), $u_{f2}(x,t)=u_{xx}(x,t)$ is a fundamental state.  
\paragraph{Dirichlet-Dirichlet boundary conditions}
We first consider Dirichlet - Dirichlet boundary conditions, defined as $u(-1,t)=h_1(t), u(1,t)=h_2(t)$, with
the boundary conditions matrix 
\begin{equation}
B=\bmat{1 & 0 & 0 & 0 \\0 &1 &0 &0}\vspace{-1mm}.
\end{equation}
With this value of $B$, the 3-PI operators $\mcl T$ and $\mcl A$ in~(\ref{eq:3-pi}) for the equation (\ref{eq:heat}) are parameterized by $G_0=0, G_1(x,s)=x-s, G_2(x,s)=\frac{1}{2}(x+1)(s-1)$, and $H_0(x)=\nu, \,H_1=H_2=0$, respectively.  Applying the discretization procedure described in~\cref{sec:numsol}, we obtain a discrete $N_d\times N_d$ matrix $M$, which, given that $n_0=n_1=0$, $n_2=1$, reduces to a $N-1\times N-1$ matrix, which, for example, for $N=7$ is computed as
\begin{equation}\label{matdir}
M=\bmat{-1/4 & 0 & 7/48 & 0 & -1/60 & 0  \\0 &-1/24 &0 &1/20 & 0 & -1/168\\ 1/4 & 0 & -1/6 & 0 & 1/24 & 0 \\  0 & 1/24 & 0 & -1/16 & 0 & 1/48  \\ 0 & 0 & 1/48 & 0 & -1/30 & 0  \\0 & 0 & 0 & 1/80 & 0 &-1/48 }\vspace{-1mm},
\end{equation}
which is a pentadiagonal matrix with the exception of the first two rows, in accordance with~\cref{lem:chebmat},
while the matrix $A=\nu \cdot I$. 
We now specify the following values for the boundary and initial conditions: \\$u(-1,t)=\sin(-9\pi/8)e^{-\nu\,\pi^2 t}, \,u(1,t)=\sin(11\pi/8)e^{-\nu\,\pi^2 t}$, $u(x,0)=\sin(5\pi/4\, x+\pi/8)$, and initial conditions on the fundamental state,\\ $u_{f2}(x,0)=u_{xx}(x,0)=-(5\pi/4)^2\sin(5\pi/4\, x+\pi/8)$, with the exact solution 
$u(x,t)=\sin(5\pi/4 \,x+\pi/8)e^{-\nu\,\pi^2 t}$.

 In this case, inhomogeneous term in the form $-K(x) B_T^{-1} d\,\mbf{h}(t)/d t$ is present in~(\ref{eq:g}), with 
\begin{equation}\label{bheat}
K(x)=\bmat{1 & x+1}; \:\:B_T^{-1}=\bmat{1 & 0 \\-1/2 &1/2}, 
\end{equation}
so that $K(x)B_T^{-1}=\bmat{\frac{1}{2}(-x+1) & \frac{1}{2}(x+1)}$.

\begin{figure}
\begin{subfigmatrix}{3}
 \setlength{\unitlength}{0.012500in}  \vspace{-1mm}
  \subfigure[Solution plot at $N=8$ for Dirichlet-Dirichlet case. Black solid line, exact solution; symbols, numerical solution (see the caption).]{\includegraphics[width=0.32\textwidth]{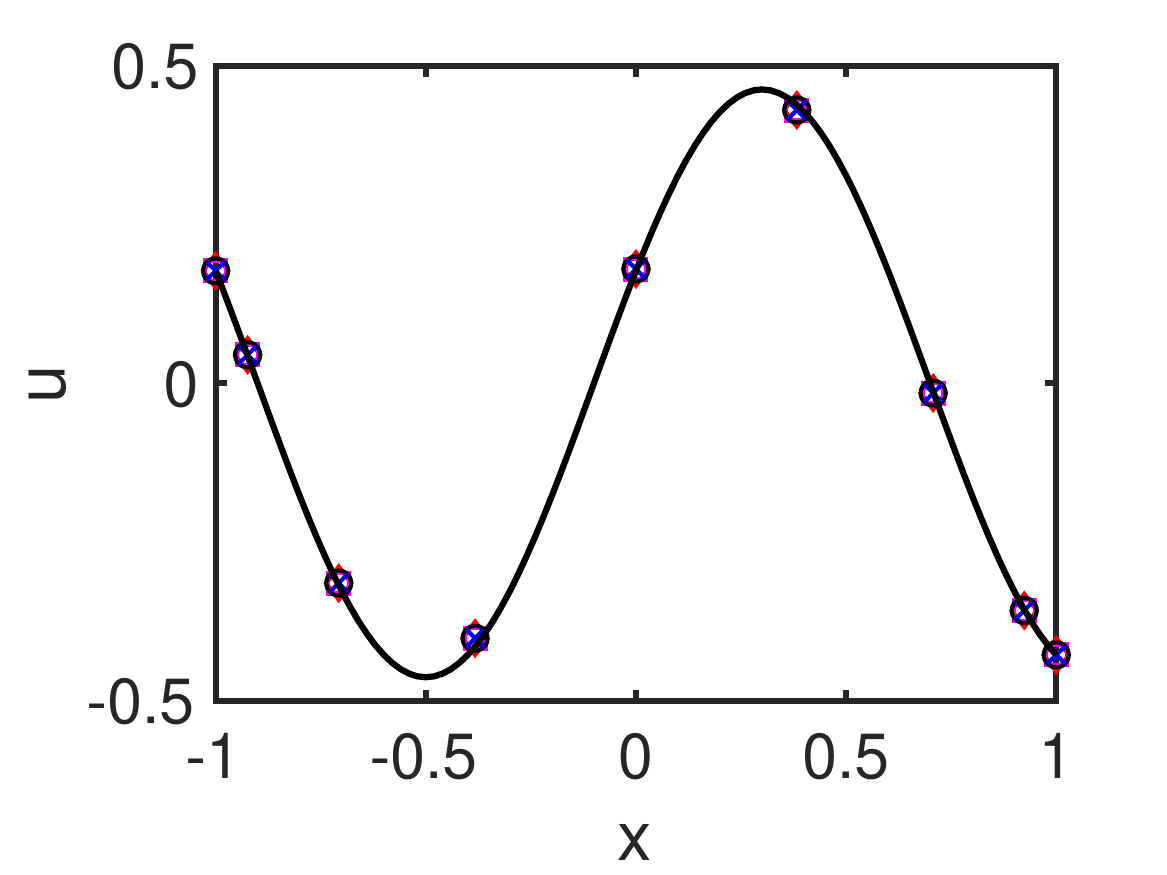}\label{fig:heatconds}}
  \subfigure[$L_2$ error versus the polynomial order $N$ for Dirichlet-Dirichlet case.]{\includegraphics[width=0.32\textwidth]{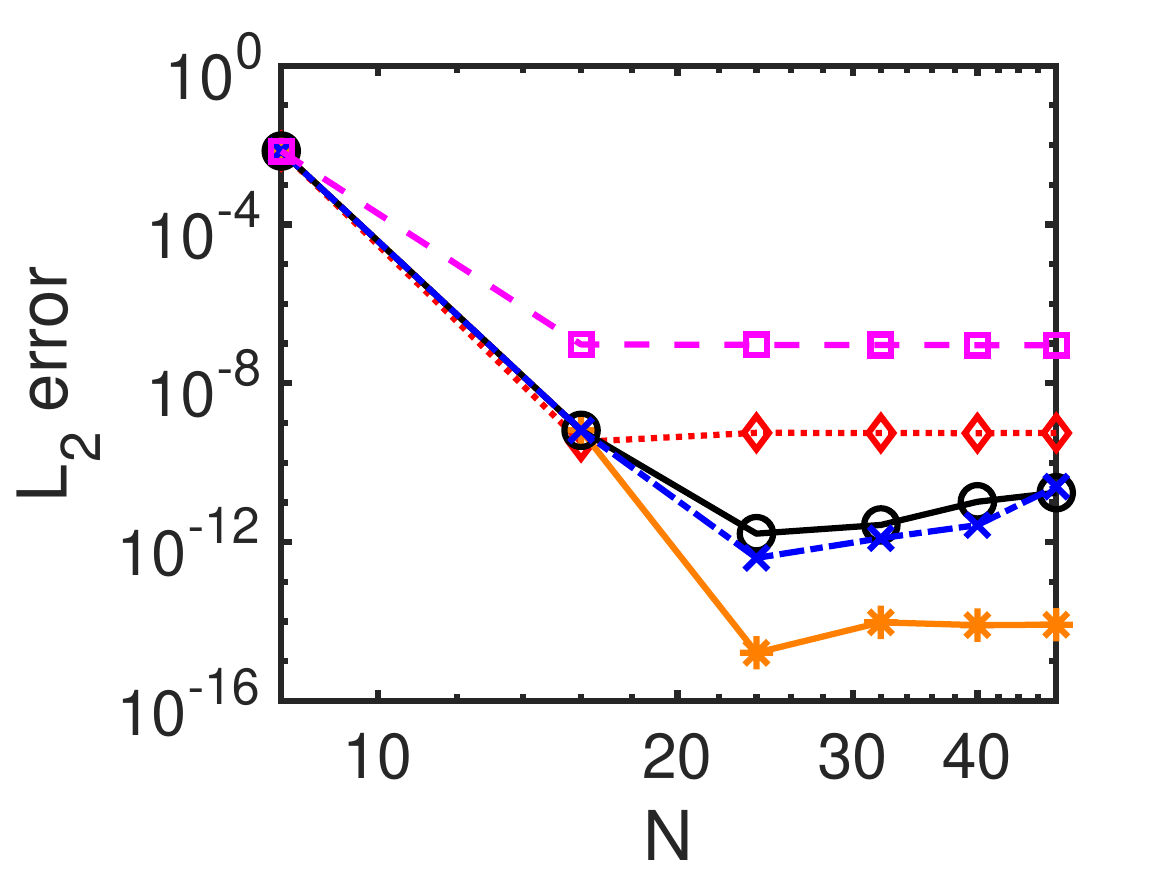}\label{fig:heatconde}}
  \subfigure[$L_2$ error versus the polynomial order $N$ for Dirichlet-Neumann case.]{\includegraphics[width=0.32\textwidth]{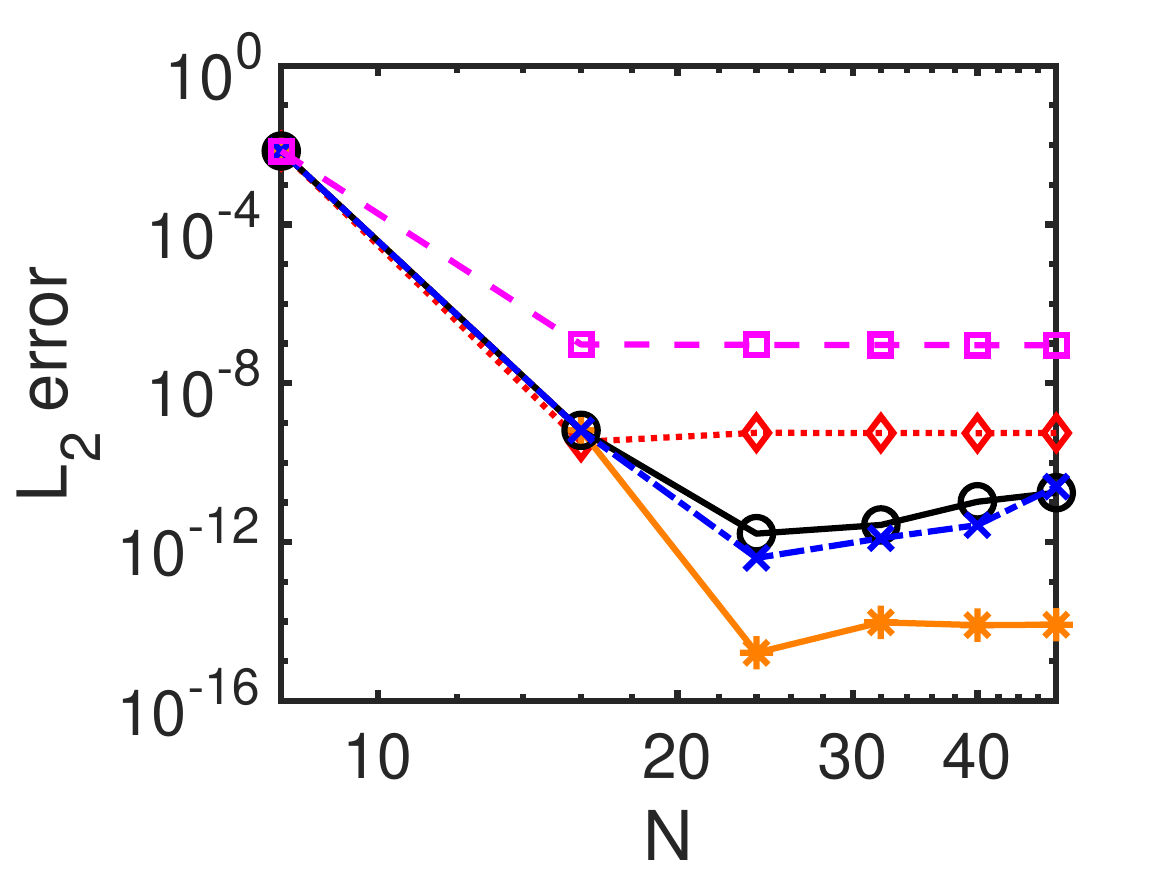}\label{fig:heatconne}}
     \end{subfigmatrix}
     \caption{Solution and convergence plots for Example 1: diffusion equation with a constant viscosity at a time $t=0.1$. Orange solid line with asterisks, analytical evaluation of Eq. (\ref{timediag}); black solid line with circles, analytical evaluation of Eq. (\ref{eq:timediag}); blue dash-dotted line with crosses, Gauss integration of Eq. (\ref{timegenorig}) with $N_g=10$ and $N_{int}=10$ non-uniform time intervals with the geometric progression ratio $r=0.25$; red dotted line with diamonds, BDF4 with $\Delta\,t=10^{-3}$; magenta dashed line with squares, BDF3 with $\Delta\,t=10^{-3}$.\label{fig:heatcon} }
      \end{figure}

The solution and the convergence plots for this test case with different time integration approaches are presented for $\nu=0.5$, time step $\Delta t=10^{-3}$, and $t=0.1$ in \cref{fig:heatconds}, \cref{fig:heatconde}.

 \paragraph{Dirichlet-Neumann boundary conditions} We now consider the case of Dirichlet-Neumann boundary conditions. We use the initial conditions and the analytical solution of the previous example, but we change the boundary condition at the right end to be of Neumann type, i.e. boundary conditions are defined as follows: \\   $u(-1,t)=\sin(-9\pi/8)e^{-\nu\,\pi^2 t}, \,u_x(1,t)=5\pi/4\cos(11\pi/8)e^{-\nu\,\pi^2 t}$. The boundary conditions matrix is now
\begin{equation}
B=\bmat{1 & 0 & 0 & 0 \\0 &0 &0 &1}\vspace{-1mm},
\end{equation}
which changes the structure of the 3-PI $\mcl T$ operator, which is now given by $G_0=0, G_1(x,s)=x-s, G_2(x,s)=-x-1$. The operator $\mcl A$ is still the same as in the previous example,
 so as $K(x)$. However, due to a different matrix $B$, we a have a different matrix $B_T^{-1}=\bmat{1 & 0\\ 0 & 1}$, and a different operator $K(x)B_T^{-1}=\bmat{1 & x+1}$. As an example, a discrete matrix $M$ for this test case for $N=7$ is given below:
\begin{equation}\label{matneum}
M=\bmat{-5/4 & -1/3 & 23/48 & 1/5 & 1/20 & 1/21  \\-1 &-3/8 &1/3 &1/4 & 1/15 & 1/24\\ 1/4 & 0 & -1/6 & 0 & 1/24 & 0 \\  0 & 1/24 & 0 & -1/16 & 0 & 1/48  \\ 0 & 0 & 1/48 & 0 & -1/30 & 0  \\0 & 0 & 0 & 1/80 & 0 &-1/48 }\vspace{-1mm}.
\end{equation}
As proven in the \cref{lem:chebmat}, the boundary conditions affect only the first two rows of the matrix $M$, while the rest of the matrix is unchanged between (\ref{matdir}) and (\ref{matneum}).
The solution and the convergence plots for this test case are presented for $\nu=0.5$, time step $\Delta t=10^{-3}$, and $t=0.1$ in \cref{fig:heatconne}.

\subsubsection{Example 2: Diffusion Equation with variable viscosity}
We now consider a diffusion equation with a variable viscosity, such that
\begin{equation}\label{eq:heatvar}
u_t=x\,u_{xx}.
\end{equation}
We use the domain $x\in[0,2]$ to ensure a non-negative value of viscosity for physical stability. We define
initial conditions as $u(x,0)=-x^2$, boundary conditions as Dirichlet-Dirichlet with $u(0,t)=0$, $u(2,t)=-4\,t-4$. In this case, an analytical solution exists, which is given by $u(x,t)=-2\,x\,t-x^2$. When the physical domain does not coincide with $[-1,1]$, a mapping of the physical domain $x\in[a,b]$ into the computational domain $x^{(c)}\in[-1,1]$ must be performed according to (\ref{eq:mapxi}), (\ref{eq:mapfromxi}), with a corresponding transformation applied to Eq.~(\ref{eq:piein}),
\begin{equation}\label{eq:pieshort}
\mcl T^{(c)}\: \frac{\partial \,\mbf u^{(c)}_f(x^{(c)},t)}{\partial \,t}=\mcl A^{(c)} \: \mbf u^{(c)}_f(x^{(c)},t)+\mbf g^{(c)} (x^{(c)},t),
\end{equation}
where the superscript $(c)$ indicates that all the space-dependent functions and the 3-PI operators are now evaluated in the computational domain.  
Since the state of the PDE and the boundary conditions matrix $B$ are the same as in Example 1 with Dirichlet-Dirichlet boundary conditions, the 3-PI operator $\mcl T^{(c)}$, $K^{(c)}x^{(c)}$ and  $K^{(c)}x^{(c)}B_T^{-1}$, are, again, the same, when evaluated in the computational domain. However, for the operator  $\mcl A^{(c)}$ we have $H^{(c)}_0(x^{(c)})=x^{(c)}+1, \,H^{(c)}_1=H^{(c)}_2=0$.

\begin{figure}
\begin{subfigmatrix}{2}
 \setlength{\unitlength}{0.012500in}  \vspace{-1mm}
   \subfigure[Solution plot at $N=8$. Solid line, exact solution; symbols, numerical solution.]{\includegraphics[width=0.49\textwidth]{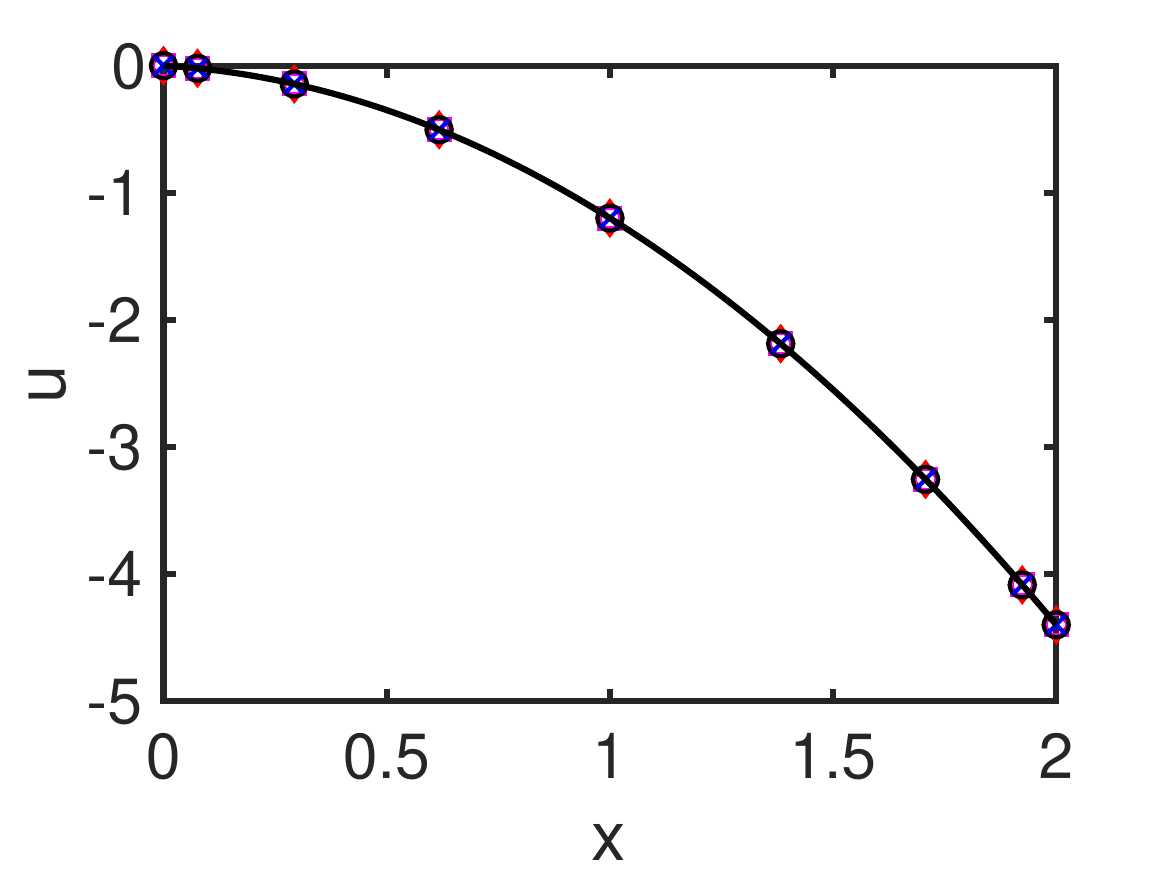}\label{fig:heatvars}}
  \subfigure[$L_2$ error versus the polynomial order $N$.]{\includegraphics[width=0.49\textwidth]{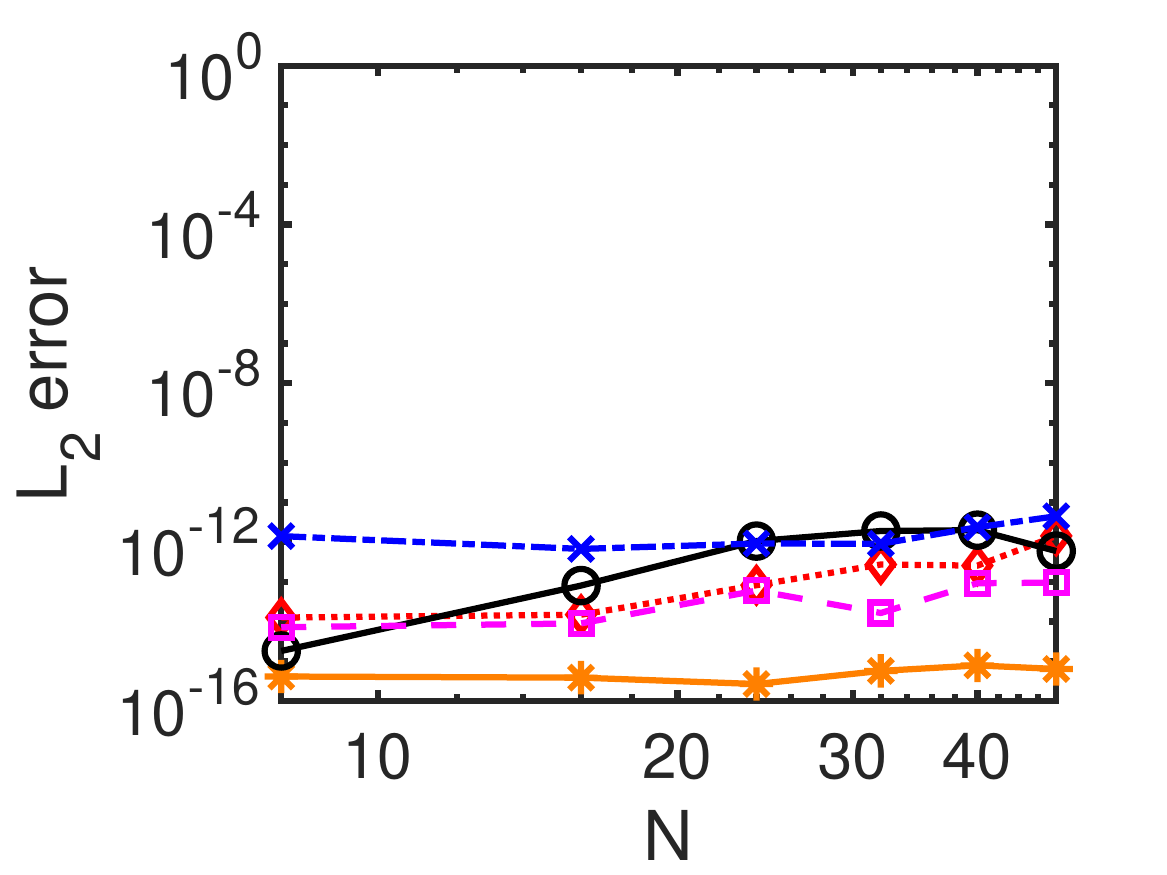}\label{fig:heatvare}}
      \end{subfigmatrix}
     \caption{Solution and convergence plots for Example 2: diffusion equation  with a variable viscosity at a time $t=0.1$. Lines and symbols are the same as in \cref{fig:heatcon}.\label{fig:heatconvar} }
      \end{figure}

The solution and the convergence plots for this test case are presented for $t=0.1$ in \cref{fig:heatconvar}. 
Note, since the exact solution is a second-order polynomial, which is resolved starting with $N=2$, the initial error is already at a machine precision in this test case.

 \subsubsection{Example 3: Convection-Diffusion Equation}
 
\begin{figure}
\begin{subfigmatrix}{2}
 \setlength{\unitlength}{0.012500in}  \vspace{-1mm}
  \subfigure[Solution plot at $N=8$. Solid line, exact solution; symbols, numerical solution.]{\includegraphics[width=0.49\textwidth]{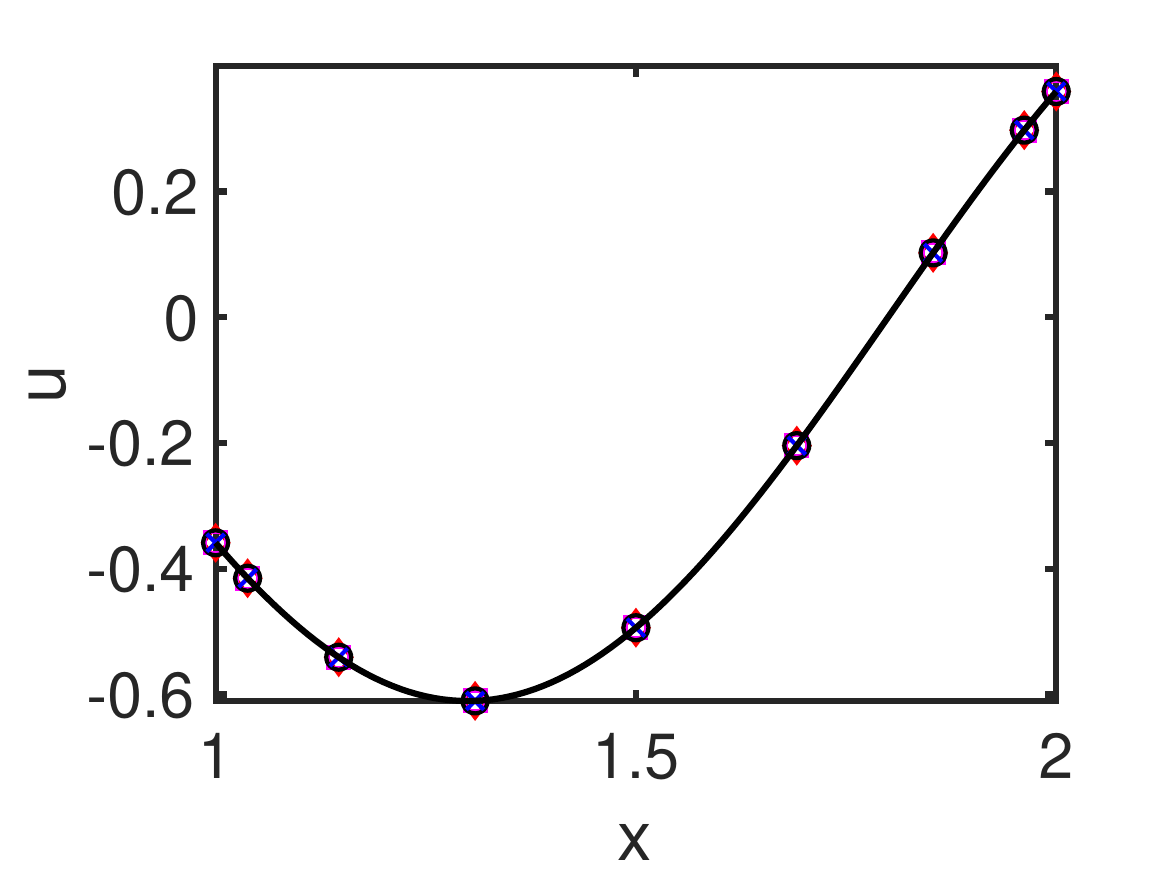}}
  \subfigure[$L_2$ error versus the polynomial order $N$.]{\includegraphics[width=0.49\textwidth]{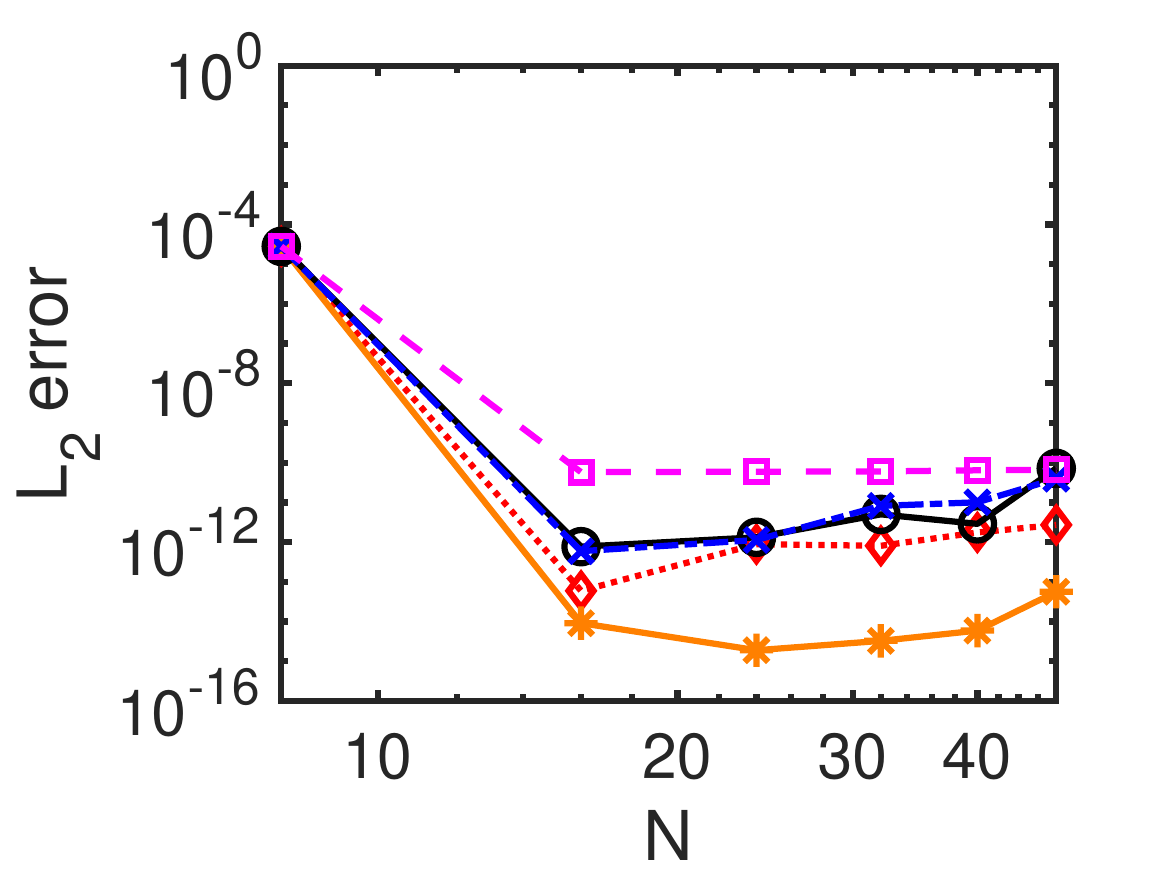}}
      \end{subfigmatrix}
     \caption{Solution and convergence plots for Example 3: convection-diffusion equation for $\nu=0.5, c=-2$  at a time $t=0.1$. Lines and symbols are the same as in \cref{fig:heatcon}. \label{fig:convdiff} }
      \end{figure}

 In this example, we consider the convection-diffusion equation  given by 
\begin{equation}\label{eq:convdiff}
u_t+c\,u_x=\nu\,u_{xx},
\end{equation}
defined on a domain $x\in[a,b]$, which has an exact solution \\ $u(x,t)=\sin\left(\pi(x-c\,t)\right)e^{-\nu\pi^2\,t}$ that satisfies the initial condition $u(x,0)=\sin(\pi x)$, and Dirichlet-Dirichlet boundary conditions $u(a,t)=\sin\left(\pi(a-c\,t)\right)e^{-\nu\pi^2\,t}, \:u(b,t)=\sin\left(\pi(b-c\,t)\right)e^{-\nu\pi^2\,t}$. As in the previous example, a coordinate transformation must be done in accordance with (\ref{eq:mapxi}), (\ref{eq:mapfromxi}), (\ref{eq:pieshort}) to represent the PIE equation in the computational domain $x^{(c)}\in[-1,1]$.
In this case, we have $A_0^{(c)}(x^{(c)})=0, \,A_1^{(c)}(x^{(c)})=-2c/(b-a), \,A_2^{(c)}(x^{(c)})=4\,\nu/(b-a)^2$. The fundamental state, the matrix $B$, the 3-PI operator $\mcl T^{(c)}$ as well as the matrix $K^{(c)}(x^{(c)})$ are, again, the same as in the previous example, but now we also have a non-zero $A_1^{(c)}(x^{(c)})$, which leads to a non-zero value of the operators $H_1^{(c)}(x^{(c)}, s^{(c)})=A_1^{(c)}(x^{(c)}), H_2^{(c)}(x^{(c)},s^{(c)})=\frac{1}{2} A^{(c)}_1(x^{(c)})(s^{(c)}-1)$, while $H_0^{(c)}(x^{(c)}, s^{(c)})=A_2^{(c)}(x^{(c)})$. Note that, in this case, we have a contribution to a non-homogeneous term due to a second term in equation (\ref{eq:g}). However, since $n_1=0, n_2=1$, and $B_T^{-1}$ given by Eq.~(\ref{bheat}) in this problem has antisymmetric entries in the second row, this term would vanish for a solution with the equal values of the boundary condition entries, therefore we consider a non-periodic domain given by $[a,b]=[1,2]$. The solution and the convergence plots are presented  in \cref{fig:convdiff} for this test case for  $\nu=0.5, c=-2$  at a time $t=0.1$. 

 \subsubsection{Example 4: Parabolic Equation with Forcing}
\begin{figure}
\begin{subfigmatrix}{2}
 \setlength{\unitlength}{0.012500in}  \vspace{-1mm}
  \subfigure[Solution plot at $N=8$.  Black solid line, exact solution; symbols, numerical solution.]{\includegraphics[width=0.49\textwidth]{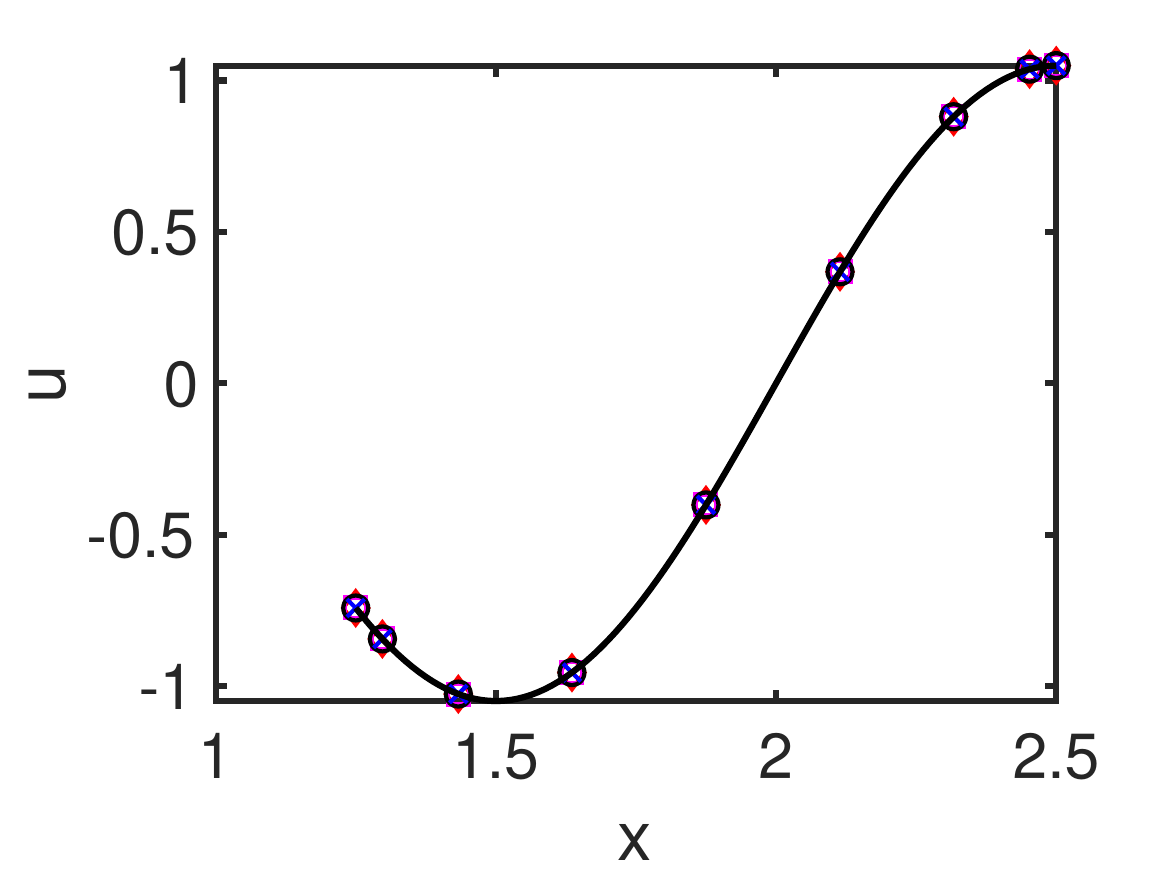}}
  \subfigure[$L_2$ error versus the polynomial order $N$.]{\includegraphics[width=0.49\textwidth]{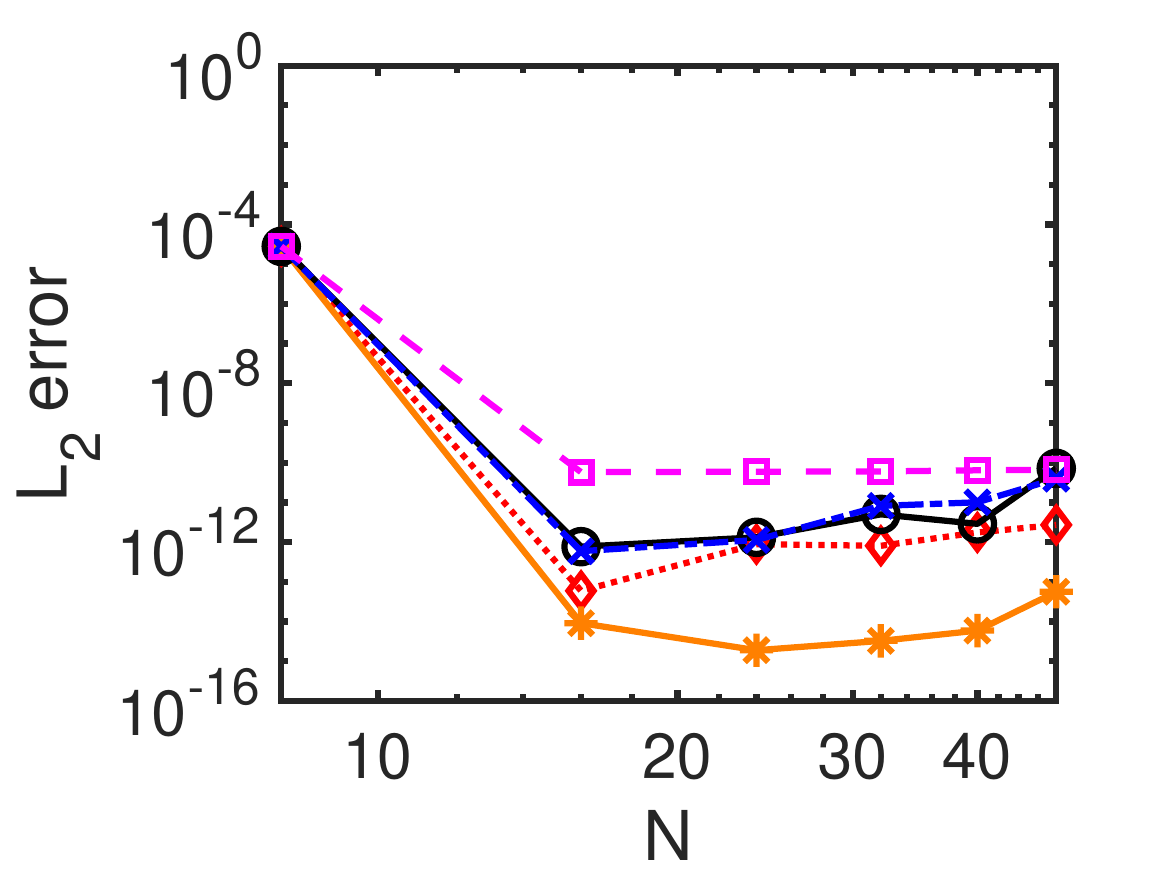}}
      \end{subfigmatrix}
     \caption{Solution and convergence plots for Example 4: parabolic equation with forcing at a time $t=0.1$. Lines and symbols are the same as in \cref{fig:heatcon}. \label{fig:parfull} }
      \end{figure}

In this example, we test a full form in the PDE representation (\ref{eq:primary}), where all three coefficients $A_0(x)$, $A_1(x)$ and $A_2(x)$ are present. We use the Method of Manufactured Solutions~\cite{oberkampf1998issues, roache2002code}  to construct an exact solution of the equation
\begin{equation}\label{eq:force}
u_t=\alpha \,u+\beta\, u_x+\gamma \,u_{xx}+f(x,t)
\end{equation}
in the form $u(x,t)=\sqrt{t+1}\,\sin(\pi\,x)$, with the corresponding right-hand side \\$f(x,t)=\frac{1}{2\sqrt{t+1}}\sin(\pi\,x) - \alpha\,\sqrt{t+1} \sin(\pi\,x)-\beta\,\pi\,\sqrt{t+1}\cos(\pi\,x)+\gamma\,\pi^2\,\sqrt{t+1}\,\sin(\pi\,x)$. In this case, $A_0(x)=\alpha$, $A_1(x)=\beta$, $A_2(x)=\gamma$, $n_0=n_1=0$, $n_2=1$. We apply the Neumann boundary condition $u_x(a,t)=-\pi\,\sqrt{t+1}\,\cos(\pi\,a)$ on the left side, and the Dirichlet boundary condition $u(b,t)=\sqrt{t+1}\,\sin(\pi\,b)$ on the right side, for which the matrix $B$ reads 
\begin{equation}
B=\bmat{0 & 0 & 1 & 0 \\0 &1 &0 &0}\vspace{-1mm}.
\end{equation}
Upon the transformation of the PDE (\ref{eq:force}) into the computational domain $x^{(c)}\in[-1,1]$, the corresponding functions are transformed as  $A^{(c)}_0(x^{(c)})=A_0(x)=\alpha$, $A^{(c)}_1(x^{(c)})=2A_1(x)/(b-a)=2\beta/(b-a)$ and $A^{(c)}_2(x^{(c)})=4A_2(x)/(b-a)^2=4\gamma/(b-a)^2$. Since the forcing function does not contain any derivatives of $u(x,t)$,no transformation of the forcing function is required. Finally, the Dirichlet boundary conditions are imposed on $u^{(c)}(x^{(c)})$ as given by $u^{(c)}(1,t)=u(b,t)=\sqrt{t+1}\,\sin(\pi\,b)$, while Neumann boundary conditions are recalculated as $u^{(c)}_{x^{(c)}}(-1,t)=u_x(a,t)\cdot (b-a)/2=-\pi\,\sqrt{t+1}\,\cos(\pi\,a)\cdot (b-a)/2$. The operators in  the computational domain thus become $G^{(c)}_0=0, G^{(c)}_1(x^{(c)},s^{(c)})=x^{(c)}-s^{(c)}, G^{(c)}_2(x^{(c)},s^{(c)})=s^{(c)}-1$, $H^{(c)}_0(x^{(c)},s^{(c)})=A^{(c)}_2(x^{(c)}), H^{(c)}_1(x^{(c)},s^{(c)})=A^{(c)}_1(x^{(c)})+A^{(c)}_0(x^{(c)})(x^{(c)}-s^{(c)}), H^{(c)}_2(x^{(c)},s^{(c)})=\\A^{(c)}_0(x^{(c)})(s^{(c)} - 1)$, $K^{(c)}(x^{(c)})=\bmat{x^{(c)}-1 & 1}$. In this case, all four components in the inhomogeneous term (\ref{eq:g}) are present. We use the following parameter values for this test case: $a=1.25, b=2.5, \alpha=4, \beta=2, \gamma=0.5$. The solution and the convergence plots are presented  in \cref{fig:parfull} for this test case at a time $t=0.1$. 

 \subsubsection{Example 5: Euler-Bernoulli Beam}
 
 Euler-Bernoulli beam model is represented by a fourth-order PDE 
 \begin{equation}
 u_{tt}=-c\,u_{xxxx},\label{eq:beam}
 \end{equation}
 where $c=EI/\mu$, $E$ is the elastic modulus, $I$ is the second moment of area of the beam's cross-section, and $\mu$ is the mass per unit length.  In a cantilevered state described by the boundary conditions 
 \begin{equation}\label{beambc}
 u(0,t)=u_x(0,t)=u_{xx}(L,t)=u_{xxx}(L,t)=0
 \end{equation}
  a free vibration solution exists given by the following harmonic modes~\cite{volterra1965dynamics}
\begin{equation} 
 u_n(x,t)={\text{Re}}\left[{\tilde{u}}_n(x)~e^{{-i\omega_n t}}\right],
 \end{equation}
 with eigenmodes
 \begin{equation}\label{eq:beamexact}
 \tilde{u}_n(x)=A_n[\cosh (\beta_n x) -\cos (\beta_n\,x)+\frac{\cos (\beta_n x) +\cosh (\beta_n\,x)}{\sin (\beta_n x) +\sinh (\beta_n\,x)}(\sin (\beta_n x) -\sinh (\beta_n\,x))],
 \end{equation}
eigenvalues $\beta_n$ being a solution of the following eigenvalue problem
 \begin{equation}\label{beameig}
 \cosh (\beta_n\,L) \cos (\beta_n\,L)+1=0,
 \end{equation}
and the vibration frequencies defined as $\omega_n=\beta_n^2\sqrt{EI/\mu}=\beta_n^2\sqrt{c}$.

To cast Equation (\ref{eq:beam}) into a state-space representation of (\ref{eq:primary}), we define the following states $v_1(x,t)=u_t(x,t)$, $v_2(x,t)=u_{xx}(x,t)$, so that (\ref{eq:beam}) transforms into 
\begin{equation}
\mbf v_t=\underbrace{\bmat{0&-c\\1&0}}_{A_2}\mbf v_{xx},\label{eq:beamstate}
\end{equation}
where the state vector $\mbf v=[v_1\:\:\: v_2]^T$, $n_0=n_1=0, n_2=2$, which represents an example of a vector-valued state. Thus, the fundamental state is $\mbf v_f=[v_{1xx}\:\:\: v_{2xx}]^T$, $A_0=A_1=0$, and $A_2$ is as given by Eq. (\ref{eq:beamstate}). For the boundary conditions defined by (\ref{beambc}), the last two equations can be restated in terms of the state $v_2(x,t)$ as $v_2(L,t)=0, v_{2x}(L,t)=0$. The first two boundary conditions can be differentiated in time to give boundary constraints for the state $v_1(x,t)$ as $v_1(0,t)=0, v_{1x}(L,t)=0$. With these, the boundary conditions matrix $B$ reads
\begin{equation}
\underbrace{\bmat{1&0&0&0&0&0&0&0\\
0&0&0&1&0&0&0&0\\
0&0&0&0&1&0&0&0\\
0&0&0&0&0&0&0&1}}_{B}{\scriptsize\bmat{v_1(0,t)\\v_2(0,t)\\v_1(L,t)\\v_2(L,t)\\v_{1x}(0,t)\\v_{2x}(0,t)\\v_{1x}(L,t)\\v_{2x}(L,t)}}=0.
\end{equation}
To reconstruct the original variable $u(x,t)$ from the state-space variables $v_1(x,t)$, $v_2(x,t)$, we can utilize Equation (\ref{tran:uxx}) to recover $u(x,t)$ from its second-derivative $u_{xx}(x,t)=v_2(x,t)$. In the PIE framework, this effectively can be done by a transformation (\ref{eq:maph}) applied to $v_2(x,t)$, with $\mcl T=\{0,x-s,0\}$, $K(x)B_T^{-1}=\bmat{1 & x-a}$, and $\mbf h(t)=\bmat{u(a,t)&u_x(a,t)}^T$, with $a=0$, which, incidentally is equivalent to a PIE transformation with $n_0=n_1=0, n_2=1$, $A_0=A_1=0, A_2=1$, and the boundary conditions given by $u(a,t)=h_1(t),  u_x(a,t)=h_2(t)$. This reconstruction approach can, therefore, be utilized methodically given different state-space representation forms and different boundary conditions.

\begin{figure}
\begin{subfigmatrix}{2}
 \setlength{\unitlength}{0.012500in}  \vspace{-1mm}
  \subfigure[Solution plot at $N=8$. First eigenmode. Exact solution is in black.]{\includegraphics[width=0.4\textwidth]{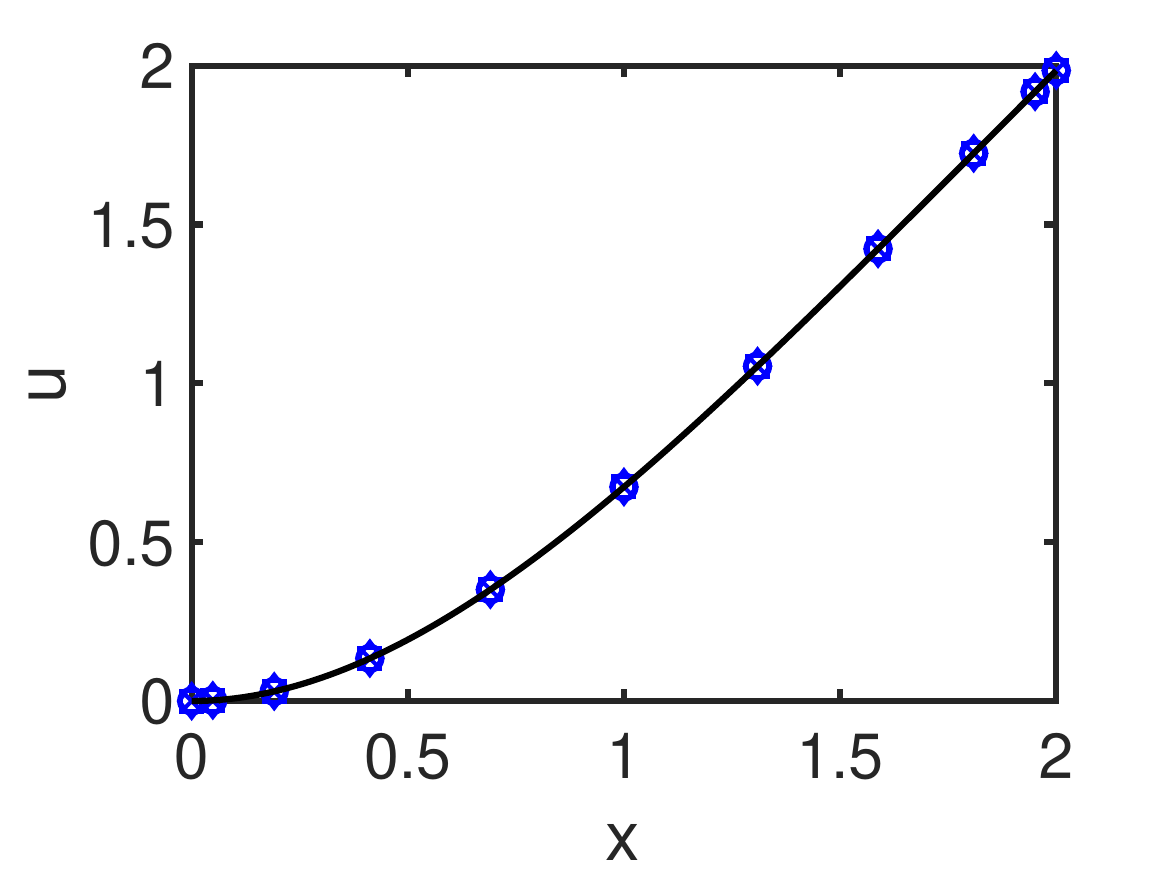}}
  \subfigure[$L_2$ error. First eigenmode.]{\includegraphics[width=0.4\textwidth]{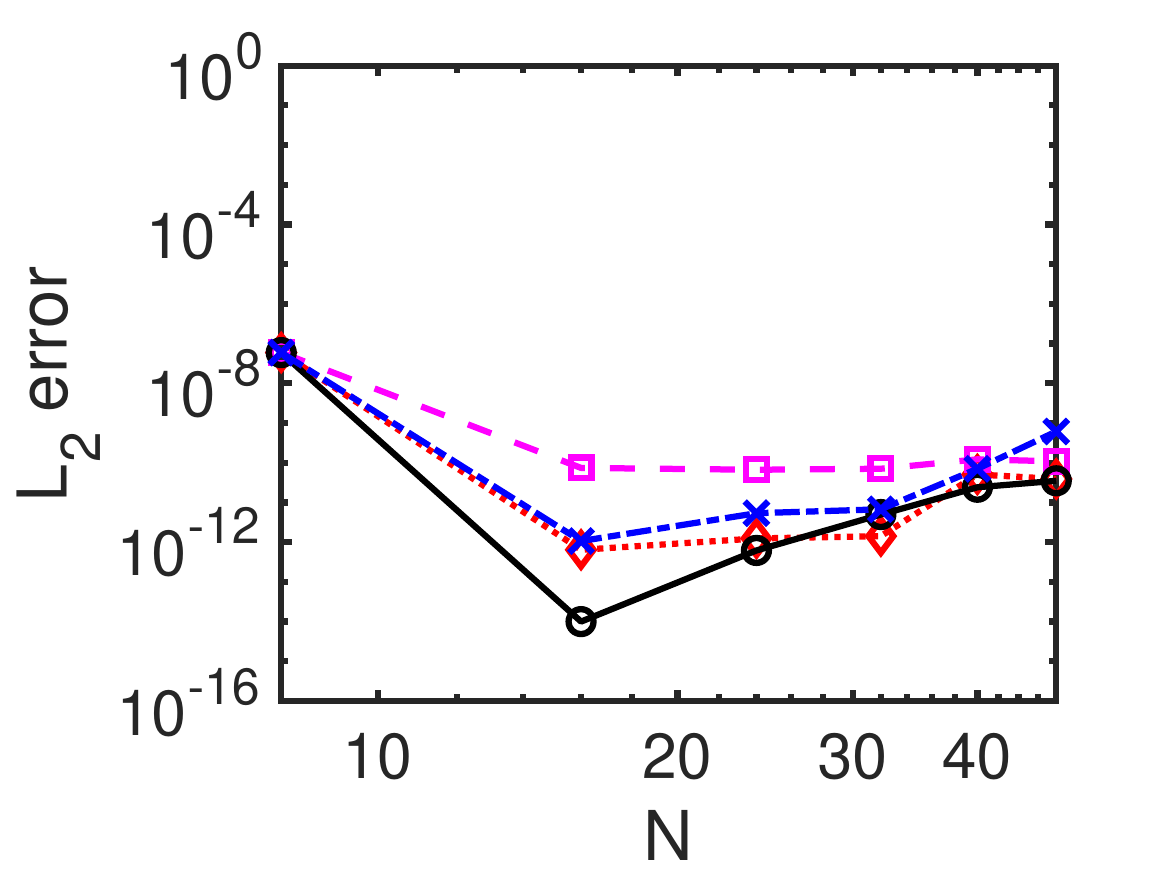}}
  \subfigure[Solution plot at $N=8$. Second eigenmode. Exact solution is in black.]{\includegraphics[width=0.4\textwidth]{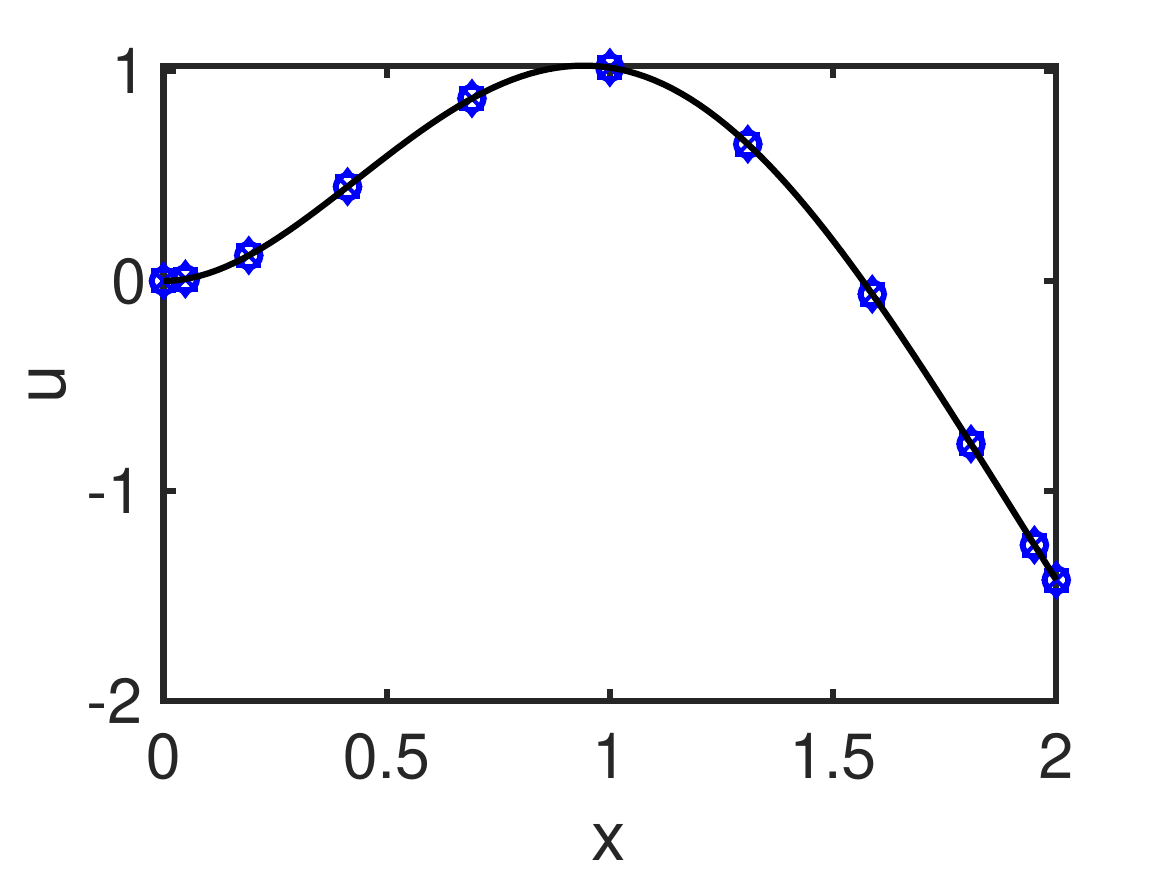}}
  \subfigure[$L_2$ error. Second eigenmode.]{\includegraphics[width=0.4\textwidth]{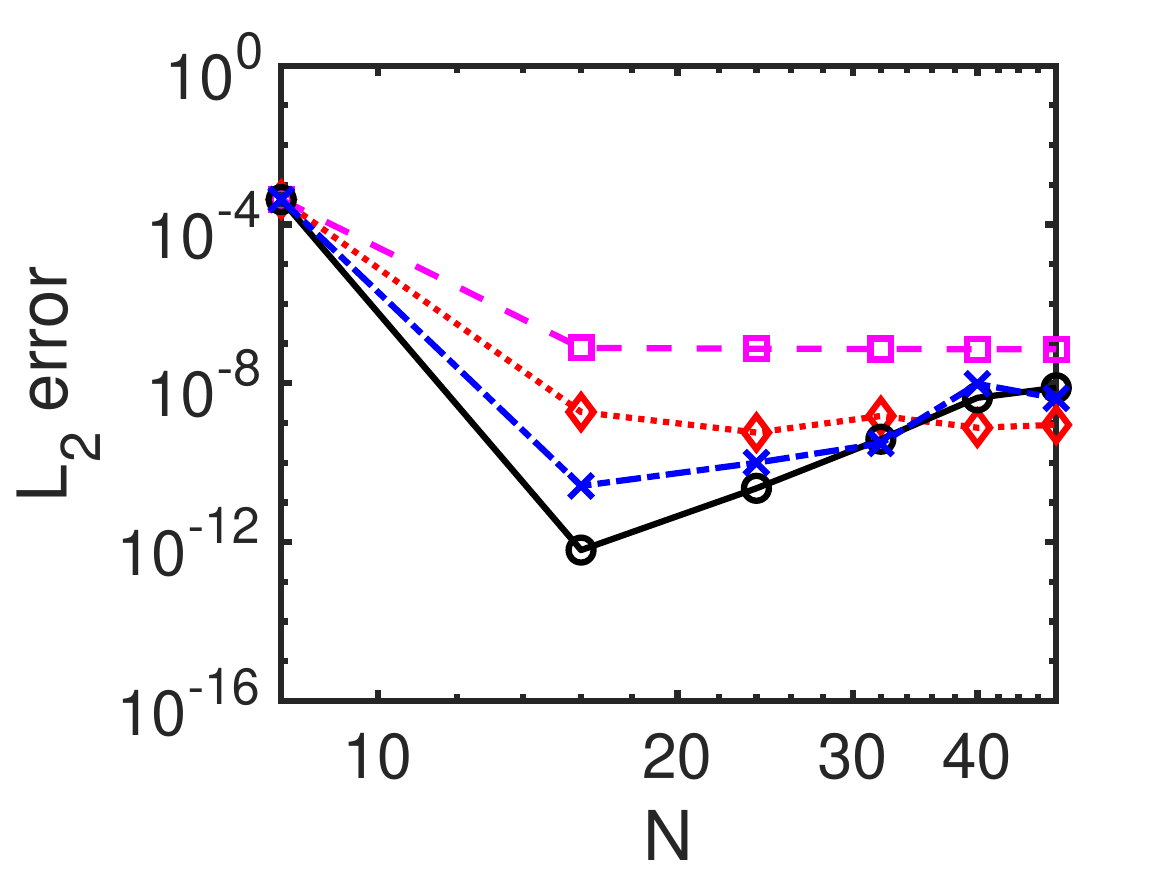}}
  \subfigure[Solution plot at $N=8$ (blue) and $N=16$ (red). Exact solution is in black. Third eigenmode.]{\includegraphics[width=0.4\textwidth]{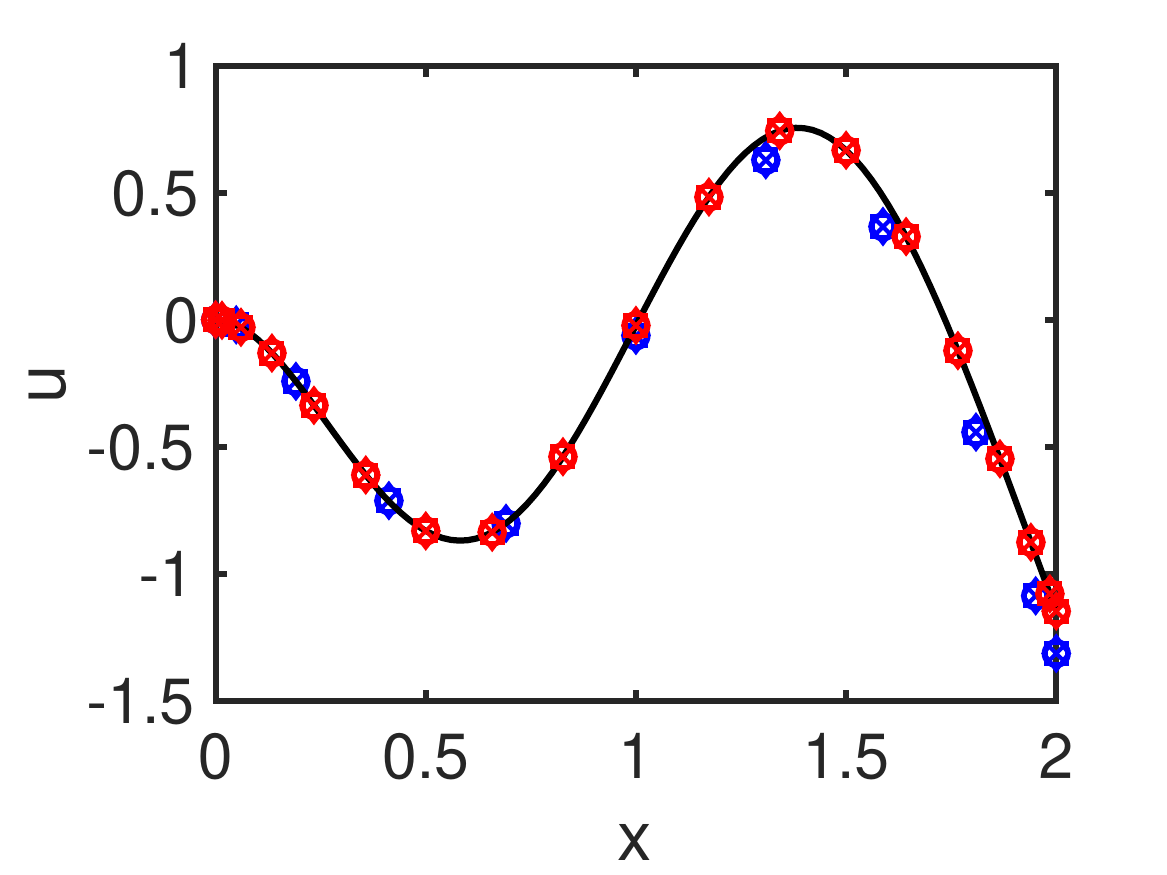}}
  \subfigure[$L_2$ error. Third eigenmode.]{\includegraphics[width=0.4\textwidth]{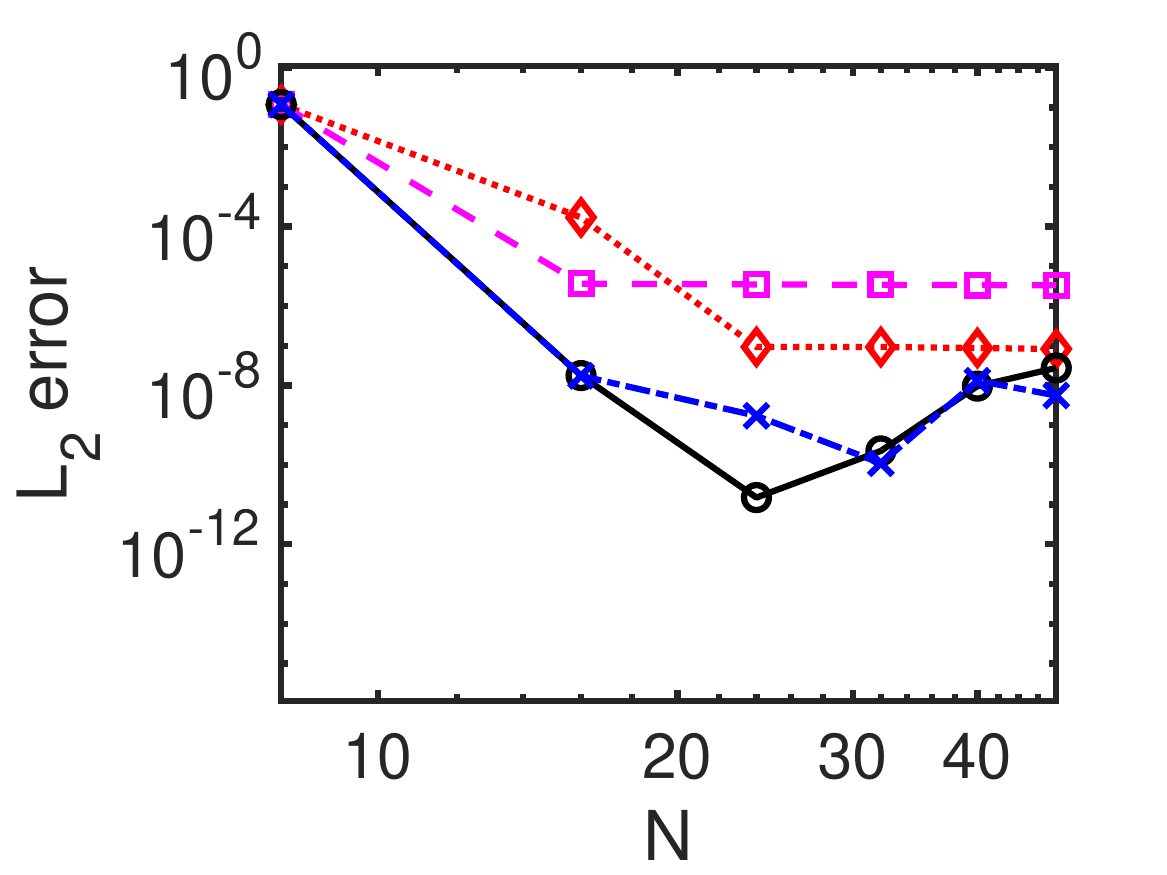}}
  \subfigure[Solution plot at $N=8$ (blue) and $N=16$ (red). Exact solution is in black. Fourth eigenmode.]{\includegraphics[width=0.4\textwidth]{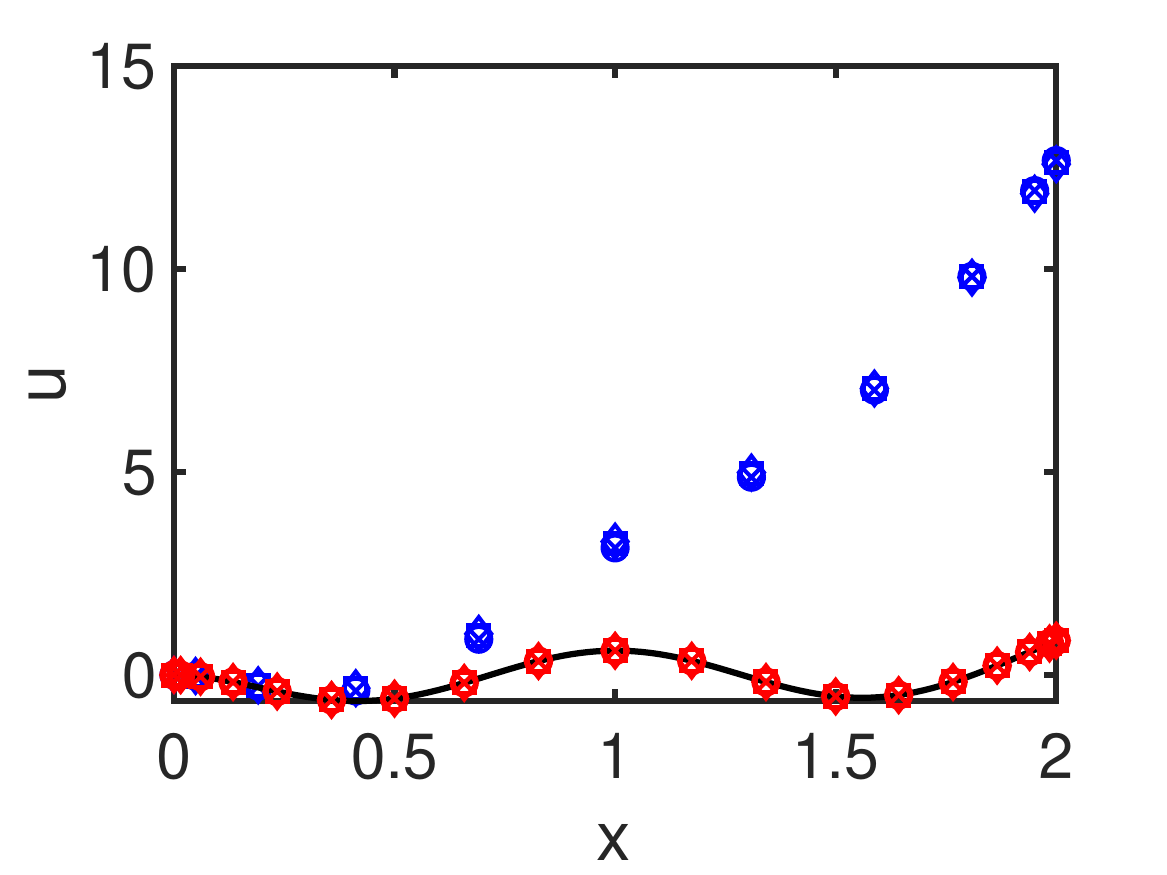}}
  \subfigure[$L_2$ error. Fourth eigenmode.]{\includegraphics[width=0.4\textwidth]{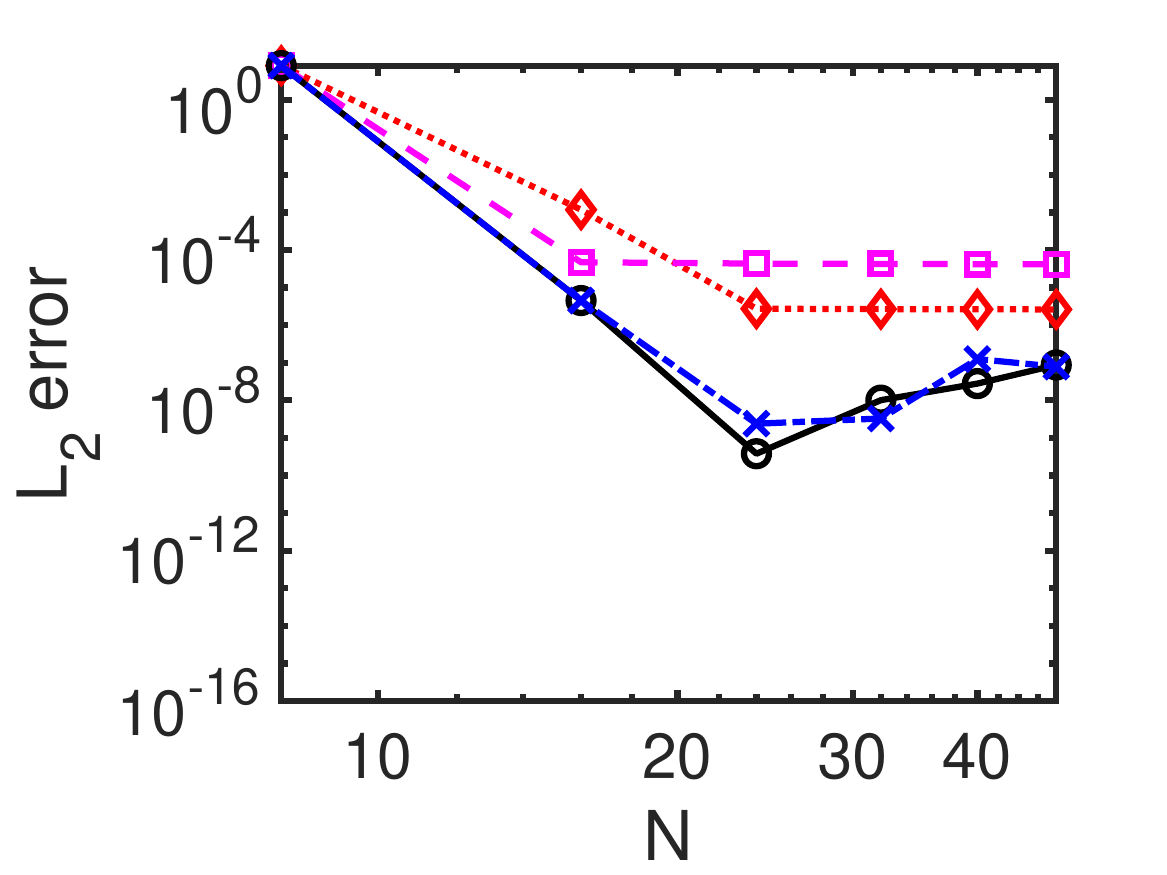}}
      \end{subfigmatrix}
     \caption{Solution and convergence plots for Example 5: Euler-Bernoulli beam equation with $c=2$ at a time $t=0.1$. Rows from one to four correspond to the first through the fourth eigenmodes, respectively. Left panel shows the solution, while the right panel illustrates convergence plots. Lines and symbols are the same as in \cref{fig:heatcon}.\label{fig:beam} }
      \end{figure}

In the following, we choose $L=2$ and keep our solution domain at $x^{(c)}\in[-1,1]$ while recovering the original solution in $x\in[0,L=2]$ by the transformation $x=x^{(c)}+1$. With this, the 3-PI operators become $G^{(c)}_0=0$, $G^{(c)}_1=\bmat{x^{(c)}-s^{(c)}&0\\0&x^{(c)}-s^{(c)}}$, $G^{(c)}_2=\bmat{0&0\\0&-x^{(c)}+s^{(c)}}$, $H^{(c)}_0=\bmat{0&-c\\1&0}$, $H^{(c)}_1=H^{(c)}_2=0$, \\ and $K^{(c)}(x^{(c)})B_T^{-1}= \bmat{1&x^{(c)}-1&0&0\\0&0 & 1& x^{(c)}+1}$. The solution and convergence plots for the first four eigenmodes of a cantilever beam are shown in \cref{fig:beam} at $t=0.1$ obtained with $c=2$, $\Delta\,t=10^{-3}$. To compute these solutions, we set the initial conditions corresponding to an eigenmode shape (\ref{eq:beamexact}) with the amplitude $A_n=1$ for each eigenmode, which is an exact solution at $t=0$. It can be seen that the first and second eigenmodes are well captured with $N=8$. The third eigenmode has a slight deviation near the free boundary at $N=8$, but a correct shape with $N=16$, while the fourth egienmode shows a vastly incorrect deflection with $N=8$, while recovering a correct shape with $N=16$. Note that the tolerance in solving a nonlinear eigenvalue problem (\ref{beameig}) must be set to a very low value ($10^{-16}$ was used in the current work) to obtain these convergence plots, otherwise the convergence will be limited by the value of the set tolerance. 

\subsection{Hyperbolic Problems}
 \subsubsection{Example 6: Transport Equation}
 Here, we consider a transport equation of the form
 \begin{equation}
 u_t+c\,u_x=0,
 \end{equation}
 on the domain $x\in[-1,1]$,
 with $A_0(x)=0, A_1(x)=-c, A_2(x)=0$. As opposed to the previous examples, here we have $n_0=n_2=0$, $n_1=1$, leading to a primary state $u_1(x,t)=u(x,t)$, and a fundamental state $u_{f1}(x,t)=u_x(x,t)$. Transport equation admits solutions in the form of right- (for $c>0$), or left- (for $c<0$) propagating waves. We consider a test case of a propagating Gaussian bump given by the exact solution $u(x,t)=\frac{1}{\sigma\sqrt{2\pi}} e^{-\frac{1}{2}(\frac{x-ct-\mu}{\sigma})^2}$, with the corresponding initial condition and a Dirichlet boundary condition. For $c>0$, we specify a Dirichlet boundary condition at the left at $x=-1$.  The matrix $B$ in this case reduces to $B=[1\:\:\:\: 0]$, $K(x)=1$, $K(x)B_T^{-1}=1$, and the 3-PI operators are $G_0=0, G_1=1, G_2=0$ for the $\mathcal{T}$ operator, and $H_0=-c, H_1=H_2=0$ for the  $\mathcal{A}$ operator. Since the transport equation involves $n_1$ state and not $n_2$ state as a fundamental state, the discrete matrix $M$ now looks different, which, for $N=7$, equals to
\begin{equation}\label{eq:tranmat}
M=\bmat{1 & -1/4 & -1/3 & 1/8 & -1/15 & 1/24 & -1/35  \\1 &0 &-1/2 &0 & 0 & 0 & 0\\ 0 & 1/4 & 0 & -1/4 & 0 & 0 & 0 \\  0 & 0 & 1/6& 0& -1/6 & 0 & 0  \\ 0 & 0 & 0 & 1/8 & 0& -1/8 & 0  \\0 & 0 & 0 & 0 & 1/10 &0 &-1/10 \\ 0 & 0 & 0 & 0 & 0 & 1/12 & 0 }\vspace{-1mm}.
\end{equation} 
In accordance with~\cref{lem:chebmat}, the matrix here is tridiagonal (with the exception of the first row), as opposed to pentadiagonal in parabolic problems with $n_2$ states.

\begin{figure}
\begin{subfigmatrix}{2}
 \setlength{\unitlength}{0.012500in}  \vspace{-1mm}
  \subfigure[Solution plot at $N=8$ (blue) and $N=16$ (red). Exact solution is in black.]{\includegraphics[width=0.49\textwidth]{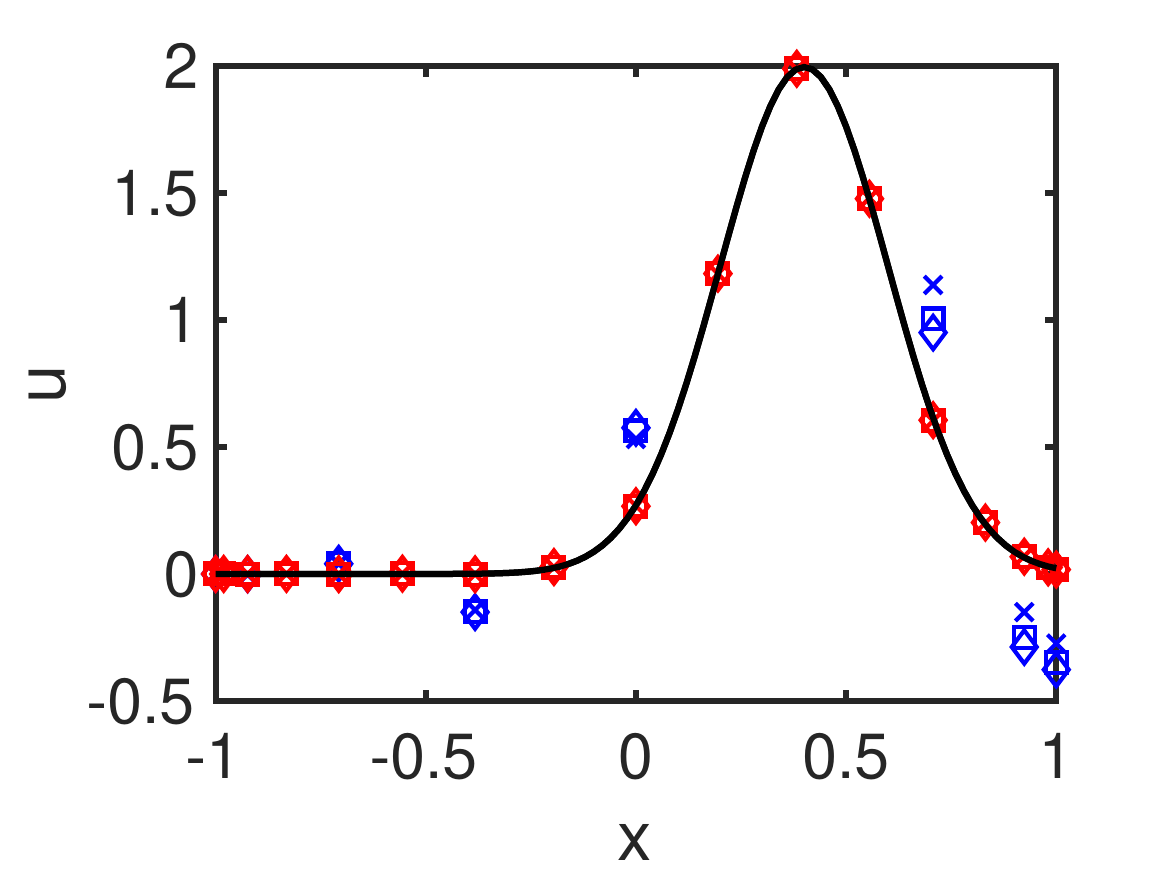}}
  \subfigure[$L_2$ error versus the polynomial order $N$.]{\includegraphics[width=0.49\textwidth]{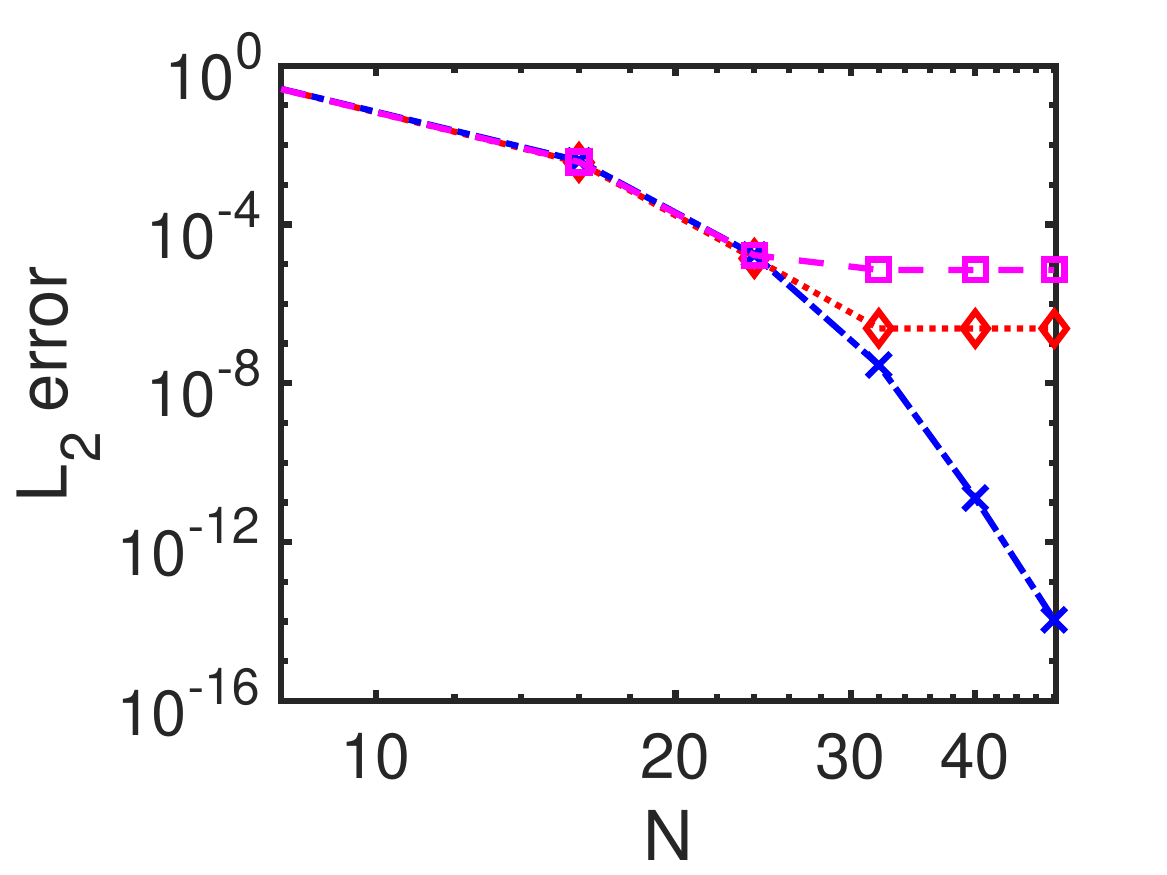}}
  \end{subfigmatrix}
  \caption{Solution and convergence plots for Example 6: transport equation for a propagating Gaussian bump with $c=4$, $\sigma=0.2$, $\mu=0$ at a time $t=0.1$. Blue dash-dotted line with crosses, Gauss integration of Eq. (\ref{timegenorig}) with $N_g=100$ and $N_{int}=1$; red dotted line with diamonds, BDF4  with $\Delta\, t=10^{-3}$; magenta dashed line with squares, BDF3  with $\Delta\, t=10^{-3}$. \label{fig:trangauss} }
  \end{figure}

Choosing $\sigma=0.2$, $\mu=0$ and $c=4$, the solution and the convergence plots are presented  in \cref{fig:trangauss} at a time $t=0.1$. As with the Euler-Bernoulli beam example, it is seen that the Gaussian bump is not well resolved with $N=8$ points, while a correct solution profile is recovered starting at $N=16$. 

\begin{figure}
\begin{subfigmatrix}{2}
 \setlength{\unitlength}{0.012500in}  \vspace{-1mm}
  \subfigure[Solution plot at $N=8$. Solid line, exact solution; symbols, numerical solution.]{\includegraphics[width=0.49\textwidth]{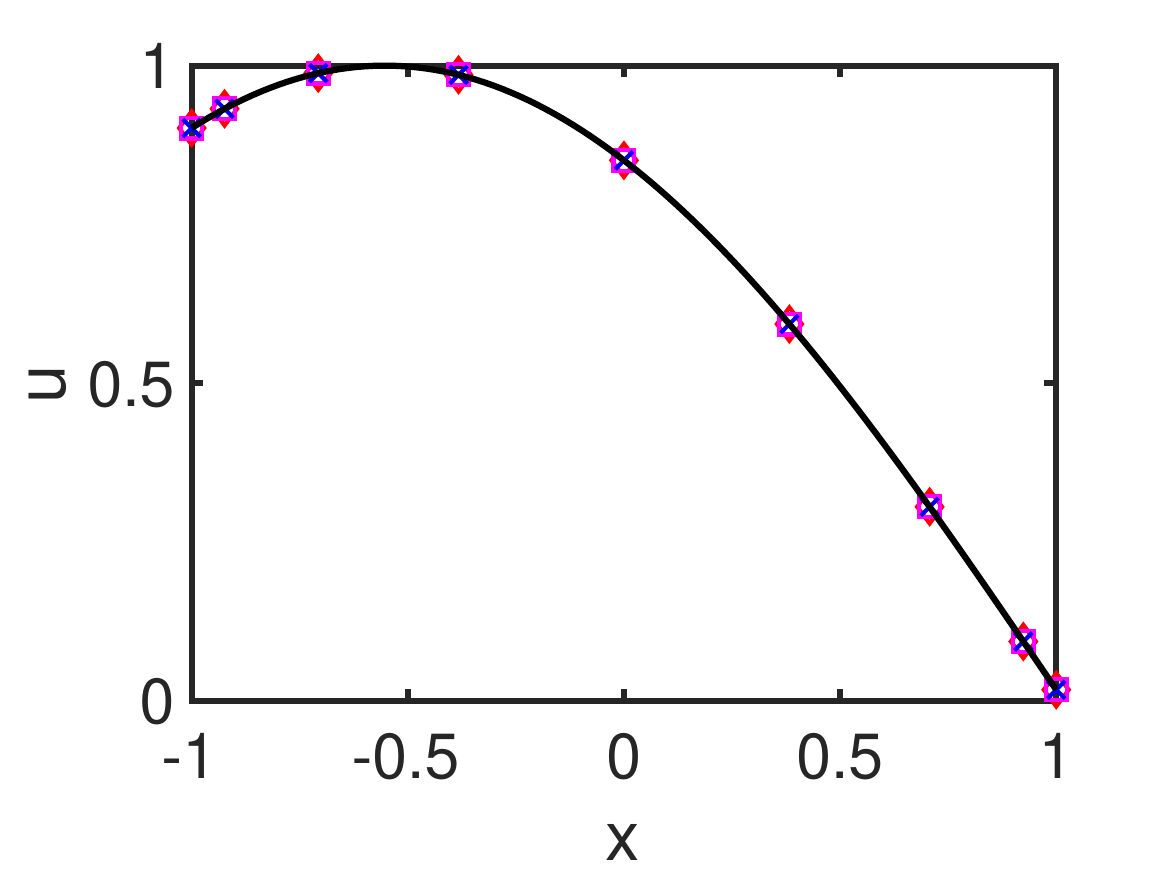}}
  \subfigure[$L_2$ error versus the polynomial order $N$.]{\includegraphics[width=0.49\textwidth]{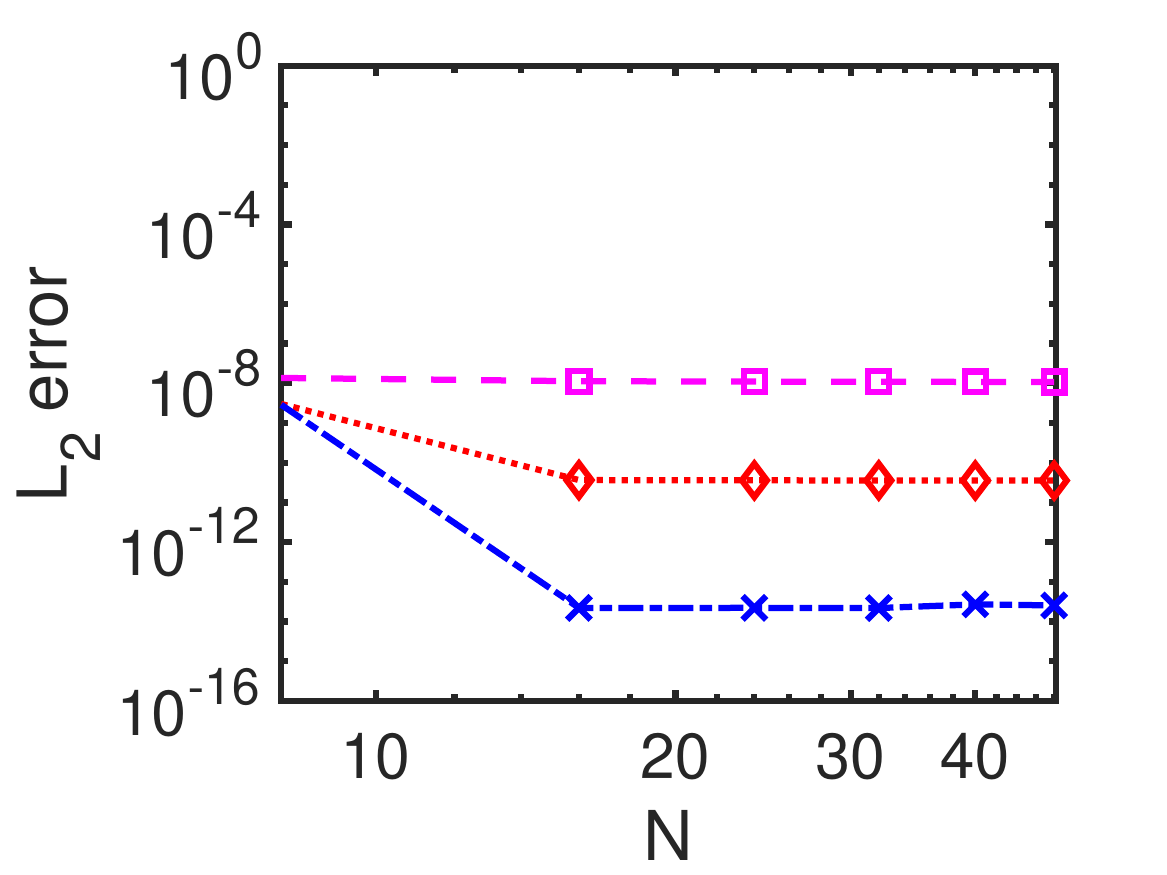}}
      \end{subfigmatrix}
     \caption{Solution and convergence plots for Example 6: transport equation for a traveling sine wave with $c=4$ at a time $t=100$. Blue dash-dotted line with crosses, Gauss integration of Eq. (\ref{timegenorig}) with $N_g=100$ and $N_{int}=100$ uniform time intervals; red dotted line with diamonds, BDF4  with $\Delta\, t=10^{-3}$; magenta dashed line with squares, BDF3  with $\Delta\, t=10^{-3}$. \label{fig:tran} }
      \end{figure}

We test long term integration and conservation properties of the methodology on the example of a traveling sine wave in the form of $u(x,t)=\sin(x-c\,t)$, where initial conditions $u(x,0)=\sin(x)$ and boundary conditions  $u(-1,t)=\sin(-1-c\,t)$ are specified. The results of a long-time integration at $t=100$ and $c=4$ are presented in \cref{fig:tran}. It is seen that the traveling sine wave is well recovered with $N=8$ points, and solution is perfectly conserved even after $t=100$ time units.
 
 \subsubsection{Example 7: Wave Equation}
 \paragraph{Dirichlet-Neumann boundary conditions}
 We now proceed to solving  a wave equation of the form
\begin{equation}\label{eq:wave}
 u_{tt}=c^2\,u_{xx}
 \end{equation}
 on the domain $x\in[-1,1]$ with Dirichlet-Neumann boundary conditions $u(-1,t)=h_1(t), u_x(1,t)=h_2(t)$ and initial conditions 
 \begin{equation}\label{wave:incon}
 u(x,0)=f(x), u_t(x,0)=g(x). 
 \end{equation}
 Exact solution to the wave equation is given by the d'Alembert's formula and depends on the initial conditions for both the function $u(x,0)$ and its time derivative $u_t(x,0)$,
 \begin{equation}\label{dalembert}
 u(x,t)=\frac{1}{2}\left[f(x-ct)+f(x+ct)\right]+\frac{1}{2c}\int_{x-ct}^{x+ct}g(\xi)\,d\,\xi,
 \end{equation}
 where the functions $f(x)$ and $g(x)$ come from the initial conditions (\ref{wave:incon}). Thus, in general, the solution to the wave equation consists of the left- and right- propagating waves. However, in certain situations, depending on the initial conditions, one of the waves can cancel out due to the contribution from the initial conditions on the time derivative, which results in a single left- or right- traveling wave solution.
 
 To reduce a wave equation to its standardized state-space form given by (\ref{eq:primary}), we introduce two states $v_1(x,t)=u_t(x,t)$, $v_2(x,t)=u_x(x,t)$, with the corresponding boundary conditions on the states $v_1(-1,t)=g^\prime_1(t), v_2(1,t)=g_2(t)$, i.e., in terms of the new state vector $\mbf v=\bmat{v_1&v_2}^T$, we have Dirichlet-Dirichlet boundary conditions on both states. With this state vector, the equation (\ref{eq:wave}) now looks
 \begin{equation}
\mbf v_t=\underbrace{\bmat{0&c^{\,2}\\1&0}}_{A_1}\mbf v_{x}.
\end{equation}
The fundamental state is, therefore, $\mbf v_f=\bmat{v_{1x}& v_{2x}}^T$, and we have $n_0=n_2=0$, $n_1=2$, which represents another example of a vector-valued state, as in the Euler-Bernoulli beam equation. The boundary conditions matrix $B$ is assembled as 
 \begin{equation}
 B=\bmat{1&0&0&0\\0&0&0&1},
 \end{equation}
 which leads to the operators $G_0=0$, $G_1=\bmat{1&0\\0&1}$, $G_2=\bmat{0&0\\0&-1}$, $H_0=\bmat{0&c^{\,2}\\1&0}$, $H_1=H_2=0$, and $K(x)B_T^{-1}=\bmat{1&0\\0 &1}$. Since the boundary conditions on the states are not coupled, the matrix $M$ represents a block-diagonal matrix, with each of the two $N\times N$ blocks having the entries identical to the matrix (\ref{eq:tranmat}) of the transport equation, apart from the first row of the second block, where some entries change sign due to boundary conditions. An example of the matrix $M$ for $N=4$ is illustrated below.
  
   \begin{equation}\label{eq:matwave}
    M=\bmat{        1    &       -1/4        &  -1/3        &   1/8       &     0       &      0      &       0        &     0      \\
          1       &      0    &       -1/2    &        0    &         0    &         0     &        0      &       0      \\
          0       &     1/4      &      0       &    -1/4         &   0       &      0        &     0     &        0      \\
          0      &       0     &       1/6      &      0      &       0      &       0      &       0    &         0      \\
          0        &     0      &       0    &       0       &     -1      &     -1/4    &       1/3     &      1/8     \\
           0      &       0      &       0         &    0       &      1      &       0    &       -1/2      &      0     \\
           0       &      0       &      0     &        0    &         0     &       1/4      &      0       &    -1/4     \\
           0        &     0         &    0     &        0    &         0      &       0       &     1/6      &      0 } . 
\end{equation} 
 
 As in the Euler-Bernoulli beam example, to recover the original variable $u(x,t)$ from a state-space variable $u_x(x,t)$, we need to perform an additional transformation $u(x,t)=\mcl T u_x(x,t)+K(x)B_T^{-1}\mbf h(t)$, with $\mcl T=\{0,1,0\}$, $K(x)B_T^{-1}=1, \mbf h(t)=u(-1,t)$ which corresponds to the formula (\ref{tran:ux}).

\begin{figure}
\begin{subfigmatrix}{2}
 \setlength{\unitlength}{0.012500in}  \vspace{-1mm}
  \subfigure[Solution plot at $N=16$ for the splitting case.]{\includegraphics[width=0.49\textwidth]{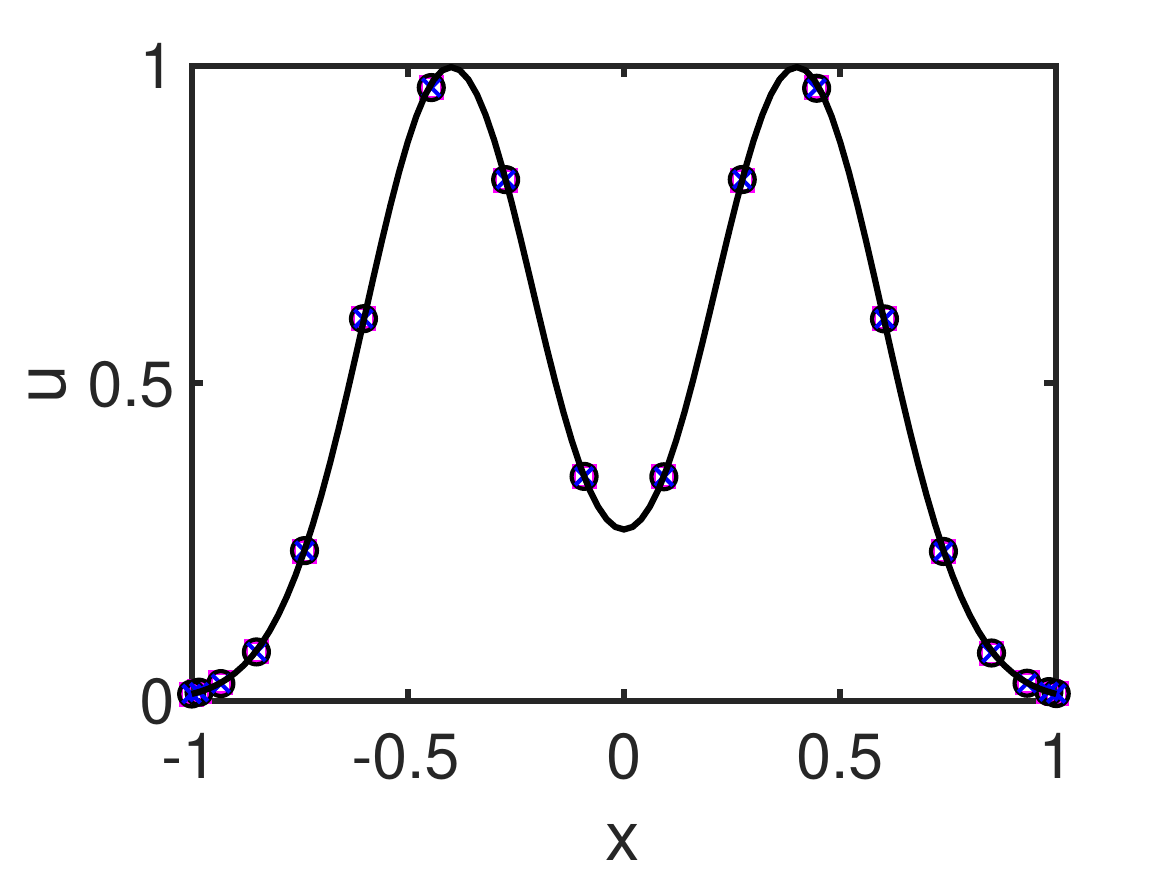}}
  \subfigure[Solution plot at $N=16$ for the right-propagating case.]{\includegraphics[width=0.49\textwidth]{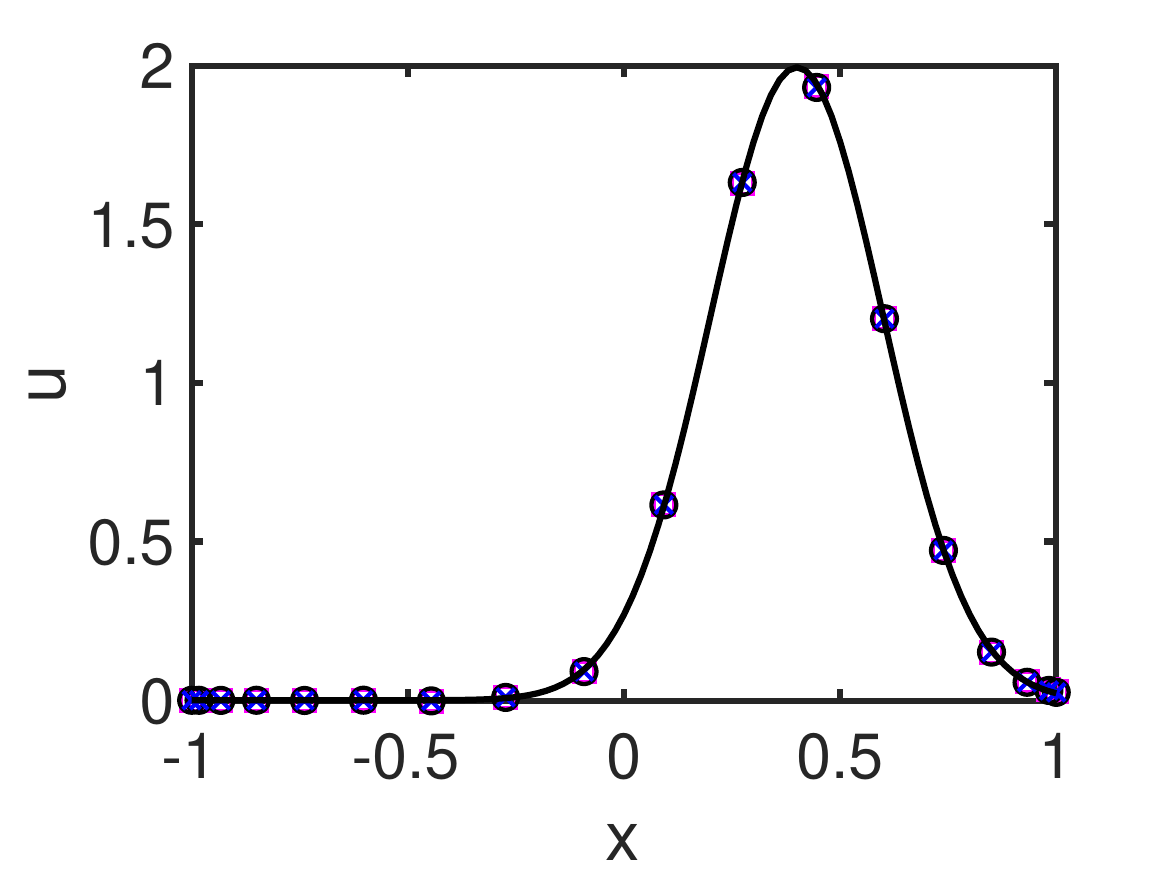}}
  \subfigure[$L_2$ error versus the polynomial order $N$ for the splitting case.]{\includegraphics[width=0.49\textwidth]{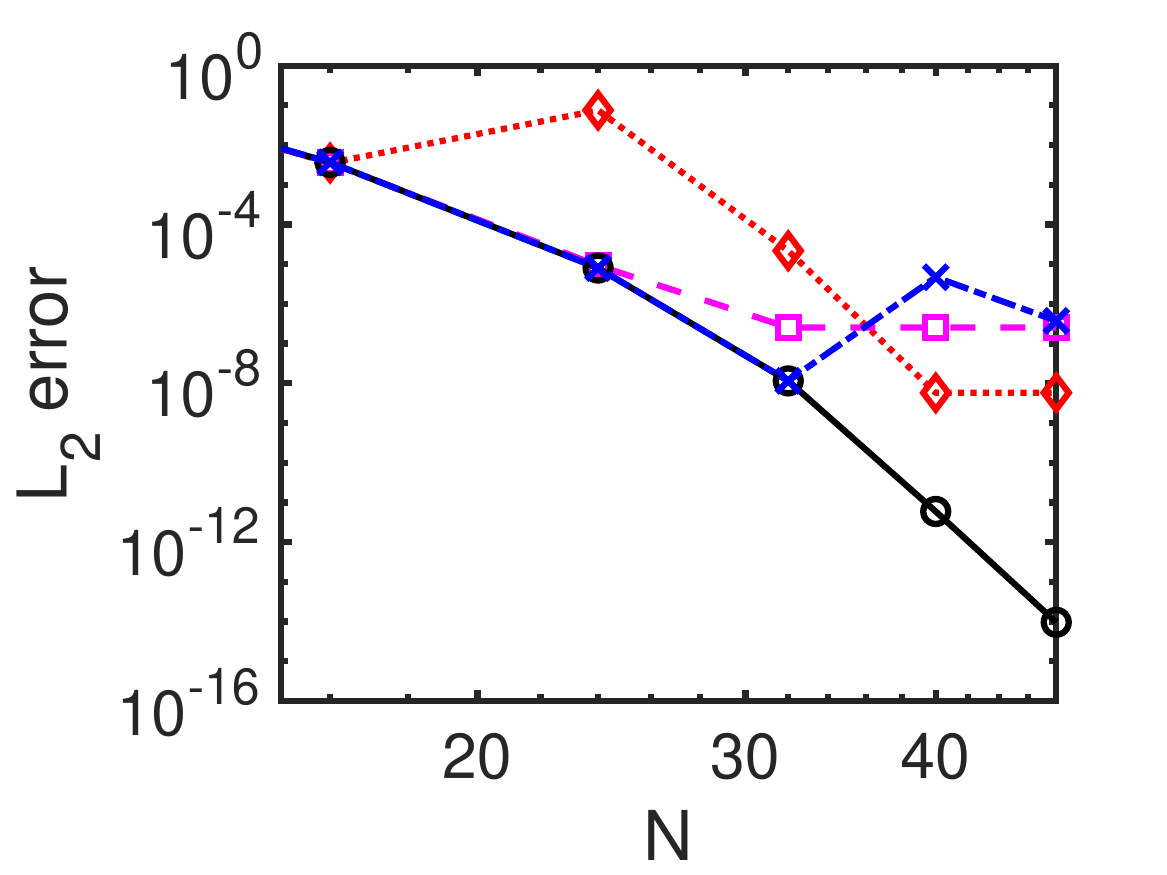}}
  \subfigure[$L_2$ error versus the polynomial order $N$ for the right-propagating case.]{\includegraphics[width=0.49\textwidth]{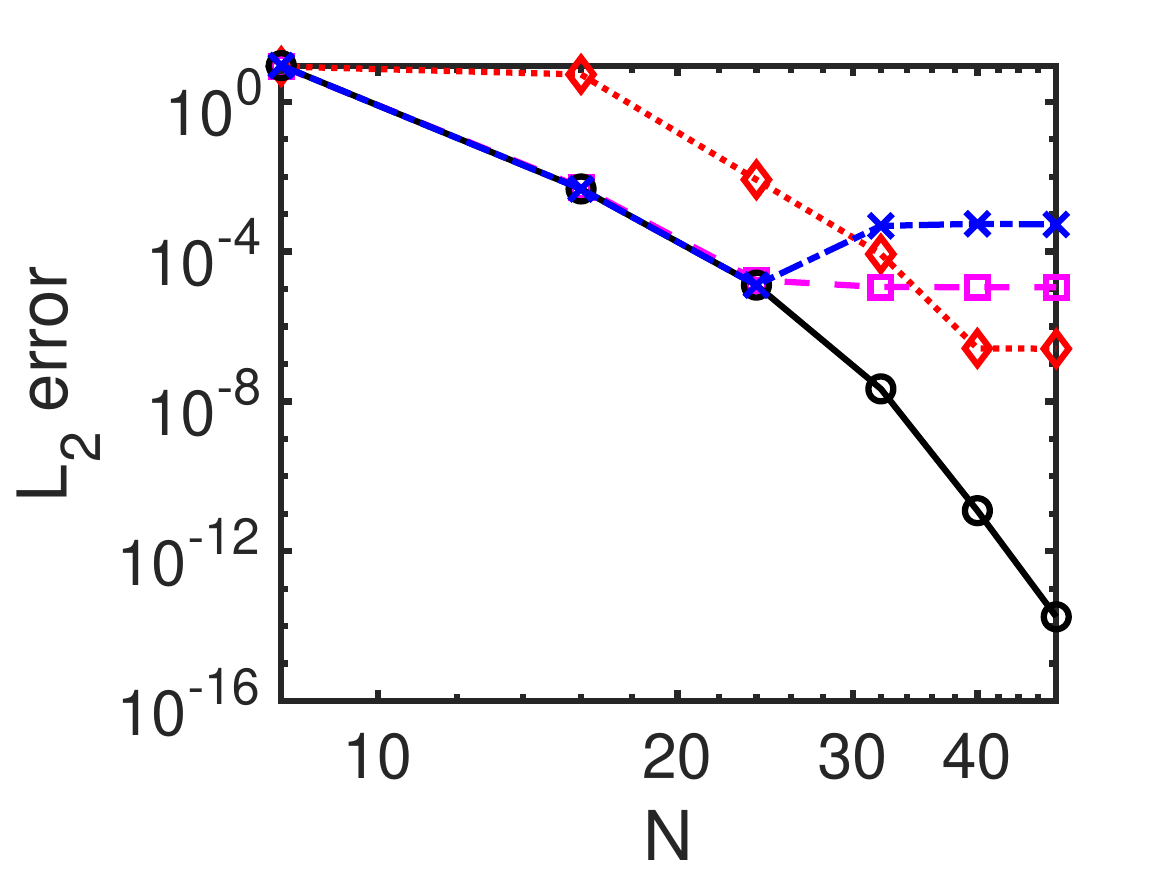}}
      \end{subfigmatrix}
     \caption{Solution and convergence plots for Example 7: wave equation for a Gaussian bump with Dirichlet-Neumann boundary conditions with $c=4$, $\sigma=0.2$, $\mu=0$ at a time $t=0.1$ for the splitting case (left), and the right-propagating case (right).  Black solid line with circles, analytical evaluation of Eq. (\ref{eq:timediag}); blue dash-dotted line with crosses, Gauss integration of Eq. (\ref{timegenorig}) with $N_g=100$ and $N_{int}=1$; red dotted line with diamonds, BDF4 with $\Delta\, t=10^{-3}$; magenta dashed line with squares, BDF3 with $\Delta\, t=10^{-3}$. \label{fig:wave} }
      \end{figure}

 As discussed above, the exact solution to the wave equation depends on the initial conditions on both the functions $u(x,t)$ and $u_t(x,t)$. We first show how, depending on the initial conditions on the derivative $u_t(x,0)$, the same initial shape in a form of a Gaussian bump given by the function $u(x,0)=\frac{1}{\sigma\sqrt{2\pi}} e^{-\frac{1}{2}(\frac{x-\mu}{\sigma})^2}$, can either propagate in one direction, or split in half and give rise to left- and right-propagating waves.
 \paragraph{\underline{Splitting case}}
 According to the d'Alembert's formula (\ref{dalembert}), a splitting case is realized if the initial time derivative $u_t(x,0)=g(x)=0$, and we have the following exact solution
 \begin{equation}
 u(x,t)=\frac{1}{2\,\sigma\sqrt{2\pi}} \left[e^{-\frac{1}{2}(\frac{x-c\,t-\mu}{\sigma})^2}+e^{-\frac{1}{2}(\frac{x+ct-\mu}{\sigma})^2}\right].
 \end{equation}
 \paragraph{\underline{Right-propagating case}}
 In this case, the initial time derivative is specified as
 \begin{equation}
  u_t(x,0)=g(x)=c\left(\frac{x-c\,t-\mu}{\sigma^2}\right) \cdot \frac{1}{\sigma\sqrt{2\pi}} e^{-\frac{1}{2}(\frac{x-ct-\mu}{\sigma})^2}, 
  \end{equation}
  and the exact solution is
 \begin{equation}
 u(x,t)=\frac{1}{\sigma\sqrt{2\pi}} e^{-\frac{1}{2}(\frac{x-ct-\mu}{\sigma})^2}.
 \end{equation} 
 Choosing $\sigma=0.2$, $\mu=0$, and $c=4$, the numerical solution obtained with the PIE-Galerkin framework, and the convergence plots are shown in \cref{fig:wave} at $t=0.1$ obtained with $\Delta\,t=10^{-3}$ for both splitting and right-propagating cases. 
      
       \begin{figure}
\begin{subfigmatrix}{2}
 \setlength{\unitlength}{0.012500in}  \vspace{-1mm}
  \subfigure[Solution plot at $N=8$. Solid line, exact solution; symbols, numerical solution.]{\includegraphics[width=0.49\textwidth]{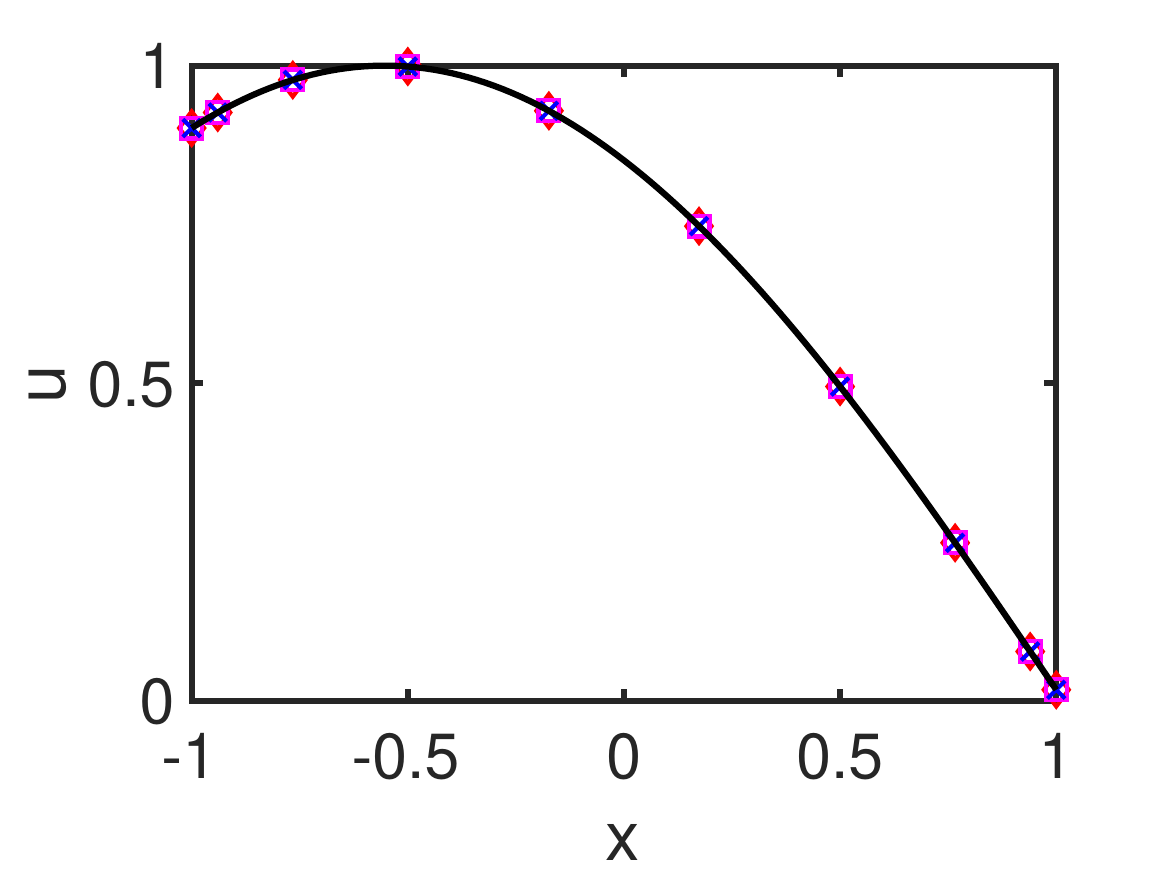}}
  \subfigure[$L_2$ error versus the polynomial order $N$.]{\includegraphics[width=0.49\textwidth]{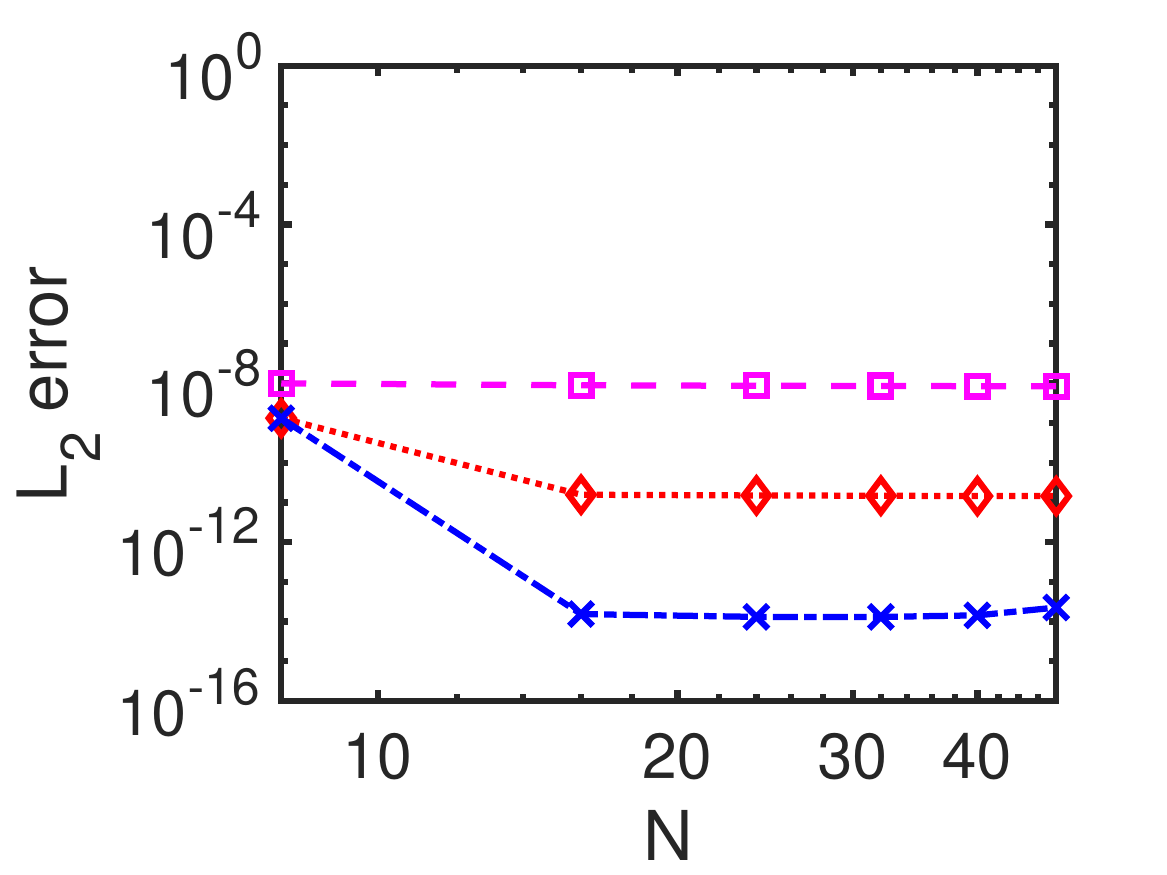}}
      \end{subfigmatrix}
     \caption{Solution and convergence plots for Example 7: wave equation for a traveling sine wave with Dirichlet-Characteristic boundary conditions with $c=4$ at a time $t=100$. Blue dash-dotted line with crosses, Gauss integration of Eq. (\ref{timegenorig}) with $N_g=100$ and $N_{int}=100$ uniform time intervals; red dotted line with diamonds, BDF4  with $\Delta\, t=10^{-3}$; magenta dashed line with squares, BDF3  with $\Delta\, t=10^{-3}$. \label{fig:wavelong} }
      \end{figure}

  \paragraph{Dirichlet-Characteristic boundary conditions}
  
  Since the exact value of the function derivative at the domain outflow is typically not available, we are now considering a characteristic, or a ``non-reflecting'', boundary condition at the right end of the domain given by a characteristics equation $u_t+c\,u_x=0$, while keeping a Dirichlet boundary condition at the left end of the domain. The advantage of the PIE framework is that this boundary condition, which is an optimum choice for an outflow boundary condition in hyperbolic problems, can now be enforced exactly in a strong form. For that, the matrix $B$ is given by  
 \begin{equation}
 B=\bmat{1&0&0&0\\0&0&1&c},
 \end{equation}
 which results in the operators $G_0=0$, $G_1=\bmat{1&0\\0&1}$, $G_2=\bmat{0&0\\-1/c&-1}$, \\ $K(x) B_T^{-1}=\bmat{1 & 0\\-1/c & 1/c}$, and the same value of $\mcl A=\{H_0, H_1, H_2\}$ as in the Dirichlet-Neumann case considered above. The states $u_t$ and $u_x$ are now coupled through the boundary condition, resulting in a matrix $M$ no longer being block-diagonal, but manifesting this coupling across the states in its $N+1^{st}$ row, 
 which is seen, for example, in the matrix $M$ for $N=4$ below,

  \begin{equation}
    M=\bmat{        1    &       -1/4        &  -1/3        &   1/8       &     0       &      0      &       0        &     0      \\
          1       &      0    &       -1/2    &        0    &         0    &         0     &        0      &       0      \\
          0       &     1/4      &      0       &    -1/4         &   0       &      0        &     0     &        0      \\
          0      &       0     &       1/6      &      0      &       0      &       0      &       0    &         0      \\
          -1        &     0      &      1/3     &       0       &     -1      &     -1/4    &       1/3     &      1/8     \\
           0      &       0      &       0         &    0       &      1      &       0    &       -1/2      &      0     \\
           0       &      0       &      0     &        0    &         0     &       1/4      &      0       &    -1/4     \\
           0        &     0         &    0     &        0    &         0      &       0       &     1/6      &      0 } . 
\end{equation} 
The implemented built-in characteristic boundary condition demonstrates a remarkable level of robustness and allows for a long time integration of the wave equation with all the time stepping schemes considered. The results for the traveling sine wave with the same setup as the one described in Example 6  are presented in \cref{fig:wavelong} for the time $t=100$. As with the transport equation, no numerical dissipation or dispersion of the solution is observed at time $t=100$. 
  
\section{Conclusion}\label{sec:conclude}
This paper presents a new theoretical and computational me-thodology to incorporate boundary constraints during a solution of  partial differential equations in a unified and consistent manner. With this methodology, a PDE or a system of PDEs is first transformed into an equivalent partial integral equation (PIE) representation, whose solution lies in a so called fundamental state that does not require boundary conditions, while the latter are analytically embedded into the dynamics of the PIE equation. Not having to enforce boundary conditions on the solution functions brings up several important advantages, such as flexibility in a choice of approximation spaces, enhanced possibilities for analysis and control, and generalizability. As opposed to a weak imposition of the boundary conditions, these advantages do not come at the expense of introducing ad-hoc penalization parameters~\cite{vymazal2019weak, freund1995weakly, juntunen2009nitsche}. In fact, the developed framework is based on firm theoretical grounds, and allows one to retain the fundamental properties of the exact PDE solution in its discrete representation, such as, e.g., the conservation laws~\cite{leveque1992numerical, tadmor2012review, bansal2020structure}.  

A new spectrally-convergent computational technique for a solution of a general class of linear PDEs with variable coefficients and non-periodic boundary conditions is introduced that is based on an expansion of a fundamental solution of the corresponding PIE equation into Chebyshev polynomials of the first kind. With this new methodology, we are able to provide an analytical solution in a form of a function approximation series (spectrally convergent for stable systems) to almost any set of PDEs in the above mentioned class. Furthermore, a general fully-automated programmatic procedure for achieving such solutions for one-dimensional problems is implemented, and is available through an open-source computational solver PIESIM.   Several computational examples that feature parabolic and hyperbolic equation systems are presented, which demonstrate an expected spatial exponential convergence with the polynomial refinement. An approximation solution in a form of a Chebyshev series can be evaluated analytically in time in many practical situations, while a numerical integration in time can also be achieved by employing time discretization techniques. In the current paper, we have evaluated several time integration options, involving an analytical integration whenever possible, and presented their comparison.

A natural further extension of the presented framework, which is currently underway, involves multi-dimensional problems. A possibility of extending to nonlinear cases includes, among other options, treating nonlinearities as non-constant coefficients at each time level, which will be explored in the future work. An extension of the methodology to treat periodic boundary conditions and higher-order PDE systems is also possible, see, e.g.,~\cite{shivakumar2021extension}. Finally, a PDE-PIE reformulation of governing equations presents new avenues for developing a theoretically-consistent treatment of interface conditions and multi-physics coupling laws, which will be explored in the future work as well.

\appendix
\section{Composition Rule for 3-PI operators}\label{sec:composition}
Here, we give the following lemma, which defines the composition rule for the 3-PI operators.
\begin{lemma}\label{thm:composition}
 For any bounded functions $B_0,N_0:[a,b]\rightarrow \R^{n \times n}$, $B_1,B_2,N_1,N_2:[a,b]^2 \rightarrow \R^{n \times n}$, we have\vspace{\eqnspace}
\begin{equation}
\mcl P_{\{R_i\}}=  \mcl P_{\{B_i\}}  \mcl P_{\{N_i\}},\vspace{\eqnspace}
\end{equation}
where\vspace{\eqnspace}
\begin{align}
&R_0(x)=B_0(x)N_0(x), \label{eqn:composition}\\
&R_1(x,s)=B_0(x)N_1(x,s)+B_1(x,s)N_0(s) +\int_{s}^x B_1(x,\xi) N_1(\xi, s)\,d \xi  \notag, \\
&R_2(x,s)=B_0(x)N_2(x,s)+B_2(x,s)N_0(s) \notag \\
&\; +\int_{a}^x B_1(x,\xi) N_2(\xi,s)\,d \xi
+\int_{s}^b B_2(x,\xi) N_1(\xi,s)\,d \xi\notag +\int_{a}^b B_2(x,\xi) N_2(\xi,s)\,d \xi. \notag\\
\notag
\end{align}
\end{lemma}
\begin{proof}
Proof follows from the proof of the Theorem 9 in~\cite{peet_arxiv_PDE} and the relation between the 3-PI operators $\mcl{P}_{\{Q_0,Q_1,Q_2\}}$ in the current paper and $\mcl{P}_{\{Q'_0,Q'_1,Q'_2\}}$ in~\cite{peet_arxiv_PDE} given by
\begin{equation}
Q_0  =Q'_0;\:
Q_1 =Q'_1-Q'_2;\:
Q_2 =Q'_2.
\end{equation}
\end{proof}

\section{Definition of  3-PI Operators in the PIE Representation}\label{sec:appendix} 
This appendix gives a definition of the functions $G_i(x,s), i=0\ldots 5$, appearing in the composition of 3-PI operators in (\ref{eq:3-pi}). 
\begin{align}
  G_0 & =\bmat{I_{n_0}&0&0\\0&0&0\\0&0&0},  G_1(x,s)  =  \bmat{0&0&0\\0&I_{n_1}&0\\0&0&(x-s)I_{n_2}}, \notag \\\notag  \\   G_2(x,s) & =  -K(x)B_T^{-1}BQ(s),\notag\\ \notag \\
  G_3 & =\bmat{0&I_{n_1}&0\\0&0&0},  G_4(s)  =\bmat{0&0&0\\0&0&I_{n_2}}, \notag \\  \notag \\ G_5(s) & = -VB_T^{-1}BQ(s)  \\ \notag \\
  T  & =  \bmat{I_{n_1}&0 &0\\I_{n_1}&0&0 \\0&I_{n_2} &0\\0&I_{n_2}&(b-a)I_{n_2}\\0&0&I_{n_2}\\0&0 &I_{n_2}},  Q(s)=\bmat{0&0&0\\0&I_{n_1} &0 \\ 0&0&0 \\0& 0&(b-s)I_{n_2} \\ 0& 0&0\\ 0&0& I_{n_2}}, \notag
\\ \notag \\  K(x) & =  \bmat{0&0&0\\I_{n_1}&0&0\\0&I_{n_2}&(x-a)I_{n_2}},  V=\bmat{0&0&0\\0&0&I_{n_2}}.\notag \label{eqn:Gdefs}
\end{align}

\section{Proof of \cref{lem:cheb}}\label{sec:app-proof}
\begin{proof}
\begin{enumerate}
\item To prove the first case: if $\mathcal{T}_{mn}$ is such that $m\leq n_0$, according to the structure of $G_0$, $G_1$ and $G_2$, it must have a form $\mathcal{T}_{mn}=\mcl P_{\{\delta_{mn},0,0\}}$, and thus it is easily computed that $\mcl T_{mn} T_{k}(x)=\delta_{mn}T_k(x)$. 
\item To prove the second case, we first need to 
recall some useful recursive relations for Chebyshev polynomials~\cite{canuto1988spectral, moin2001fundamentals}:
\begin{equation}\label{eq:recint}
\int T_k(x)\,dx=
\begin{cases}
T_1(x)+C_0, & k=0 \\
\frac{1}{4}\left[T_0(x)+T_2(x)\right]+C_1, & k=1 \\
\frac{1}{2}\left[\frac{T_{k+1}(x)}{k+1}-\frac{T_{k-1}(x)}{k-1}\right]+C_k, & k\ge 2
\end{cases}
\end{equation}

\begin{equation}\label{eq:recursive}
x\,T_k(x)=
\begin{cases}
T_1(x), & k=0 \\
\frac{1}{2}\left[T_{k-1}(x)+T_{k+1}(x)\right], & k\ge 1
\end{cases}
\end{equation}

Let us now consider $\mcl T_{mn}$ such that $n_0 < m\leq n_0+n_1$. According to the structure of $G_0$, $G_1$ and $G_2$, it has a form of $\mcl T_{mn}=\mcl P_{\{0,\delta_{mn},G_{2mn}\}}$, such that
\begin{equation}\label{g1}
\mcl P_{\{0,\delta_{mn},G_{2mn}\}} T_k(x)=\delta_{mn}\int_{-1}^x T_k(s)\,d\,s+\int_{-1}^1 G_{2mn}(x,s) T_k(s)\,d\,s\\.
\end{equation}
The first integral in the right-hand side can be evaluated according to (\ref{eq:recint}). Let us now consider the second integral. According to the composition of the operator $G_2$, its general entry $G_{2mn}$ would be of the form $G_{2mn}=\beta_{0mn}+\beta_{1mn}\,s+\beta_{2mn}\,x+\beta_{3mn}\,x s$, where $\beta_{jmn}, j=0\ldots 3$, are some real constants. Taking an integral yields
\begin{align}
\int_{-1}^1 G_{2mn}(x,s) T_k(s)\,d\,s=\int_{-1}^1 \left(\beta_{0mn}+\beta_{1mn}\,s+\beta_{2mn}\,x+\beta_{3mn}\,x s\right) T_k(s)\,d\,s \notag \\=
\int_{-1}^1 \left(\beta_{0mn}+\beta_{1mn}\,s\right) T_k(s)\,d\,s+x\int_{-1}^1 \left(\beta_{2mn}+\beta_{3mn}\,s\right) T_k(s)\,d\,s.\label{eq:twoint}
\end{align}
The two integrals in (\ref{eq:twoint}) evaluate to $\gamma_{jkmn} T_0(x)$,  due to the constant limits of integration, where $\gamma_{jkmn}$, $j=0,1$, are some real constants. The muliplication by $x$ in the second integral produces the result $x\cdot\gamma_{1kmn} T_0(x)=\gamma_{1kmn} T_1(x)$. Combining the two integral contributions, (\ref{g1}) can be rewritten as 
\begin{align}
\mcl P_{\{0,\delta_{mn},G_{2mn}\}} T_k(x)=\gamma_{0kmn} T_0(x)+\gamma_{1kmn} T_1(x)\\+\delta_{mn}
\begin{cases}
T_1(x)-T_1(-1),& k=0 \\
\frac{1}{4}\left[T_0(x)+T_2(x)\right]-\frac{1}{4}\left[T_0(-1)+T_2(-1)\right], & k=1 \\
\frac{1}{2}\left[\frac{T_{k+1}(x)}{k+1}-\frac{T_{k-1}(x)}{k-1}\right]-\frac{1}{2}\left[\frac{T_{k+1}(-1)}{k+1}-\frac{T_{k-1}(-1)}{k-1}\right], & k\ge 2
\end{cases}
\\=
b^{(1)}_{0\,kmn} T_0(x)+b^{(1)}_{1kmn} T_1(x)+\delta_{mn}
\begin{cases}
\frac{1}{2}\left[\frac{T_{k+1}(x)}{k+1}\right], & k=1, 2 \\
\frac{1}{2}\left[\frac{T_{k+1}(x)}{k+1}-\frac{T_{k-1}(x)}{k-1}\right], & k\ge 3,
\end{cases}
\end{align}
since $T_k(-1)=(-1)^k=(-1)^k T_0(x)$, leading to (\ref{eq:n1tran}), (\ref{eq:n1coef}).
\item For the third case, we have that $\mcl T_{mn}$, $m>n_0+n_1$ has the form of $\mcl T_{mn}=\mcl P_{\{0,\delta_{mn}(x-s),G_{2mn}\}}$ and 
\begin{equation}\label{g2}
\mcl P_{\{0,\delta_{mn}(x-s),G_{2mn}\}} T_k(x)=\delta_{mn}\int_{-1}^x (x-s) T_k(s)\,d\,s+\int_{-1}^1 G_{2mn}(x,s) T_k(s)\,d\,s\\.
\end{equation}
The last integral in \cref{g2} is evaluated analogously to the previous case. The first integral yields
\begin{equation}
\int_{-1}^x (x-s) T_k(s)\,d\,s=x\int_{-1}^x T_k(s)\,d\,s-\int_{-1}^x s \,T_k(s)\,d\,s.
\end{equation}
Considering the first contribution, we have
\begin{align}
x\int_{-1}^x T_k(s)\,d\,s & & \notag \\=x\begin{cases}
T_1(x)-T_1(-1),& k=0 \notag \\
\frac{1}{4}\left[(T_0(x)+T_2(x)\right]-\frac{1}{4}\left[T_0(-1)+T_2(-1)\right], & k=1 \\
\frac{1}{2}\left[\frac{T_{k+1}(x)}{k+1}-\frac{T_{k-1}(x)}{k-1}\right]-\frac{1}{2}\left[\frac{T_{k+1}(-1)}{k+1}-\frac{T_{k-1}(-1)}{k-1}\right], & k\ge 2
\end{cases} & &\\=\tilde{\alpha}_{1k} T_1(x)+x\begin{cases}
T_1(x),& k=0 \label{contr1}\\
\frac{1}{4}\left[(T_0(x)+T_2(x)\right], & k=1 \\
\frac{1}{2}\left[\frac{T_{k+1}(x)}{k+1}-\frac{T_{k-1}(x)}{k-1}\right], & k\ge 2,
\end{cases}\\=\tilde{\alpha}_{1k} T_1(x)+\begin{cases}
\frac{1}{2}\left[T_{0}(x)+T_{2}(x)\right],& k=0 \\
\frac{1}{4}\left[T_1(x)+\frac{1}{2}\left[T_{1}(x)+T_{3}(x)\right]\right], & k=1 \\
\frac{1}{2}\left[\frac{\frac{1}{2}\left[T_{k}(x)+T_{k+2}(x)\right]}{k+1}-\frac{\frac{1}{2}\left[T_{k-2}(x)+T_{k}(x)\right](x)}{k-1}\right], & k\ge 2,
\end{cases} \notag\\= \tilde{\tilde{\alpha}}_{0k}T_0(x)+\tilde{\tilde{\alpha}}_{1k} T_1(x)+\begin{cases}
\frac{1}{2}\left[\frac{T_{k+2}(x)}{k+1}\right],& k=0 \\
\frac{1}{4}\left[\frac{T_{k+2}(x)}{k+1}\right], & k=1 \\
\frac{1}{4}\left[\frac{T_{k+2}(x)}{k+1}-\frac{2\,T_k(x)}{k^2-1}\right], & k=2, 3\notag\\
\frac{1}{4}\left[\frac{T_{k+2}(x)}{k+1}-\frac{2\,T_k(x)}{k^2-1}-\frac{T_{k-2}(x)}{k-1}\right], & k\ge 4.\notag
\end{cases}
\end{align}

Considering the second contribution, we have
\begin{align}
-\int_{-1}^x s\,T_k(s)\,d\,s=-\int_{-1}^x  ds \begin{cases}
T_1(s), & k=0 \\
\frac{1}{2}\left[T_{k-1}(s)+T_{k+1}(s)\right], & k\ge 1
\end{cases} \notag \\=\tilde{\beta}_{0k} T_0(x)-
\begin{cases}
\frac{1}{4}\left[T_0(x)+T_2(x)\right],& k=0 \\
\frac{1}{2}T_1(x)+\frac{1}{4}\left[\frac{T_{3}(x)}{3}-T_{1}(x)\right], & k=1 \\
\frac{1}{8}\left[T_0(x)+T_2(x)\right]+\frac{1}{4}\left[\frac{T_{4}(x)}{4}-\frac{T_{2}(x)}{2}\right], & k=2 \\
\frac{1}{4}\left[\frac{T_{k}(x)}{k}-\frac{T_{k-2}(x)}{k-2}\right]+\frac{1}{4}\left[\frac{T_{k+2}(x)}{k+2}-\frac{T_{k}(x)}{k}\right], & k\ge 3
\end{cases} \label{contr2}\\
=\tilde{\tilde{\beta}}_{0k} T_0(x)+\tilde{\tilde{\beta}}_{1k} T_1(x)-
\begin{cases}
\frac{1}{2}\left[\frac{T_{k+2}(x)}{k+2}\right]& k=0 \\
\frac{1}{4}\left[\frac{T_{k+2}(x)}{k+2}\right], & 1\le k\le 3\\
\frac{1}{4}\left[\frac{T_{k+2}(x)}{k+2}-\frac{T_{k-2}(x)}{k-2}\right], & k\ge 4\notag
\end{cases}
\end{align}
Combining \cref{eq:twoint},  \cref{g2}, \cref{contr1} and \cref{contr2} 
yields \cref{eq:n2tran} with \cref{eq:n2coef}. 
\end{enumerate}
Dependence of the constants $b^{(i)}_{jkmn}$, $i=1,2,\,j=0,1$, on the boundary conditions comes from the dependence of the operator entries $G_{2mn}$ on the boundary conditions defined by the matrix $B$.
This concludes the proof.
\end{proof}
\section{Proof of \cref{lem:chebmat}}\label{sec:app-proofmat}
\begin{proof}
Evaluating the inner products on both sides of the equation (\ref{eq:test}) with $\pmb \phi_{mn}(x)$ produces the $l^{th}$ out of $N_d$ algebraic equations for the $a_{ik}(t)$ Chebyshev coefficients, where $l=(m-1)n_s+n+1$, which will correspond to the $l^{th}$ row in the associated discrete matrices $M$ and $A$.
Evaluating $\left(\mcl T\: \frac{\partial \hat{\mbf u}_f(x,t)}{\partial \,t}, \pmb \phi_{mn}(x)\right)$, $m=1\ldots n_s, n=0\ldots N-p(m)$, gives, according to (\ref{eq:chebcol}),
\begin{align}\label{eq:proofmat1}
\left(\mcl T\: \frac{\partial \hat{\mbf u}_f(x,t)}{\partial \,t}, \pmb \phi_{mn}(x)\right)=\left(\sum_{i=1}^{ns} \sum_{k=0}^{N-p(i)} \textrm{Col}_{\,i} (\mcl T) \,T_k(x)\,\dot{a}_{ik}(t),\pmb \phi_{mn}(x)\right)\\= \left(\sum_{i=1}^{n_s} \sum_{k=0}^{N-p(i)} \, \mcl T_{mi}\, T_k(x)\, \dot{a}_{ik}(t), T_{n}(x)\right) \notag,
\end{align}
where $\dot{a}_{ik}(t)$ denotes temporal derivative of $a_{ik}(t)$. Equation (\ref{eq:proofmat1}) can be expanded as 
\begin{align}\label{eq:proofmat2}
\left(\mcl T\: \frac{\partial \hat{\mbf u}_f(x,t)}{\partial \,t}, \pmb \phi_{mn}(x)\right)=\left(\sum_{i=1}^{n_0} \sum_{k=0}^{N} \, \mcl T_{mi}\, T_k(x)\, \dot{a}_{ik}(t), T_{n}(x)\right)\\+\left(\sum_{i=n_0+1}^{n_0+n_1} \sum_{k=0}^{N-1} \, \mcl T_{mi}\, T_k(x)\, \dot{a}_{ik}(t), T_{n}(x)\right)+\left(\sum_{i=n_0+n_1+1}^{n_s} \sum_{k=0}^{N-2} \, \mcl T_{mi}\, T_k(x)\, \dot{a}_{ik}(t), T_{n}(x)\right),\notag
\end{align}
Considering the first term in the right-hand side of (\ref{eq:proofmat2}), and according to (\ref{eq:mult}), we can write
\begin{align}
\left(\sum_{i=1}^{n_0} \sum_{k=0}^{N} \, \mcl T_{mi}\, T_k(x)\, \dot{a}_{ik}(t), T_{n}(x)\right)=\left(\sum_{i=1}^{n_0} \sum_{k=0}^{N} \, \delta_{mi}\, T_k(x)\, \dot{a}_{ik}(t), T_{n}(x)\right)\\= \sum_{i=1}^{n_0} \delta_{mi}  \:\Arrowvert T_n(x)\Arrowvert ^2 \dot{a}_{in}(t)= 
\begin{cases}
\Arrowvert T_n(x)\Arrowvert ^2 \dot{a}_{mn}(t),  &m\le n_0\\
0, & \textrm{otherwise}
\end{cases}\notag
\end{align}
Considering the second term in the right-hand side of (\ref{eq:proofmat2}), and according to (\ref{eq:n1tran}), we can write
\begin{align}\label{eq:proofmat3}
\Bigg(\sum_{i=n_0+1}^{n_0+n1} \sum_{k=0}^{N-1} \, \mcl T_{mi}\, T_k(x)\, \dot{a}_{ik}(t), T_{n}(x)\Bigg) =\Bigg(\sum_{i=n_0+1}^{n_0+n1} \sum_{k=0}^{N-1} (b^{(1)}_{0kmi}T_0(x)+b^{(1)}_{1kmi}T_1(x) \notag \\ +\delta_{mi}(c^-_{k-1}T_{k-1}(x)+c^+_{k+1}T_{k+1}(x)) )\dot{a}_{ik}(t), T_{n}(x)\Bigg) \notag\\=\sum_{i=n_0+1}^{n_0+n1} \sum_{k=0}^{N-1}\left(b^{(1)}_{0kmi}\delta_{n0}\Arrowvert T_0(x)\Arrowvert ^2+b^{(1)}_{1kmi}\delta_{n1}\Arrowvert T_1(x)\Arrowvert ^2\right)  \dot{a}_{ik}(t) \\+\sum_{i=n_0+1}^{n_0+n1} \delta_{mi} \Arrowvert T_n(x)\Arrowvert ^2  \left(c^-_{n}\dot{a}_{i(n+1)}(t)+c^+_{n} \dot{a}_{i(n-1)}(t) \right) \notag\\=\sum_{i=n_0+1}^{n_0+n1} \sum_{k=0}^{N-1}\left(b^{(1)}_{0kmi}\delta_{n0}\Arrowvert T_0(x)\Arrowvert ^2+b^{(1)}_{1kmi}\delta_{n1}\Arrowvert T_1(x)\Arrowvert ^2\right)  \dot{a}_{ik}(t) & \notag\\+\begin{cases}
\Arrowvert T_n(x)\Arrowvert ^2\left(c^-_{n}\dot{a}_{m(n+1)}(t)+c^+_{n} \dot{a}_{m(n-1)}(t) \right) & n_0<m\le n_0+n_1\\ 0  & \textrm{otherwise}\end{cases}  \notag
\end{align}
Performing similar manipulations for the third term in the right-hand side of (\ref{eq:proofmat2}), and according to (\ref{eq:n2tran}), one has
\begin{align}\label{eq:proofmat4}
\left(\sum_{i=n_0+n_1+1}^{n_s} \sum_{k=0}^{N-2} \, \mcl T_{mi}\, T_k(x)\, \dot{a}_{ik}(t), T_{n}(x)\right) \notag \\=\sum_{i=n_0+n_1+1}^{n_s} \sum_{k=0}^{N-2}\left(b^{(2)}_{0kmi}\delta_{n0}\Arrowvert T_0(x)\Arrowvert ^2+b^{(2)}_{1kmi}\delta_{n1}\Arrowvert T_1(x)\Arrowvert ^2\right)  \dot{a}_{ik}(t) \\+\begin{cases}
\Arrowvert T_n(x)\Arrowvert ^2\left(d^-_{n}\dot{a}_{m(n+2)}(t)+d_{n} \dot{a}_{mn}(t)+d^+_{n} \dot{a}_{m(n-2)}(t) \right) & m> n_0+n_1\\ 0  & \textrm{otherwise}\end{cases} & \notag
\end{align}
Collecting the corresponding entries multiplying the Chebyshev coefficients $\dot{a}_{mn}(t)\rightarrow \dot{a}_{(m-1)n_s+n+1}(t)$ in the formulas (\ref{eq:proofmat2}), (\ref{eq:proofmat3}), (\ref{eq:proofmat4}) into their respective column positions in the $l^{th}$ row of the matrix $\tilde{M}$, it is easy to see that the structure of the matrix $\tilde{M}=\Lambda M$, where $\Lambda$ is a diagonal matrix consisting of $\Arrowvert T_n(x)\Arrowvert ^2$ in the corresponding diagonal entries $\Lambda_{ll}$,  $l=(m-1)n_s+n+1$. Since the same matrix will be multiplying the matrix $A$ in the right-hand side of equation (\ref{eq:test}), we can multiply both sides of the equation by $\Lambda^{-1}$, which exists due to the entries $\Arrowvert T_n(x)\Arrowvert ^2$ of a diagonal matrix $\Lambda$ being non-zero norms of the Chebyshev polynomials.  The structure of the matrix $M$ described in the proposition of this lemma is now easily deducible from (\ref{eq:proofmat2}), (\ref{eq:proofmat3}) and (\ref{eq:proofmat4}).
\end{proof}


\begin{thebibliography}{10}

\bibitem{encyclopaedia}
{\em Fundamental solution}, in Encyclopaedia of Mathematics, Kluwer, 1994.
\newblock Hazewinkel, M. (Ed.).

\bibitem{greenlibrary}
{\em Green's function library}.
\newblock http://www.greensfunction.unl.edu/home/index.html, 2020.

\bibitem{atkinson1997numerical}
{\sc K.~E. Atkinson}, {\em The numerical solution of boundary integral
  equations}, Clarendon Press, Oxford.
\newblock State of the Art in Numer. Anal., ed. by I. Duff and G. Watson, 1997,
  pp. 223-259.

\bibitem{bansal2020structure}
{\sc H.~Bansal, S.~Weiland, L.~Iapichino, W.~H. Schilders, and N.~van~de Wouw},
  {\em Structure-preserving spatial discretization of a two-fluid model}, in
  IEEE-CDC, 2020, pp.~5062--5067.

\bibitem{bazilevs2007weakimp}
{\sc Y.~Bazilevs and T.~J.~R. Hughes}, {\em Weak imposition of {Dirichlet}
  boundary conditions in fluid mechanics}, Comp. Fluids, 36 (2007), pp.~12--26.

\bibitem{bostrom2017boundary}
{\sc E.~Bostr{\"{o}}m}, {\em Boundary Conditions for Spectral Simulations of
  Atmospheric Boundary Layers}, PhD thesis, KTH Royal Institute of Technology,
  Stockholm, 2017.

\bibitem{bressan1986analysis}
{\sc N.~Bressan and A.~Quarteroni}, {\em Analysis of {Chebyshev} collocation
  methods for parabolic equations}, SIAM Journal on Numerical Analysis, 23
  (1986), pp.~1138--1154.

\bibitem{brezis1998partial}
{\sc H.~Brezis and F.~Browder}, {\em Partial differential equations in the 20th
  century}, Adv. Math., 135 (1998), pp.~76--144.

\bibitem{canuto1986boundary}
{\sc C.~Canuto}, {\em Boundary conditions in {Chebyshev and Legendre} methods},
  SIAM J. Numer. Anal., 23 (1986), pp.~815--831.

\bibitem{canuto1988spectral}
{\sc C.~Canuto, M.~Y. Hussaini, A.~Quarteroni, and T.~A. Zang}, {\em Spectral
  Methods in Fluid Dynamics}, Springer--Verlag, 1988.

\bibitem{canuto1982error}
{\sc C.~Canuto and A.~Quarteroni}, {\em Error estimates for spectral and
  pseudospectral approximations of hyperbolic equations}, SIAM Journal on
  Numerical Analysis, 19 (1982), pp.~629--642.

\bibitem{carvalho2020asymptotic}
{\sc C.~Carvalho, S.~Khatri, and A.~D. Kim}, {\em Asymptotic approximations for
  the close evaluation of double-layer potentials}, SIAM J. Sci. Comp., 42
  (2020), pp.~A504--A533.

\bibitem{das2019h}
{\sc A.~Das, S.~Shivakumar, S.~Weiland, and M.~M. Peet}, {\em {$H_{\infty}$}
  optimal estimation for linear coupled {PDE} systems}, in 58th IEEE Conf.
  Decision and Control (CDC), 2019, pp.~262--267.

\bibitem{deconinck2014method}
{\sc B.~Deconinck, T.~Trogdon, and V.~Vasan}, {\em The method of {Fokas} for
  solving linear partial differential equations}, SIAM Review, 56 (2014),
  pp.~159--186.

\bibitem{deville2002high}
{\sc M.~O. Deville, P.~F. Fischer, and E.~H. Mund}, {\em High-Order Methods for
  Incompressible Fluid Flow}, Cambridge University Press, Cambridge, UK, 2002.

\bibitem{driscoll2010automatic}
{\sc T.~A. Driscoll}, {\em Automatic spectral collocation for integral,
  integro-differential, and integrally reformulated differential equations}, J.
  Comp. Phys., 229 (2010), pp.~5980--5998.

\bibitem{fischer1997overlapping}
{\sc P.~F. Fischer}, {\em An overlapping {Schwarz} method for spectral element
  solution of the incompressible {Navier–Stokes} equations}, J. Comp. Phys.,
  133 (1997), pp.~84--101.

\bibitem{fokas1997unified}
{\sc A.~Fokas}, {\em A unified transform method for solving linear and certain
  nonlinear {PDEs}}, Proc. Roal Soc. Lond. A, 453 (1997), pp.~1411--1443.

\bibitem{fokas1998lax}
{\sc A.~Fokas}, {\em Lax pairs and a new spectral method for linear and
  integrable nonlinear {PDEs}}, Selecta Mathematica, 4 (1998), pp.~31--68.

\bibitem{freund1995weakly}
{\sc J.~Freund and R.~Stenberg}, {\em On weakly imposed boundary conditions for
  second order problems}, in Proceedings of the Ninth Int. Conf. Finite
  Elements in Fluids, 1995, pp.~327--336.

\bibitem{fridman2009lmi}
{\sc E.~Fridman and Y.~Orlov}, {\em An {LMI} approach to {$H^{\infty}$}
  boundary control of semilinear parabolic and hyperbolic systems}, Automatica,
  45 (2009), pp.~2060--2066.

\bibitem{gottlieb1977numerical}
{\sc D.~Gottlieb and S.~A. Orszag}, {\em Numerical Analysis of Spectral
  Methods: Theory and Applications}, SIAM Press, 1977.

\bibitem{greengard1991spectral}
{\sc L.~Greengard}, {\em Spectral integration and two-point boundary value
  problems}, SIAM J. Numer. Anal., 28 (1991), pp.~1071--1080.

\bibitem{grigoryan2010partial}
{\sc V.~Grigoryan}, {\em Partial differential equations}.
\newblock web.math.ucsb.edu/\texttildelow grigoryan/124A.pdf, 2010.

\bibitem{guo2009generalized}
{\sc B.-Y. Guo, J.~Shen, and L.-L. Wang}, {\em Generalized {Jacobi}
  polynomials/functions and their applications}, Applied Numer. Math., 59
  (2009), pp.~1011--1028.

\bibitem{haidvogel1979accurate}
{\sc D.~B. Haidvogel and T.~A. Zang}, {\em The accurate solution of
  {Poisson’s} equation by expansion in {Chebyshev} polynomials}, J. Comput.
  Phys., 30 (1979), pp.~167--180.

\bibitem{hiegemann1997chebyshev}
{\sc M.~Hiegemann}, {\em Chebyshev matrix operator method for the solution of
  integrated forms of linear ordinary differential equations}, Acta Mechanica,
  122 (1997), pp.~231--242.

\bibitem{johnson2010notes}
{\sc S.~G. Johnson}, {\em Notes on {Green’s} functions in inhomogeneous
  media}.
\newblock math.mit.edu/\texttildelow stevenj/18.303/inhomog-notes.pdf, 2010.

\bibitem{jovanovic2016ritz}
{\sc V.~Jovanovic and S.~Koshkin}, {\em The {Ritz} method for boundary problems
  with essential conditions as constraints}, Adv. Math. Physics, 3 (2016),
  pp.~7058017:1--12.

\bibitem{juntunen2009nitsche}
{\sc M.~Juntunen and R.~Stenberg}, {\em Nitsche’s method for general boundary
  conditions}, Math. Comp., 78 (2009), pp.~1353--1374.

\bibitem{karniadakis2005spectral}
{\sc G.~E. Karniadakis and S.~Sherwin}, {\em Spectral/hp Element Methods for
  Computational Fluid Dynamics}, Oxford Science Publications, 2005.

\bibitem{kesici2018numerical}
{\sc E.~Kesici, B.~Pelloni, T.~Pryer, and D.~Smith}, {\em A numerical
  implementation of the unified {Fokas} transform for evolution problems on a
  finite interval}, EJAM, 29 (2018), pp.~543--567.

\bibitem{kythe2012fundamental}
{\sc P.~Kythe}, {\em Fundamental Solutions for Differential Operators and
  Applications}, Springer Sci. Business Media, 2012.

\bibitem{lehotzky2016pseudospectral}
{\sc D.~Lehotzky and T.~Insperger}, {\em A pseudospectral tau approximation for
  time delay systems and its comparison with other weighted-residual-type
  methods}, Int. J. Numer. Meth. Eng., 108 (2016), pp.~588--613.

\bibitem{leveque1992numerical}
{\sc R.~J. LeVeque}, {\em Numerical methods for conservation laws}, vol.~3,
  Springer, 1992.

\bibitem{marin2012highly}
{\sc O.~Marin, K.~Gustavsson, and A.-K. Tornberg}, {\em A highly accurate
  boundary treatment for confined {Stokes} flow}, Comp. Fluids, 66 (2012),
  pp.~2015--230.

\bibitem{moin2001fundamentals}
{\sc P.~Moin}, {\em Fundamentals of Engineering Numerical Analysis}, Cambridge
  University Press, 2001.

\bibitem{nitsche1971variationsprinzip}
{\sc J.~Nitsche}, {\em {{\"U}ber ein Variationsprinzip zur L{\"o}sung von
  Dirichlet-Problemen bei Verwendung von Teilr{\"a}umen, die keinen
  Randbedingungen unterworfen sind}}, in Abhandlungen aus dem mathematischen
  Seminar der Universit{\"a}t Hamburg, vol.~36, 1971, pp.~9--15.

\bibitem{oberkampf1998issues}
{\sc W.~L. Oberkampf and F.~G. Blottner}, {\em Issues in computational fluid
  dynamics code verification and validation}, AIAA journal, 36 (1998),
  pp.~687--695.

\bibitem{peet2018new}
{\sc M.~M. Peet}, {\em A new state-space representation for coupled {PDEs} and
  scalable {Lyapunov} stability analysis in the {SOS} framework}.
\newblock Proc. IEEE Conf. on Decision and Control, 2018.

\bibitem{peet_arxiv_PDE}
{\sc M.~M. Peet}, {\em A partial integral equation {(PIE)} representation of
  coupled linear {PDEs} and scalable stability analysis using {LMIs}}, 2018.
\newblock arxiv.org/abs/1812.06794.

\bibitem{peet2020partial}
{\sc M.~M. Peet}, {\em A partial integral equation representation of coupled
  linear {PDEs} and scalable stability analysis using {LMIs}}, Automatica, 125
  (2021), pp.~109473: 1--14.

\bibitem{roach1982green}
{\sc G.~F. Roach}, {\em Green's Functions, 2nd Edition}, Cambridge University
  Press, Cambridge, Great Britain, 1982.

\bibitem{roache2002code}
{\sc P.~J. Roache}, {\em Code verification by the method of manufactured
  solutions}, J. Fluids Eng., 124 (2002), pp.~4--10.

\bibitem{ruess2013weakly}
{\sc M.~Ruess, D.~Schillinger, Y.~Bazilevs, V.~Varduhn, and E.~Rank}, {\em
  Weakly enforced essential boundary conditions for {NURBS}‐embedded and
  trimmed {NURBS} geometries on the basis of the finite cell method}, Int. J.
  Numer. Methods Eng., 95 (2013), pp.~811--846.

\bibitem{sesma2001approximate}
{\sc F.~J. S{\'{a}}nchez‐Sesma, R.~Madariaga, and K.~Irikura}, {\em An
  approximate elastic two‐dimensional {Green’s} function for a
  constant‐gradient medium}, Geophys. J. Int., 146 (2001), pp.~237--248.

\bibitem{shen1994efficient}
{\sc J.~Shen}, {\em Efficient spectral-{Galerkin} method {I. Direct} solvers
  for the second and fourth order equations using {Legendre} polynomials}, SIAM
  J. Numer. Anal., 15 (1994), pp.~1489--1505.

\bibitem{shen2003new}
{\sc J.~Shen}, {\em A new dual-{Petrov–Galerkin} method for third and higher
  odd-order differential equations: application to the {KDV} equation}, SIAM J.
  Numer. Anal., 41 (2003), pp.~1489--1505.

\bibitem{shivakumar2020pietools}
{\sc S.~Shivakumar, A.~Das, and M.~M. Peet}, {\em {PIETOOLS: A MATLAB} toolbox
  for manipulation and optimization of partial integral operators}.
\newblock In Proceedings of 2020 American Control Conference (ACC), Denver, CO,
  USA, 2020, pp. 2667-2672.

\bibitem{shivakumar2021extension}
{\sc S.~Shivakumar, A.~Das, S.~Weiland, and M.~Peet}, {\em An extension of
  {PIE} representation of coupled linear {ODE-PDE} systems}, SIAM J. Control
  Optimiz., to be submitted,  (2021).

\bibitem{shivakumar2020duality}
{\sc S.~Shivakumar, A.~Das, S.~Weiland, and M.~M. Peet}, {\em Duality and
  {$H_{\infty}$} optimal control of coupled {ODE-PDE} systems}.
\newblock Proc. $59^{th}$ Conference on Decision in Control (CDC), 2020.

\bibitem{siyyam1997accurate}
{\sc H.~I. Siyyam and M.~I. Syam}, {\em An accurate solution of the {Poisson}
  equation by the {Chebyshev-Tau} method}, J. Comp. Appl. Math., 85 (1997),
  pp.~1--10.

\bibitem{smyshlyaev2005backstepping}
{\sc A.~Smyshlyaev and M.~Krstic}, {\em Backstepping observers for a class of
  parabolic {PDEs}}, Systems \& Control Letters, 54 (2005), pp.~613--625.

\bibitem{stakgold1979green}
{\sc I.~Stakgold}, {\em Green's Functions and Boundary Value Problems},
  Wiley-Interscience Publications, New York, USA, 1979.

\bibitem{tadmor1994spectral}
{\sc E.~Tadmor}, {\em Spectral methods for hyperbolic problems}, 1994.
\newblock Lecture Notes delivered at Ecole Des Ondes, Inria-Rocquencort,
  France.

\bibitem{tadmor2012review}
{\sc E.~Tadmor}, {\em A review of numerical methods for nonlinear partial
  differential equations}, Bull. Amer. Math. Soc., 49 (2012), pp.~507--554.

\bibitem{treharne2007initial}
{\sc P.~A. Treharne and A.~Fokas}, {\em Initial-boundary value problems for
  linear {PDEs} with variable coefficients}, Math. Proc. Camb. Phyl. Soc., 143
  (2007), pp.~221--242.

\bibitem{manen2005modeling}
{\sc D.~J. van Manen, J.~O.~A. Robertson, and A.~Curtis}, {\em Modeling of wave
  propagation in inhomogeneous media}, Phys. Rev. Letters, 94 (2005),
  p.~164301.

\bibitem{volterra1965dynamics}
{\sc E.~Volterra and E.~Zachmanoglou}, {\em Dynamics of Vibrations}, Charles E.
  Merrill Books, 1965.

\bibitem{vymazal2019weak}
{\sc M.~Vymazal, D.~Moxey, C.~D. Cantwell, S.~J. Sherwin, and R.~M. Kirby},
  {\em On weak {Dirichlet} boundary conditions for elliptic problems in the
  continuous {Galerkin} method}, Journal of Computational Physics, 394 (2019),
  pp.~732--744.

\bibitem{yu2019jacobi}
{\sc X.~Yu, Z.~Wang, and H.~Li}, {\em Jacobi-{Sobolev} orthogonal polynomials
  and spectral methods for elliptic boundary value problems}, Comm. Applied
  Math. Comp., 1 (2019), p.~283.

\end{thebibliography}
\end{document}